\pdfoutput=1
\documentclass[12pt]{amsart}
\usepackage{etex} 
\usepackage{amscd}
\usepackage{amssymb}
\usepackage[noadjust]{cite}
\usepackage{booktabs}
\usepackage{url}
\usepackage{hyphenat}
\usepackage{epsf}
\usepackage{mathtools}
\usepackage{ifpdf}
\usepackage{mathrsfs}
\usepackage{array}
\usepackage[curve]{xypic}
\usepackage[T1]{fontenc}
\usepackage[utf8]{inputenc}
\usepackage{xr-hyper}
\ifpdf
  \usepackage[pdftex,final]{graphicx}
  \usepackage[pdftex,lmargin=1in,rmargin=1in,tmargin=1in,bmargin=1in]{geometry}
  \usepackage[bookmarks=true, bookmarksopen=true,%
    bookmarksdepth=3,bookmarksopenlevel=2,%
    colorlinks=true,%
    linkcolor=blue,%
    citecolor=blue,%
    filecolor=blue,%
    menucolor=blue,%
    urlcolor=blue]{hyperref}
  \hypersetup{pdftitle={Heegaard Floer homology as morphism spaces}}
  \hypersetup{pdfauthor={Robert Lipshitz, Peter S. Ozsváth, and Dylan
      P. Thurston}}
\else
  \usepackage[dvips,final]{graphicx}
  \usepackage[dvips,lmargin=1in,rmargin=1in,tmargin=1in,bmargin=1in]{geometry}
  \usepackage[draft]{hyperref}
\fi
\usepackage{tikz}
\usetikzlibrary{matrix,arrows}
\tikzset{cdlabel/.style={above,sloped,
    execute at begin node=$\scriptstyle,execute at end node=$}}
\tikzset{algarrow/.style={->, thick}}
\tikzset{blgarrow/.style={->, thick}}
\tikzset{clgarrow/.style={->, thick}}
\tikzset{tensoralgarrow/.style={double, double equal sign distance, -implies}}
\tikzset{tensorblgarrow/.style={double, double equal sign distance, -implies}}
\tikzset{tensorclgarrow/.style={double, double equal sign distance, -implies}}
\tikzset{modarrow/.style={->, dashed}}
\tikzset{othmodarrow/.style={->, thick}}
\tikzset{Amodar/.style={->, dashed}}
\tikzset{Dmodar/.style={->, dashed}}

\newread\testin

\graphicspath{{draws/}{mpdraws/}}
\makeatletter
\def\input@path{{}{draws/}}
\makeatother

\def\mathcenter#1{%
  \vcenter{\hbox{$#1$}}%
}

\DeclareRobustCommand{\widebar}[1]{\overline{#1}{}}

\makeatletter
\newcommand\mi@kern[1]{%
  \settowidth\@tempdima{$\mi@obj^{#1}$}
  \kern-\@tempdima
  #1
  \settowidth\@tempdima{$\mi@obj$}
  \kern\@tempdima
}

\newtoks\mi@toksp
\newtoks\mi@toksb
\DeclareRobustCommand{\manyindices}[5]{
  \def\mi@obj{#5}
  \mi@toksp\expandafter{\mi@kern{#2}}
  \mi@toksb\expandafter{\mi@kern{#1}}
  \@mathmeasure4\textstyle{#5_{#1}^{#2}}
  \@mathmeasure6\textstyle{#5_{#3}^{#4}}
  \dimen0-\wd6 \advance\dimen0\wd4
  \@mathmeasure8\textstyle{\hphantom{{}_{#1}^{#2}}#5^{\the\mi@toksp#4}_{\the\mi@toksb#3}}
  \hbox to \dimen0{}{\kern-\dimen0\box8}
}
\makeatother 

\usepackage{ifpdf}
\ifpdf
  \let\textalt\texorpdfstring
\else
  \newcommand{\textalt}[2]{#1}
\fi

\newcommand{\RR}{\mathbb R}

\newcommand{\ZZ}{\mathbb Z}
\newcommand{\QQ}{\mathbb Q}

\newcommand{\FF}{\mathbb F}

\newcommand{\bD}{\mathbb{D}}

\newcommand{\co}{\colon}

\newcommand{\superset}{\supset}

\newcommand{\bdy}{\partial}

\newcommand{\lbracket}{[}
\newcommand{\rbracket}{]}

\newcommand{\Hyph}{\text{-}}

\newcommand{\spinc}{\mathfrak s}

\DeclareMathOperator{\Hom}{Hom}

\DeclareMathOperator{\Ext}{Ext}
\DeclareMathOperator{\Tor}{Tor}

\DeclareMathOperator{\spin}{spin}
\newcommand{\SpinC}{\spin^c}

\DeclareMathOperator{\ind}{ind}

\DeclareMathOperator{\gr}{gr}

\newcommand{\Barop}{{\mathrm{Bar}}}
\newcommand{\rBarop}{{\mathrm{Bar}_{r}}}
\DeclareMathOperator{\Cobarop}{Cob}

\newcommand\interior{\mathrm{int}}
\newcommand\dr{\mathrm{dr}}

\theoremstyle{plain}
\newtheorem{theorem}{Theorem}
\newtheorem{thmcor}[theorem]{Corollary} 
\numberwithin{equation}{section}
\newtheorem{citethm}[equation]{Theorem}
\newtheorem{proposition}[equation]{Proposition}
\newtheorem{lemma}[equation]{Lemma}
\newtheorem{corollary}[equation]{Corollary}

\newtheorem{definition}[equation]{Definition}
\newtheorem{construction}[equation]{Construction}

\theoremstyle{definition}

\theoremstyle{remark}

\newtheorem{remark}[equation]{Remark}

\newcommand{\HF}{\mathit{HF}}
\newcommand{\HFa}{\widehat {\HF}}

\newcommand{\CFa}{\widehat {\mathit{CF}}}
\newcommand{\HFKa}{\widehat{\mathit{HFK}}}
\newcommand{\x}{\mathbf x}
\newcommand{\y}{\mathbf y}
\newcommand{\z}{\mathbf z}

\newcommand\HH{\mathit{HH}}
\newcommand\HC{\mathit{HC}}
\newcommand\Hochschild\HH

\newcommand{\Ainf}{\mathcal A_\infty}
\newcommand{\Alg}{\mathcal{A}}

\newcommand\Blg{\mathcal{B}}

\newcommand{\alphas}{{\boldsymbol{\alpha}}}
\newcommand{\betas}{{\boldsymbol{\beta}}}

\newcommand{\rhos}{{\boldsymbol{\rho}}}

\newcommand{\bSigma}{\widebar{\Sigma}}

\newcommand{\cM}{\mathcal{M}}
\newcommand{\Mod}{\cM}

\newcommand{\DD}{\textit{DD}}
\newcommand{\DA}{\textit{DA}}

\newcommand{\AD}{\textit{AD}}
\newcommand{\AAm}{\textit{AA}} 
\newcommand{\CFD}{\mathit{CFD}}
\newcommand{\CFDD}{\mathit{CFDD}}
\newcommand{\CFA}{\mathit{CFA}}
\newcommand{\CFAD}{\mathit{CFAD}}
\newcommand{\CFADa}{\widehat{\CFAD}}
\newcommand{\CFDA}{\mathit{CFDA}}
\newcommand{\CFDAa}{\widehat{\CFDA}}
\newcommand{\CFAA}{\mathit{CFAA}}
\newcommand{\CFAAa}{\widehat{\CFAA}}
\newcommand{\CFDa}{\widehat{\CFD}}

\newcommand{\CFK}{\mathit{CFK}}
\newcommand{\CFKa}{\widehat{\CFK}}

\newcommand{\CFDDa}{\widehat{\CFDD}}

\newcommand{\CFAa}{\widehat{\CFA}}

\newcommand{\cZ}{\mathcal{Z}}
\newcommand{\PtdMatchCirc}{\cZ}
\newcommand{\PMC}{\PtdMatchCirc}

\newcommand{\dg}{\textit{dg} }

\newcommand\Id{\mathbb{I}}
\newcommand\Ground{\mathbf k}
\newcommand\Groundl{\mathbf l}

\newcommand\DTP{\mathop{\widetilde\otimes}\nolimits}
\newcommand\DT{\boxtimes}
\newcommand\Gen{\mathfrak{S}}

\newcommand\Tensor{\mathcal T}

\newcommand{\Field}{{\FF_2}}

\newcommand{\Heegaard}{\mathcal{H}}
\newcommand{\HD}{\Heegaard}

\newcommand{\bigGroup}{G'}
\newcommand{\smallGroup}{G}

\newcommand{\DBigGrSet}{S'_D}
\newcommand{\DSmallGrSet}{S_D}

\newcommand{\ModCat}{\mathsf{Mod}}

\DeclareMathOperator{\Mor}{Mor}

\newcommand{\Denis}{\mathsf{AZ}}
\newcommand{\oSmallBar}{\lsup{\Alg}\overline{\textit{bar}}{}^{\Alg}}
\newcommand{\SmallBar}{\lsup{\Alg}\textit{bar}^{\Alg}}
\newcommand{\MirrorDenis}{\overline{\mathsf{AZ}}}

\newcommand{\op}{\mathrm{op}}
\newcommand\PunctF{F^\circ}
\newcommand{\Chord}{\mathrm{Chord}}
\newcommand{\bOne}{\mathbf{1}}
\newcommand{\cG}{\mathcal{G}}
\newcommand{\ol}[1]{\overline{#1}{}}
\newcommand{\arcz}{\mathbf{z}}

\makeatletter
\newcommand\honestalg[3]{\bigl\lbracket
\begin{smallmatrix} #1\@ifempty{#3}{}{&#3} \\ #2 \end{smallmatrix}
\bigr\rbracket}

\makeatother
\newcommand{\lab}[1]{$\scriptstyle #1$}

\newcommand{\sos}[3]{\mathbin{{}_{#1}\mathord#2_{#3}}}
\newcommand{\lsub}[2]{{}_{#1}#2}
\newcommand{\lsup}[2]{{}^{#1}\mskip-.6\thinmuskip#2}
\newcommand{\lsupv}[2]{{}^{#1}#2}

\newcommand{\lsubsupv}[3]{\manyindices{#1}{\mskip.2\thinmuskip#2\mskip-.2\thinmuskip}{}{}{\mathord{#3}}}

\begin{document}
\title{Heegaard Floer homology as morphism spaces}

\author[Lipshitz]{Robert Lipshitz}
\thanks{RL was supported by an NSF Mathematical Sciences Postdoctoral
  Fellowship, NSF Grant number DMS-0905796 and a Sloan
  Research Fellowship.}
\address{Department of Mathematics, Columbia University\\
  New York, NY 10027}
\email{lipshitz@math.columbia.edu}

\author[Ozsv\'ath]{Peter Ozsv\'ath}
\thanks{PSO was supported by NSF grant number DMS-0804121 and a 
  Guggenheim Fellowship.}
\address {Department of Mathematics, Columbia University\\ New York, NY 10027}
\email {petero@math.columbia.edu}

\author[Thurston]{Dylan~P.~Thurston}
\thanks {DPT was supported by a Sloan Research Fellowship.}
\address{Department of Mathematics,
         Barnard College,
         Columbia University\\
         New York, NY 10027}
\email{dthurston@barnard.edu}

\begin{abstract}
  In this paper we prove another pairing theorem for bordered Floer
  homology. Unlike the original pairing theorem, this one is
  stated in terms of homomorphisms, not tensor products. The present
  formulation is closer in spirit to the usual TQFT framework, and
  allows a more direct comparison with Fukaya-categorical
  constructions. The result also leads to various dualities in
  bordered Floer homology.
\end{abstract}

\maketitle

\tableofcontents

\section{Introduction}
In~\cite{LOT1}, we introduced extensions of the Heegaard Floer
homology group $\HFa(Y)$ (with coefficients in $\Field=\ZZ/2\ZZ$) to
$3$-manifolds with boundary. To a surface $F$, together
with a handle decomposition $\PMC$ of $F$ and a little extra data (in
the form of a basepoint), we associated a differential graded algebra
$\Alg(\PMC)$. To a $3$-manifold $Y$ with boundary parameterized by $\PMC$,
we associated a right $\Ainf$-module $\CFAa(Y)$ over $\Alg(\PMC)$ and a
left differential graded module $\CFDa(Y)$ over $\Alg(-\PMC)$ (where
$-$ denotes orientation reversal), each of
which is well-defined up to homotopy equivalence in the corresponding category. These relate to the
closed invariants via the \emph{pairing theorem}, which states that if
$Y_1$ and $Y_2$ are $3$-manifolds with boundaries parameterized by $F$
and $-F$ respectively then
\[
\HFa(Y_1\cup_F Y_2)\cong H_*(\CFAa(Y_1)\DTP\CFDa(Y_2))\eqqcolon\Tor_{\Alg(\PMC)}(\CFAa(Y_1),\CFDa(Y_2));
\]
see~\cite[Theorem~\ref*{LOT:thm:TensorPairing}]{LOT1}.
(We review these constructions a little more thoroughly in
Section~\ref{sec:HF-basics}.)

In this paper, we prove a different pairing theorem, formulated in
terms of the $\Hom$ functor rather than the tensor product functor.  This version
has the advantage that it allows one to work exclusively with
$\CFDa(Y)$ (or, if one prefers, exclusively with $\CFAa(Y)$); this is
of interest since $\CFDa$ is typically easier to compute.
(See~\cite{LOT4}.) 
The present pairing theorem also meshes well with the ``Fukaya-categorical'' formulation of
Lagrangian Floer homology, providing a direct comparison of our
pairing result with Auroux's construction of bordered Floer
homology~\cite{AurouxBordered}. Indeed, our first result is the
following, which also appears as~\cite[``Theorem''
1.5]{AurouxBordered}\footnote{``Scare quotes'' his.}:

\begin{theorem}\label{thm:hom-pair}
  Let $Y_1$ and $Y_2$ be bordered $3$-manifolds with $\bdy Y_1=\bdy
  Y_2=F(\PMC)$. Then
  \begin{equation}
    \begin{split}
      \HFa(-Y_1\cup_\bdy Y_2)&\cong H_*(\Mor_{\Alg(-\PMC)}(\CFDa(Y_1),\CFDa(Y_2)))
      \eqqcolon \Ext_{\Alg(-\PMC)}(\CFDa(Y_1),\CFDa(Y_2))\\
      &\cong H_*(\Mor_{\Alg(\PMC)}(\CFAa(Y_1),\CFAa(Y_2)) )\eqqcolon \Ext_{\Alg(\PMC)}(\CFAa(Y_1),\CFAa(Y_2)).
    \end{split}
  \end{equation}
\end{theorem}

The $\Hom$ pairing theorem (Theorem~\ref{thm:hom-pair}) follows from 
the behavior of the bordered Floer invariants under orientation reversal,
stated as a 
a duality theorem relating $\CFDa(Y)$ and $\CFAa(-Y)$
(Theorem~\ref{thm:or-rev} below). Note that this is a different kind
of duality from the relationship between $\CFDa(Y)$ and $\CFAa(Y)$
from~\cite[Corollary~\ref*{LOT2:cor:ConvertDintoA} and
Theorem~\ref*{LOT2:thm:Duality-precise}]{LOT2}.  In addition to
studying orientation reversal, we also a prove duality theorems in two
other contexts: one corresponding to conjugation of $\SpinC$ structure
(Theorem~\ref{thm:D-is-A}) and another corresponding to reversing the
Morse function on the surface (Theorem~\ref{thm:KoszulDual}).  We also
prove analogues of Theorems~\ref{thm:hom-pair} and~\ref{thm:or-rev}
for bimodules (see Section~\ref{subsec:Analogues}).

As a tool for establishing Theorem~\ref{thm:or-rev}, we use a Heegaard
diagram discovered independently by Auroux~\cite{AurouxBordered} and
Zarev~\cite{Zarev:JoinGlue}
(see Section~\ref{sec:denis}). Studying this diagram gives
several algebraic results, including an algebraic Serre duality
theorem (Theorem~\ref{thm:serre}), and an interpretation of Hochschild
cohomology as a knot Floer homology group
(Corollary~\ref{cor:Hochschild}). It also leads to an in
interpretation of the duality results
from~\cite[Corollary~\ref*{LOT2:cor:ConvertDintoA}]{LOT2} in terms of
Koszul duality (Sections~\ref{sec:dualizing-bimodules}
and~\ref{sec:koszul}).

In spite of its aesthetic appeal, this $\Hom$ version of the pairing
theorem is less economical than the original tensor product pairing theorem: the complex
$\Hom(\CFDa(Y_1),\CFDa(Y_2))$ is typically much larger than
$\CFa(-Y_1\cup_\bdy Y_2)$. (By contrast, for any Heegaard diagram
respecting the decomposition $Y=-Y_1\cup_\bdy Y_2$, the dimensions over
$\Field$ of the complex $\CFa(Y)$ from~\cite{OS04:HolomorphicDisks}
and the complex $\CFAa(-Y_1)\DT\CFDa(Y_2)$ from~\cite{LOT1} are the
same.)

We now explain the duality theorems in more detail.

\subsection{Dualities for bordered Floer modules}

We start with the effect of orientation-reversal on the bordered Floer
invariants. To state it, it is convenient to work at the level of
chain complexes, not homology. In particular, we let $\Mor_A(M,N)$
denote the chain complex of $\Ainf$-morphisms from
$M$ to $N$, whose homology is $\Ext_A(M,N)$
(cf.~Section~\ref{sec:RHom-Ext}).
\begin{theorem}\label{thm:or-rev}
  Let $Y$ be a bordered $3$-manifold, with boundary parameterized by
  $\phi\co F(\PMC)\to \bdy Y$. Let $-Y$ denote $Y$ with its
  orientation reversed and boundary parameterized by the same $\phi$,
  viewed as a map $F(-\PMC)\to \bdy(-Y)$. Then there are homotopy
  equivalences:
  \begin{align}
  \Mor_{\Alg(-\PMC)}(\lsub{\Alg(-\PMC)}\CFDa(Y),\Alg(-\PMC))&\simeq \CFAa(-Y)_{\Alg(-\PMC)} 
  \label{eq:ReverseTypeD}\\
  \Mor_{\Alg(\PMC)}(\CFAa(Y)_{\Alg(\PMC)},\Alg(\PMC))&\simeq \lsub{\Alg(\PMC)}\CFDa(-Y).
  \label{eq:ReverseTypeA}
  \end{align}
\end{theorem}

The above result gives a direct relation between $\CFDa$ and
$\CFAa$, with no orientation reversal, as follows.
The algebra $\Alg(-\PMC)$ is the opposite algebra to $\Alg(\PMC)$. So,
we can regard the left module $\lsub{\Alg(-\PMC)}\CFDa(Y)$ over
$\Alg(-\PMC)$ as a right module $\CFDa(Y)_{\Alg(\PMC)}$. Recall also
that the invariants $\CFDa(Y)$ and $\CFAa(Y)$ decompose according to
(absolute) $\SpinC$-structures on $Y$,
\[
\CFDa(Y)=\bigoplus_{\spinc\in\SpinC(Y)}\CFDa(Y,\spinc)\qquad 
\CFAa(Y)=\bigoplus_{\spinc\in\SpinC(Y)}\CFAa(Y,\spinc).
\]
There is a $\ZZ/2$-action on $\SpinC$-structures, called
\emph{conjugation}, and written $\spinc\mapsto \overline{\spinc}$.
\begin{theorem}\label{thm:D-is-A}
  Let $Y$ be a bordered $3$-manifold, with boundary parameterized by
  $F(\PMC)$, and let $\CFDa(Y,\spinc)_{\Alg(\PMC)}$ denote the
  $\spinc$-summand of $\lsub{\Alg(-\PMC)}\CFDa(Y)$, viewed as a right
  module over $\Alg(\PMC)$. Then
  \begin{equation}
    \CFDa(Y,\spinc)_{\Alg(\PMC)}\cong
    \CFAa(Y,{\overline\spinc})_{\Alg(\PMC)}.\label{eq:D-is-A}
\end{equation}
\end{theorem}
(We are grateful to Denis Auroux for suggesting
Theorem~\ref{thm:D-is-A} to us.)

Heegaard Floer homology for closed three-manifolds satisfies a
conjugation invariance property~\cite[Theorem
2.4]{OS04:HolDiskProperties}.  Theorem~3 gives the
following version of conjugation invariance in the bordered theory:
\begin{thmcor}\label{cor:conj}
  Let $Y$ be a bordered $3$-manifold with boundary
  parameterized by $F(\PMC)$. Then 
  \begin{align}
    \CFAa(Y,\overline{\spinc})_{\Alg(\PMC)}&\simeq (\CFAa(Y,\spinc)\DT_{\Alg(\PMC)} \CFDDa(\Id))_{\Alg(\PMC)}\label{eq:cor-conj1}\\
    \lsub{\Alg(-\PMC)}\CFDa(Y,\overline{\spinc})&\simeq \lsub{\Alg(-\PMC)}(\CFAAa(\Id)\DT_{\Alg(-\PMC)}\CFDa(Y,\spinc)).\label{eq:cor-conj2}
  \end{align}
\end{thmcor}
(Recall that, as in~\cite{LOT1} and~\cite{LOT2}, we use $\DT$ to
denote a particularly convenient model for the $\Ainf$ tensor
product. Also, $\Id$ denotes the identity map of $F(\PMC)$, and
$\CFDDa(\Id)$ and $\CFAAa(\Id)$ the associated type \DD\ and \AAm\ bimodules respectively. See~Theorem~\ref{thm:PreciseDD} for an explicit 
description of $\CFDDa(\Id)$.) 

\subsection{Analogues for bimodules}
\label{subsec:Analogues}
There are several analogues of these theorems for bimodules. When
working with bimodules, many of the theorems require correcting by a
boundary Dehn twist. Recall from \cite{LOT2} that to a bordered
$3$-manifold with two
boundary components $F(\PMC_L)$ and $F(\PMC_R)$, \emph{together with a
  framed arc connecting the boundary components,} one can associate a
bimodule. We can, thus, talk about performing a \emph{boundary Dehn
  twist} on such a strongly bordered $3$-manifold, a Dehn twist along
a loop surrounding the framed arc. (A boundary
Dehn twist $\tau_\bdy$ \emph{decreases} the framing on the arc by $1$.)

Before we can state the bimodule variants, we need one more algebraic
digression. Suppose that $M$ is an $A$-$B$ bimodule which is free (or
projective) as a bimodule, i.e., as a left $A\otimes B^{\op}$ module. Then we can
dualize $M$ over either one or both of the actions. That is, we can
consider both of the bimodules $\Hom_A(M,A)$ and $\Hom_{A\otimes
  B}(M,A\otimes B)$. These are, in general, different bimodules.
Analogous constructions exist in the \dg setting; see
Section~\ref{sec:dualizing}. These two algebraic operations  lead to two different versions of
the $\Hom$ pairing theorem for bimodules. The first of these is:

\begin{theorem}\label{thm:bimod-rev-1}
  Suppose $Y$ is a strongly bordered $3$-manifold with two boundary
  components $F(\PMC_1)$ and $F(\PMC_2)$. Then
  \begin{align}
    \Mor_{\Alg_1'}(\lsub{\Alg_1',\Alg_2'}\CFDDa(Y), \Alg_1'))&\simeq \CFAAa(-Y)_{\Alg_1',\Alg_2'}\label{eq:bimod-rev-1-1}\\
    \Mor_{\Alg_1}(\lsub{\Alg_2'}\CFDAa(Y)_{\Alg_1},\Alg_1)&\simeq
    \lsub{\Alg_1}\CFDAa(-Y)_{\Alg_2'}\label{eq:bimod-rev-1-2}\\
    \Mor_{\Alg_1'}(\lsub{\Alg_1'}\CFDAa(Y)_{\Alg_2},\Alg_1')&\simeq
    \lsub{\Alg_2}\CFDAa(-Y)_{\Alg_1'}\label{eq:bimod-rev-1-3}\\
    \Mor_{\Alg_1}(\CFAAa(Y)_{\Alg_1,\Alg_2},\Alg_1)&\simeq \lsub{\Alg_1,\Alg_2}\CFDDa(-Y).\label{eq:bimod-rev-1-4}
    \end{align}
\end{theorem}
(As explanation for the strange-looking notation, recall that, as
defined in~\cite{LOT2}, the bimodule $\CFDDa(Y)$ has two left actions,
while the bimodule $\CFAAa(Y)$ has two right actions.  For
convenience, we have written $\Alg_i$ for $\Alg(\PMC_i)$ and $\Alg_i'$
for $\Alg(-\PMC_i)$.)

The second type of dualizing leads to the following theorem. 

\begin{theorem}\label{thm:bimod-rev-2}
  Suppose $Y$ is a strongly bordered $3$-manifold with two boundary
  components $F(\PMC_1)$ and $F(\PMC_2)$. Then
  \begin{align}
    \Mor_{\Alg_1'\otimes\Alg_2'}(\lsub{\Alg_1',\Alg_2'}\CFDDa(Y),
    \Alg_1'\otimes_\Field\Alg_2') 
    &\simeq
    \CFAAa(-\tau_\bdy^{-1}(Y))_{\Alg_1',\Alg_2'}\label{eq:bimod-rev-2-1}\\
    \Mor_{\Alg_1\otimes\Alg_2}(\CFAAa(Y)_{\Alg_1,\Alg_2}),\Alg_1\otimes_\Field\Alg_2)
    &\simeq \lsub{\Alg_1,\Alg_2}\CFDDa(-\tau_\bdy^{-1}(Y)).\label{eq:bimod-rev-2-2}
  \end{align}
\end{theorem}
(Here, in the notation we use the quasi-equivalence of categories
between the category of right-right $\Ainf$ $\Alg_1\Hyph\Alg_2$-bimodules and the
category of right $\Ainf$ $\Alg_1\otimes\Alg_2$-modules. This equivalence is
not as obvious as for ordinary bimodules; see, for
instance,~\cite[Section~\ref*{LOT2:sec:A-Bop-vs-bimods}]{LOT2} for
further discussion. Alternately, one could replace
$\Mor_{\Alg_1\otimes\Alg_2}$ with the chain complex of
$\Ainf$-bimodule morphisms, as defined
in~\cite[Section~\ref*{LOT2:sec:bimod-var-types}]{LOT2}, say.)

One can obtain a version of
Theorem~\ref{thm:bimod-rev-2} for $\CFDAa$ by tensoring with
$\CFDDa(\Id)$ and using the fact that this gives an equivalence of
categories (cf.~the proof of Equation~\eqref{eq:bimod-rev-1-2} in
Section~\ref{sec:Consequences}); we leave this to the interested reader.

As in the one boundary component case, Theorems~\ref{thm:bimod-rev-1}
and~\ref{thm:bimod-rev-2} lead to various pairing theorems. For
example:
\begin{thmcor}\label{cor:bimod-hom-pair}
  If $Y_1$ and $Y_2$ are strongly bordered $3$-manifolds, with $\bdy
  Y_1$ parameterized by $F(\PMC_1)$ and $F(\PMC_2)$ and $\bdy Y_2$
  parameterized by $F(\PMC_1)$ and $F(\PMC_3)$, then
  \begin{align}
    \lsub{\Alg_3'}\CFDAa(-Y_1\cup _{F(\PMC_1)}
    Y_2)_{\Alg_2'}&\simeq\Mor_{\Alg_1'}(\lsub{\Alg_2',\Alg_1'}\CFDDa(Y_1),\lsub{\Alg_1',\Alg_3'}\CFDDa(Y_2))\label{bimod-cor-1-1}\\
    \lsub{\Alg_2,\Alg_3'}\CFDDa(\tau_\bdy(-Y_1\cup_\bdy
    Y_2))&\simeq
    \Mor_{\Alg_1\otimes\Alg_1'}(\lsub{\Alg_1,\Alg_1'}\CFDDa(\Id),\label{bimod-cor-1-2}\\
    &\qquad\qquad\lsub{\Alg_1,\Alg_2}\CFDDa(-Y_1)\otimes
    \lsub{\Alg_1',\Alg_3'}\CFDDa(Y_2)).\nonumber
  \end{align}
\end{thmcor}
(Here, as earlier, $\Id$ denotes the identity map of $F(\PMC)$, and
$\CFDDa(\Id)$ the associated type \DD\ module.)

\begin{thmcor}\label{cor:bimod-mod-hom-pair}
  If $Y_1$ is a strongly bordered $3$-manifold with boundary
  parameterized by $F(\PMC_1)$ and $F(\PMC_2)$ and $Y_2$ is a bordered
  $3$-manifold with boundary parameterized by $F(\PMC_2)$ then
  \begin{equation}
    \label{eq:MorDDtoD}
    \begin{split}
      \CFAa(-Y_1\cup_{F(\PMC_2)}Y_2)_{\Alg_1'}&\simeq \Mor_{\Alg_2'}(\lsub{\Alg_1',\Alg_2'}\CFDDa(Y_1),\lsub{\Alg_2'}\CFDa(Y_2))\\
      \lsub{\Alg_1'}\CFDa(-Y_2\cup_{F(\PMC_2)}Y_1)&\simeq
      \Mor_{\Alg_2'}(\lsub{\Alg_2'}\CFDa(Y_2),\lsub{\Alg_1',\Alg_2'}\CFDDa(Y_1)).
    \end{split}
  \end{equation}
  In particular, if $\psi\co F(\PMC_2)\to F(\PMC_1)$ then
  \begin{equation}
    \label{eq:D-to-DD}
    \begin{split}
      \CFAa(\psi(Y_2))_{\Alg_1} &\simeq
        \Mor_{\Alg_2'}(\lsub{\Alg_1,\Alg_2'}\CFDDa(\psi^{-1}),
          \lsub{\Alg_2'}\CFDa(Y_2)\\
      \lsub{\Alg_1'}\CFDa(\psi(Y_2)) &\simeq
        \Mor_{\Alg_2}(\lsub{\Alg_2}\CFDa(-Y_2),
          \lsub{\Alg_1',\Alg_2}\CFDDa(\psi)).
    \end{split}
  \end{equation}
\end{thmcor}

There are also versions of Theorem~\ref{thm:D-is-A} for bimodules.

\begin{theorem} 
  \label{thm:DA-bimod}
  Suppose $Y$ is a
  strongly bordered $3$-manifold with two boundary components
  $F(\PMC_1)$ and $F(\PMC_2)$. Then viewing
  $\lsub{\Alg_1',\Alg_2'}\CFDDa(Y)$ as a right-right
  module $\CFDDa(Y)_{\Alg_1,\Alg_2} $ over
  $\Alg_1$ and $\Alg_2$,
  \begin{equation}
    \CFDDa(Y, \spinc)_{\Alg_1,\Alg_2}\simeq \CFAAa(\tau_\bdy(Y),\overline{\spinc})_{\Alg_1,\Alg_2}.\label{eq:DA-bimod-1}
  \end{equation}

  Similarly, viewing $\CFAAa(Y,\spinc)$ as a left-left module over
  $\Alg_1'$ and $\Alg_2'$, we have:
  \begin{equation}
  \lsub{\Alg_1',\Alg_2'}\CFAAa(Y, \spinc)\simeq
  \lsub{\Alg_1',\Alg_2'}\CFDDa(\tau_\bdy^{-1}(Y),\overline{\spinc}).\label{eq:DA-bimod-2}
 \end{equation}
  Finally, there are two versions of $\CFDAa(Y)$, depending on whether
  we treat $F(\PMC_1)$ or $F(\PMC_2)$ as the type $D$ side. Denote
  these two modules by
  $\lsub{\Alg_1'}\CFDAa(Y,\spinc)_{\Alg_2}$ and
  $\lsub{\Alg_2'}\CFDAa(Y,\overline{\spinc})_{\Alg_1}$,
  respectively. These two
  options are related by conjugation of the $\SpinC$-structure:
  \begin{equation}
  \lsub{\Alg_1'}\CFDAa(Y,\spinc)_{\Alg_2}\simeq
  \lsub{\Alg_2'}\CFDAa(Y,\overline{\spinc})_{\Alg_1}.\label{eq:DA-bimod-3}
\end{equation}
  Here, we mean that the modules are homotopy equivalent if we
  exchange the sidedness of the actions on either one of the two.
\end{theorem}

\subsection{Algebraic consequences}
Finally, these techniques can be used to prove several more algebraic results
about the category of $\Alg(\PMC)$-modules:
\begin{theorem}\label{thm:serre}
  Given right $\Ainf$-modules $M$ and $N$ over $\Alg(\PMC)$,
  \begin{equation}
  \Mor_{\Alg(\PMC)}(N,M\DT \CFDAa(\tau_\bdy^{-1}))\simeq \Mor_{\Alg(\PMC)}(M,N)^*,
\end{equation}
  naturally. Here, $*$ denotes the dual vector space.
\end{theorem}
In other words, Theorem~\ref{thm:serre} says that tensoring with
$\CFDAa(\tau_\bdy^{-1})$ is the Serre functor for the category of
$\Alg(\PMC)$-modules.

In~\cite{LOT2}, we identified the Hochschild homology of $\CFDAa$ with
a certain self-gluing operation. The algebraic results from this paper
allow us to translate this result into one about Hochschild
cohomology. (The non-specialist is reminded that Hochschild cohomology
is typically not dual to Hochschild homology.)
\begin{thmcor}\label{cor:Hochschild}
 Suppose that $Y$ is a strongly bordered $3$-manifold with boundary
 $F(\PMC)\amalg F(-\PMC)$.  Let $\tau_\bdy(Y)$ denote the result of
 decreasing the framing on the arc $\arcz$ in $Y$  by one and
 $\tau_\bdy(Y)^\circ$ the manifold obtained by gluing the two
 boundary components of $\tau_\bdy(Y)$ together and performing surgery
 on the framed knot $K$ coming from the arc $\arcz$.
 Let $K'$ be the knot in $\tau_\bdy(Y)^\circ$ coming from $K$.
 Then the Hochschild cohomology
 $\HH^*(\CFDAa(Y))$ is isomorphic to $\HFKa(\tau_\bdy(Y)^\circ,K')$.
\end{thmcor}

The following theorem was first proved in \cite[Theorem~\ref*{LOT2:thm:Id-is-Id}]{LOT2}, using
computations of the homology of the algebra associated to a pointed matched circle.
The techniques of this paper
lead to another proof, which is independent of those calculations.
\begin{theorem}\label{thm:id-is-id}
  The type \DA\ module $\lsub{\Alg(\PMC)}\CFDAa(\Id)_{\Alg(\PMC)}$
  associated to the identity map of $\Id\co F(\PMC)\to F(\PMC)$ is
  isomorphic to the ``identity bimodule''
  $\lsub{\Alg(\PMC)}\Alg(\PMC)_{\Alg(\PMC)}$.
\end{theorem}

The algebras $\Alg(\PMC)$ have yet more symmetries. To state them, we
give one more operation on pointed matched circles. Given a pointed
matched circle $\PMC$, we can form another pointed matched circle
$\PMC_*$ by turning the Morse function inducing $\PMC$ upside down;
see Construction~\ref{construct:HalfIdentity} for more details. 
\begin{theorem}\label{thm:KoszulDual}
  The algebra $\Alg(\PMC,i)$ is Koszul dual to $\Alg(\PMC,-i)$ and
  also to $\Alg(\PMC_*,i)$. In particular, $\Alg(\PMC,-i)$ is
  quasi-isomorphic to $\Alg(\PMC_*,i)$.
\end{theorem}
(Our algebras have differentials and are not strictly quadratic, so
the definition of Koszul duality in our setting is a modest extension
of the classical one. See Section~\ref{sec:koszul-frame}.)

Theorem~\ref{thm:KoszulDual} explains some seeming coincidences in
the dimensions of $H_*(\Alg(\PMC,i))$; see
Section~\ref{sec:koszul-bordered} for some examples.

As mentioned earlier, a key tool for establishing these results is a Heegaard
diagram discovered independently by Auroux~\cite{AurouxBordered} and
Zarev (see Section~\ref{sec:denis}).
This is a nice
diagram~\cite{SarkarWang07:ComputingHFhat}, so its holomorphic disks
can be understood explicitly; moreover, its combinatorial structure
is closely tied to the bordered Floer algebra. This allows us to
describe differentials in this diagram in a reasonably conceptual
way. Exploiting these properties, we can explicitly describe some modules
which play an important role in the theory. For instance, in
Theorem~\ref{thm:PreciseDD} we give a simple description of the
dualizing bimodule $\CFDDa(\Id)$  (which is also computed, by
different techniques, in~\cite{LOT4}). In a similar vein, we can give
a conceptual description of the Serre functor appearing above and 
a finite dimensional model for the bar complex of $\Alg$ (see
Proposition~\ref{prop:BoundedBar} below).

\subsection{Gradings}
We have stated the results in the introduction without explicitly
discussing the gradings on the modules and bimodules. Typically,
gradings are somewhat subtle in bordered Floer theory; in particular,
the algebras are graded by noncommutative groups and the modules by
$G$-sets. It turns out, however, that these issues do not introduce
any novel features for the results in this paper, beyond those already
present in the pairing theorems from~\cite{LOT1} and~\cite{LOT2}. We
review these issues, in Section~\ref{sec:Gradings}. There, we
also give a detailed statement of how the gradings work in
Theorem~\ref{thm:hom-pair}; graded statements of the other theorems
are similar, and we leave these to the reader.

\subsection{Further remarks}

It is natural to ask what operation in Heegaard Floer homology
corresponds to the composition of homomorphisms. That is, suppose we
have bordered $3$-manifolds $Y_1$, $Y_2$ and $Y_3$ with boundaries
parametrized by some surface $F$. Then there is a composition map
\[
\Ext(\CFDa(Y_1),\CFDa(Y_2))\otimes \Ext(\CFDa(Y_2),\CFDa(Y_3))\to \Ext(\CFDa(Y_1),\CFDa(Y_3)
\]
which corresponds to some homomorphism
\[
\HFa(-Y_1\cup_\bdy Y_2)\otimes \HFa(-Y_2\cup_\bdy Y_3)\to
\HFa(-Y_1\cup_\bdy Y_3).
\]
Generalizing a notion from~\cite{OS04:HolomorphicDisks}, we can use
$Y_1$, $Y_2$ and $Y_3$ to construct a $4$-manifold as follows. Let $T$
denote a triangle, with edges $e_1$, $e_2$ and $e_3$. Then let
\[
W_{Y_1,Y_2,Y_3}=(T\times F)\cup_{e_1\times F}(e_1\times
Y_1)\cup_{e_2\times F}(e_2\times Y_2)\cup_{e_3\times F}(e_3\times Y_3).
\]
Following constructions from~\cite{OS06:HolDiskFour}, this four-manifold
induces a map on Floer homology
  \[
  \hat{F}_{W_{Y_1,Y_2,Y_3}}\co \HFa(-Y_1\cup_\bdy Y_2)\otimes \HFa(-Y_2\cup_\bdy Y_3)\to
  \HFa(-Y_1\cup_\bdy Y_3).
  \]

Under the identifications from Theorem~\ref{thm:hom-pair}, one can show
that this four-manifold invariant corresponds to the composition map  
\[
\Ext(\CFDa(Y_1),\CFDa(Y_2))\otimes \Ext(\CFDa(Y_2),\CFDa(Y_3))\to \Ext(\CFDa(Y_1),\CFDa(Y_3)).
\]
Presumably a similar story holds for the multiplication operations on
Hochschild cohomology.
We return to this point in a future paper~\cite{LOTCobordisms}. 

The duality results for modules and bimodules can be seen as special cases
of results for bordered sutured
manifolds~\cite{Zarev09:BorSut}. In that, more general, context the
presence of boundary Dehn twists can be understood as follows.
If $Y$ is a sutured manifold then orientation reversal of $Y$ also
reverses the roles of $R_+$ and
$R_-$. To accommodate this difference, one must also introduce a
``half Dehn twist'' along the bordered part of the boundary of $Y_1$
(compare Definition~\ref{def:HalfDehn}) before
gluing it to $Y_2$. Though we will generally not discuss this case, see
Remark~\ref{rmk:borsut} for a little more on subtleties in the
bordered sutured context.

In a different direction, one must proceed with care in adapting
Theorem~\ref{thm:hom-pair} to the case of bordered sutured
manifolds~\cite{Zarev09:BorSut}.  This extension is developed (in a
slightly different language) in~\cite{Zarev:JoinGlue}; see also
Remark~\ref{rmk:borsut}.

\subsection{Outline of the paper}
In Section~\ref{sec:review} we
review some basic facts about bordered Floer theory and homological
algebra which are used in the rest of this paper.  A key point in the
present paper is to generalize the arced bordered Heegaard diagrams
from~\cite{LOT2} to the case where both $\alpha$ and $\beta$
curves go out to the boundary, and to give a topological meaning to
these objects.  This is done in Section~\ref{sec:ab}. 
In Section~\ref{sec:denis} we describe a particular Heegaard
diagram whose associated bimodule is the algebra itself.  This lets us
give Heegaard-diagrammatic interpretations to some of the algebraic
operations on the various modules in bordered Floer theory.  In
Section~\ref{sec:Consequences}, we collect the consequences of these
interpretations (together with the traditional pairing theorem) to
prove the results stated in the introduction. In
Section~\ref{sec:Gradings} we give a brief discussion of how gradings
can be added to the present context. In Section~\ref{sec:Examples} we
illustrate some of the above discussion with some examples. Finally,
in Section~\ref{sec:koszul}, we give a description of Koszul duality
relevant for our algebras, and prove Theorem~\ref{thm:KoszulDual}.

A summary of some of the conventions employed in the paper can be
found in Appendix~\ref{sec:conventions}.

\subsection{Acknowledgements}
The authors thank Denis Auroux and Rumen Zarev for helpful
conversations. In particular, Auroux's paper~\cite{AurouxBordered} led
us to a dramatically simplified argument. The relevant diagram of
Auroux was discovered (and communicated to us) independently by
Zarev. As noted above, Theorem~\ref{thm:hom-pair} was discovered
independently by Auroux, and a form of Theorem~\ref{thm:D-is-A} was
also suggested to us by him.  The authors also thank MSRI for its
hospitality during the completion of this project. Finally, we thank
the referees for a detailed reading and many helpful and interesting comments.


\section{Background}\label{sec:review}
\subsection{Basic structure of bordered Floer theory}\label{sec:HF-basics}
We start reviewing the key ingredients from~\cite{LOT1} and~\cite{LOT2}.

A \emph{matched circle} is an oriented circle $Z$, $2k$ points
$\{a_1,\dots,a_{2k}\}=\mathbf{a}\subset Z$ and a fixed-point free
involution $M\co \mathbf{a}\to\mathbf{a}$. A matched circle
$(Z,\mathbf{a},M)$ specifies a surface-with-boundary
$F^\circ(Z,\mathbf{a},M)$ by filling in $Z$ with a disk $D_0$ and
attaching $2$-dimensional $1$-handles at each pair
$\{a_i,a_{M(i)}\}\subset \mathbf{a}$. We shall only be interested
in matched circles which specify surfaces with a single boundary
component. We can fill in the boundary component of
$F^\circ(Z,\mathbf{a},M)$ to give a closed surface $F$, with a
distinguished disk $D=F\setminus F^\circ$. We orient $F$ so that the
orientation of $D_0$ induces the orientation of $Z$.

A \emph{pointed matched circle} is a matched circle together with a
basepoint $z\in Z\setminus \mathbf{a}$. We shall use the notation
$\PMC$ to denote a pointed matched circle $(Z,\mathbf{a},M,z)$. A
pointed matched circle $\PMC$ specifies a closed surface $F(\PMC)$
together with a distinguished disk $D\subset F(\PMC)$ and basepoint
$z\in\bdy D$. (This construction will be expanded slightly in
Section~\ref{subsec:BetaMatched}.)

Bordered Floer homology associates a \dg algebra $\Alg(\PMC)$ to each
pointed matched circle~\cite{LOT1}. If $F(\PMC)\cong F(\PMC')$ then the
algebras $\Alg(\PMC)$ and $\Alg(\PMC')$ are derived equivalent,
according
to~\cite[Theorem~\ref*{LOT2:thm:AlgebraDependsOnSurface}]{LOT2}.

A \emph{bordered $3$-manifold} is a quadruple
$(Y_1,\Delta_1,z_1,\psi_1)$, where $Y_1$ is a
three-manifold-with-boundary, $\Delta_1$ is a disk in $\partial Y_1$,
$z_1$ is a point on $\partial \Delta_1$, and
$$\psi\co (F(\PMC),D,z) \to (\partial Y_1,\Delta_1,z_1),$$
is a homeomorphism from the surface $F(\PMC)$ (for some pointed
matched circle $\PMC$) to $\partial Y_1$ sending $D$ to $\Delta_1$ and
$z$ to $z_1$. We will often suppress the preferred disk $D$ and
basepoint $z$ from the notation for a bordered $3$-manifold, and will
sometimes also suppress the homeomorphism $\psi$.

As explained in~\cite{LOT1},
bordered Floer homology associates to a bordered $3$-manifold with
boundary parameterized by $F(\PMC)$ a right $\Ainf$ $\Alg(\PMC)$-module
$\CFAa(Y)$ and a left \dg $\Alg(-\PMC)$-module $\CFDa(Y)$, each
well-defined up to quasi-isomorphism. (Here, $-\PMC$ denotes $\PMC$
with its orientation reversed.)

Bordered Floer homology also works for $3$-manifolds with more than
one boundary component; see~\cite{LOT2}. More precisely, a {\em
  strongly bordered three-manifold with boundary $F(\PMC_1)\amalg
  F(\PMC_2)$} is an oriented three manifold $Y_{12}$ with boundary,
equipped with
\begin{itemize}
  \item preferred disks $\Delta_1$ and $\Delta_2$ on its two boundary
    components,
  \item basepoints $z_i'\in \partial \Delta_i$,
  \item a homeomorphism $\psi\co 
    (F(\PMC_1)\amalg F(\PMC_2),D_1\amalg D_2,z_1\amalg z_2)
    \to (\partial Y_{12},\Delta_1\amalg \Delta_2,z_1'\amalg z_2')$,
  \item an arc $\gamma$ connecting $z_1'$ to $z_2'$, and
  \item a framing of $\gamma$, pointing into $\Delta_i$ at 
    $z_i'$ for $i=1,2$.
\end{itemize}
We will often denote the two boundary components of a strongly
bordered three-manifold $Y$ by $\bdy_LY$ and $\bdy_RY$ (for ``left''
and ``right''), but the choice of which boundary component is $\bdy_L$
and which is $\bdy_R$ is arbitrary.

For every strongly bordered $3$-manifold
with boundary $F(\PMC_1)\amalg F(\PMC_2)$, bordered Floer homology
associates an $\Ainf$-bimodule
$\CFAAa(Y)$ with right actions by $\Alg(\PMC_1)$ and $\Alg(\PMC_2)$;
an $\Ainf$-bimodule $\CFDAa(Y)$ with a left action by $\Alg(-\PMC_1)$
and a right action by $\Alg(\PMC_2)$; and a \dg bimodule $\CFDDa(Y)$
with left actions by $\Alg(-\PMC_1)$ and $\Alg(-\PMC_2)$. Each of
$\CFAAa(Y)$, $\CFDAa(Y)$ and $\CFDDa(Y)$ is well-defined up to
quasi-isomorphism. As a special case, if $\PMC$ is a pointed matched
circle, $F(\PMC)\times[0,1]$ is naturally a strongly bordered
three-manifold, which we will denote $\Id_\PMC$ or just $\Id$.

The bordered Floer modules relate to each other and to the closed
invariant $\CFa(Y)$ by \emph{pairing theorems}. The prototypical
pairing theorem (see~\cite[Theorem~\ref{LOT:thm:TensorPairing}]{LOT1})
states that if $Y_1$ and $Y_2$ are bordered $3$-manifolds with $\bdy
Y_1=F(\PMC)=-\bdy Y_2$ then
\[
\CFa(Y_1\cup_F Y_2)\simeq \CFAa(Y_1)\DTP_{\Alg(\PMC)}\CFDa(Y_2).
\]
Here, $\DTP$ denotes the derived (or $\Ainf$) tensor product. The
analogues for bimodules are listed in~\cite{LOT2}; the mnemonic is
that one can cancel expressions of the form $A\DTP D$. For instance,
if $Y_1$ is a bordered $3$-manifold with boundary $F(\PMC_1)$ and
$Y_{12}$ and $Y_{23}$ are strongly bordered $3$-manifolds with
boundaries $-F(\PMC_1)\amalg F(\PMC_2)$ and $-F(\PMC_{2})\amalg
F(\PMC_3)$ respectively then
\begin{align*}
  \CFDa(Y_1\cup_{F(\PMC_1)}Y_{12})\simeq
  \CFAa(Y_1)\DTP_{\Alg(\PMC_1)}\CFDDa(Y_{12}),\\
  \CFDAa(Y_{12}\cup_{F(\PMC_2)}Y_{23})\simeq \CFAAa(Y_{12})\DTP_{\Alg(\PMC_2)}\CFDDa(Y_{23}),
\end{align*}
and so on.

The details of the construction of the algebras $\Alg(F)$ and the
modules $\CFDa(Y)$ and $\CFAa(Y)$ can be found in~\cite{LOT1}, and the
generalization to the case of more boundary components is
in~\cite{LOT2}. Much of this paper can be read with merely a cursory
understanding of~\cite{LOT1} and~\cite{LOT2}, keeping the following
points in mind:
\begin{itemize}
\item Suppose that $\PMC$ is a pointed matched circle and $-\PMC$
  denotes the same data except with the orientation of the circle
  reversed. Then $F(-\PMC)=-F(\PMC)$, and the algebras are related by:
  \begin{equation}
    \Alg(\PMC)^{\mathrm{op}}\cong \Alg(-\PMC).\label{eq:OpAlg}
  \end{equation}
\item The modules $\CFDa(Y)$ and $\CFAa(Y)$ are not associated
  directly to the $3$-manifold $Y$ but rather to a
  \emph{bordered Heegaard diagram} for $Y$, i.e., a Heegaard diagram
  \[
  (\Sigma'_g,\alphas^c=\{\alpha^c_1,\dots,\alpha^c_{g-k}\},\betas^c=\{\beta_1,\dots,\beta_g\},z)
  \]
  for $Y$ together with $2k$ disjoint, embedded arcs
  $\alphas^a=\{\alpha^a_1,\dots,\alpha^a_{2k}\}$ in
  $\overline{\Sigma}=\Sigma'\setminus\bD^2$ with boundary on
  $\bdy\overline{\Sigma}$ giving a basis for $\pi_1(\bdy Y)$; and a
  basepoint $z\in\bdy\overline{\Sigma}$ not lying on any
  $\alpha_i^a$. The boundary of such a diagram is a pointed matched circle.

  Similarly, the bimodules $\CFDDa(Y)$, $\CFDAa(Y)$ and
  $\CFAAa(Y)$ associated to a strongly bordered $3$-manifold with two
  boundary components are associated to \emph{arced bordered Heegaard
    diagrams} with two boundary components; see~\cite{LOT2} and also
  Section~\ref{sec:abDiagrams} below.
\item The module $\CFDa(Y)$ has a special form: it is a \emph{type $D$
    structure}: 
  \begin{definition}\label{def:type-D-str}
    Let $\Alg$ be a \dg algebra over a ring
    $\Ground=\bigoplus_{i=1}^N\FF_2$. A (left) \emph{type $D$ structure}
    over
    $\Alg$ is a $\Ground$-module $X$
    equipped with a map
    $\delta^1\co X\to \Alg\otimes_{\Ground} X$
    satisfying the structure equations which ensure
    that $\delta^1$ extends via the Leibniz rule to give 
    $\Alg\otimes_{\Ground} X$ the structure of a differential
    $\Alg$-module.
  \end{definition}
  (Type $D$ structures can be thought of as differential comodules or
  twisted complexes, see
  Remarks~\cite[Remark~\ref*{LOT2:rmk:Comodule}]{LOT2}
  and~\cite[Remark~\ref*{LOT2:rmk:TwistedComplex}]{LOT2},
  respectively.)

  There is a convenient model $\DT$ for the $\Ainf$-tensor product of
  an $\Ainf$-module $L$ and a type $D$ structure $M=(X,\delta^1)$,
  with the property that the vector space underlying $L\DT M$ is just
  $L\otimes_{\Ground}X$; see~\cite[Definition~\ref*{LOT1:def:DT}]{LOT1}. (In the
  case that $L$ is an ordinary module, $L\DT M$ agrees with the
  na\"ive tensor product $L\otimes_{\Alg} M$. In particular,
  when $L$ is $\Alg$ viewed as a bimodule, $\Alg\DT M$
  is the module associated to $X$.) There are analogues of the
  operation $\DT$ for bimodules, as well;
  see~\cite[Section~\ref*{LOT2:sec:tensor-products}]{LOT2}. (Although
  $\DT$ has
  a purely algebraic definition, it arises naturally in the analysis
  of pseudoholomorphic 
  curves, as seen in the proof of
  the pairing theorem~\cite[Chapter~\ref*{LOT1:chap:tensor-prod}]{LOT1}.)

  Similarly, the bimodule $\CFDDa(Y)$ is a type \DD\ structure, i.e., a
  type $D$ structure over $\Alg(\PMC_L)\otimes\Alg(\PMC_R)$; and the
  bimodule $\CFDAa(Y)$ is a type \DA\ structure, as defined
  in~\cite[Definition~\ref*{LOT2:def:DA-structure}]{LOT2}. The
  operation $\DT$ works when
  tensoring bimodules, as long as one tensors a type $D$ side with a
  type $A$ (i.e., $\Ainf$) side.
\end{itemize}

We sometimes blur the distinction between a type $D$ structure
$(X,\delta^1)$ and its induced differential module $\Alg\DT X$.  When
it is important to distinguish them, we include a superscript in the
notation for a type $D$ structure, $\lsup{\Alg}X$. Ordinary modules
are indicated with a subscript; so we sometimes use the notation
$\lsub{\Alg}X$ to denote the associated module $\Alg\DT
\lsup{\Alg}X$. This operation has an inverse (``raising the
subscript'') which associates to a module $M_{\Alg}$ the type $D$
structure $M_{\Alg}\DT \lsupv{\Alg}\Barop(\Alg)^{\Alg}$, the {\em bar
  resolution} of $M$. These are inverses in the derived category of
modules satisfying suitable boundedness conditions;
see~\cite[Proposition~\ref*{LOT2:prop:D-to-A-and-back}]{LOT2}.

\subsection{Review of \textalt{$\Mor$}{Mor} and \textalt{$\Ext$}{Ext}}\label{sec:RHom-Ext}
Suppose that $C_*$ and $D_*$ are chain complexes (or differential
modules) over an algebra $A$ (possibly with differential) which we assume to have characteristic $2$.  Two ways to compute
$\Ext_R(C_*,D_*)$ are:
\begin{enumerate}
\item Find a complex $C'_*$ of projective modules quasi-isomorphic to
  $C_*$ and compute the homology of the complex
  $\Hom_A(C'_*,D_*)_k=\bigoplus_i \Hom(C'_i,D_{k+i})$ of maps
  respecting the module structure (but not necessarily the differential)
  from $C'_*$ to $D_*$, or
\item Take the homology of the chain complex $\Mor_A(C_*,D_*)$ of
  $\Ainf$-morphisms from $C_*$ to~$D_*$.
\end{enumerate}
(The second option is (under some finite-dimensionality assumptions) a
special case of the first: the complex of $\Ainf$-morphisms from $C_*$
to $D_*$ is exactly the chain complex of homomorphisms from the bar
resolution of $C_*$ to $D_*$; see Section~\ref{sec:dualizing}.)

Of course, if $C_*$ is already projective, one can take $C'_*$ to just
be $C_*$ itself.  Given type $D$ structures $\lsup{\Alg}M$ and
$\lsup{\Alg}N$, define
\[
\Mor^\Alg(\lsup{\Alg}M,\lsup{\Alg}N)_k=\bigoplus_i \Hom (\lsub{\Alg}M_i,\lsub{\Alg}N_{k+i}),
\]
with its obvious differential. Type $D$ structures over our algebras correspond to
projective modules
(compare~\cite[Corollary~\ref*{LOT2:cor:tensor-projective}]{LOT2}),
so
$\Ext_*(\lsub{\Alg}M,\lsub{\Alg}N)\cong H_*(\Mor^\Alg(\lsup{\Alg}M,\lsup{\Alg}N))$.

Again: the notation $\Mor_\Alg$, with a subscript, denotes the complex
of $\Ainf$-morphisms between $\Ainf$-modules, while the notation
$\Mor^\Alg$, with a superscript, denotes the complex of module maps of
type $D$ structures. In either case, the homology of the
$\Mor$ complex computes $\Ext$.

The following result of~\cite{LOT2} (an easy consequence
of~\cite[Theorems~\ref*{LOT2:thm:Id-is-Id}
and~\ref*{LOT2:thm:GenComposition}]{LOT2},
see~\cite[Lemma~\ref*{LOT2:lemma:AADDquasi-equiv}]{LOT2}) will
reduce our work by roughly half: 
\begin{citethm}\label{thm:DDId}Fix a pointed matched circle $\PMC$. Then the functors
  \begin{align*}
    \cdot\DT\CFDDa(\Id)\co H_*(\ModCat_{\Alg(\PMC)})&\to
    H_*(\lsup{\Alg(-\PMC)}\ModCat)\\
    \CFAAa(\Id)\DT\cdot \co H_*(\lsup{\Alg(-\PMC)}\ModCat)&\to
    H_*(\ModCat_{\Alg(\PMC)})
  \end{align*}
  are inverse equivalences of categories, exchanging $\CFAa(Y)$ and
  $\CFDa(Y)$. 
\end{citethm}
Here, $H_*(\ModCat_{\Alg(\PMC)})$ (respectively $
H_*(\lsup{\Alg(-\PMC)}\ModCat)$) denotes the homotopy category
  of right, $\Ainf$ (respectively left, type $D$) modules over
  $\Alg(\PMC)$ (respectively
  $\Alg(-\PMC)$). Recall also that the homotopy categories of
  $\Ainf$-modules and projective
  modules are both equivalent to the derived category; see \cite[Section
  \ref*{LOT2:sec:models-der-cat}]{LOT2}.

\begin{corollary}\label{cor:D-hom-is-A-hom}
  Fix bordered $3$-manifolds $Y_1$ and $Y_2$ with $\bdy Y_1=F=\bdy
  Y_2$. Then there is a 
  quasi-isomorphism
  \[
  \Mor_{\Alg(F)}(\CFAa(Y_1),\CFAa(Y_2))\simeq \Mor^{\Alg(-F)}(\CFDa(Y_1),\CFDa(Y_2)).
  \]
\end{corollary}

Although we do not need it for our present purposes, the module
$\CFDDa(\Id)$ is described explicitly in Section~\ref{sec:dualizing-bimodules}
or, via a different method, in~\cite{LOT4}.
\subsection{Duals of modules and type $D$
  structures}\label{sec:dualizing}
As for finite-dimensional vector spaces, where $\Hom(V,W)\cong
W\otimes V^*$, we can interpret our $\Mor$ complexes in terms of
tensor products and duals. We spell this out explicitly.
\begin{definition}\label{def:dual-type-D}
  Let $\Alg$ be a finite-dimensional \dg algebra over $\Ground=\bigoplus_{i=1}^N\FF_2$
  and $M=(X,\delta^1)$ a left type $D$ structure over $\Alg$. Let
  $\overline{X}=\Hom_\Field(X,\Field)$ denote the dual of $X$. The
  transpose of $\delta^1$ is a map
  \[
  (\delta^1)^T\co \overline{X}\otimes\Hom_\Field(\Alg,\Field)\to \overline{X}.
  \]
  We can interpret this instead as a map 
  \[
  \overline{\delta}^1\co \overline{X}\to \overline{X}\otimes\Alg.
  \]
  The \emph{dual type $D$ structure} to $M$,
  $\overline{M}$, is the right type $D$ structure induced by
  $(\overline{X},\overline{\delta}^1)$.
\end{definition}
If we draw the operation on a type $D$ structure $M$ like this:
\[
\begin{tikzpicture}
  \node at (0,0) (tblank) {};
  \node at (0,-1) (delta) {$\delta^1$};
  \node at (-1,-2) (blblank) {};
  \node at (0,-2) (bcblank) {};
  \draw[Dmodar] (tblank) to node[right]{\lab{X}} (delta);
  \draw[Dmodar] (delta) to node[right]{\lab{X}} (bcblank);
  \draw[algarrow, bend right=30] (delta) to node[above,sloped]{\lab{\Alg}} (blblank);
\end{tikzpicture}
\]
(compare~\cite[Section~\ref*{LOT2:sec:algebra-modules}]{LOT2}) then 
the type $D$ structure $\overline{M}$ is
\[
\begin{tikzpicture}
  \node at (0,0) (tblank) {};
  \node at (0,-1) (delta) {$\overline{\delta}^1$};
  \node at (1,-2) (brblank) {};
  \node at (0,-2) (bcblank) {};
  \draw[Dmodar] (bcblank) to (delta);
  \draw[Dmodar] (delta) to (tblank);
  \draw[algarrow, bend left=30] (delta) to (brblank);
\end{tikzpicture}.
\]
(As is standard in such graphical calculus, arrows pointing up
represent the dual modules of arrows pointing down.)

\begin{lemma}
  If $M=(X,\delta^1)$ is a type $D$ structure then
  $\overline{M}=(\overline{X},\overline{\delta}^1)$ is also a type $D$
  structure.
\end{lemma}
\begin{proof}
  This is a straightforward exercise in the properties of duals;
  alternately, it is clear from the graphical description.
\end{proof}

\begin{proposition}\label{prop:D-Mor-is-DT}
  Suppose that $M$ and $N$ are type $D$ structures and at least one of
  $M$ or $N$ is finite-dimensional. Then the chain complex of module
  homomorphisms
  from $M$ to $N$, $\Mor^{\Alg}(M,N)$, is isomorphic to
  $\overline{M}\DT_{\Alg}\Alg\DT_\Alg N$.
  In particular, $\Mor_{\Alg}(\lsub{\Alg}M,\Alg)\simeq \overline{M}\DT_\Alg\Alg$.
\end{proposition}
\begin{proof}
  The first part is straightforward from the definitions.  For the
  last statement, consider the type $D$ structure $\lsup{\Alg}T =
  (\Ground, 0)$. (That is, $T$ is rank one and $\delta^1_T = 0$).  Then
  $\lsub{\Alg}T = \Alg$ and $\overline{M} \DT_\Alg \Alg \cong
  \Mor^{\Alg}(\lsup{\Alg}M, T) \simeq \Mor_{\Alg}(\lsub{\Alg}M, \Alg)$.
\end{proof}

We next turn to duals of $\Ainf$-modules.
\begin{definition}\label{def:dual-type-A}
  Let $\Alg$ be a finite-dimensional \dg algebra over $\Ground=\bigoplus_{i=1}^N\FF_2$ and $M$
  a right $\Ainf$-module over $\Alg$. Let
  $\overline{M}=\Hom_\Field(M,\Field)$. The higher multiplications
  $m_{i+1}\co M\otimes \Alg^{\otimes i}\to M$ dualize to give maps $m_{i+1}^T\co \overline{M}\to
  \Hom(\Alg,\Field)^{\otimes i}\otimes \overline{M}$, which we can
  interpret as maps 
  \[
  \overline{m}_{i+1}\co \Alg^{\otimes i}\otimes \overline{M}\to \overline{M}.
  \]
  Then the data $(\overline{M},\{\overline{m}_{i+1}\})$ is the
  \emph{dual $\Ainf$-module to $M$}.
\end{definition}
\begin{lemma}
  If $(M,\{m_{i+1}\})$ is a right $\Ainf$-module then
  $(\overline{M},\{\overline{m}_{i+1}\})$ also satisfies the
  $\Ainf$-module relation to make $\overline{M}$ a left $\Ainf$-module.
\end{lemma}
\begin{proof}
  Again, this is a straightforward exercise in the properties of duals.
\end{proof}

Morphism spaces of type $A$ modules can also be described in terms
of~$\DT$. Before giving the definition, we recall some notation. An
\emph{augmentation} of a \dg algebra $\Alg$ is a homomorphism
$\epsilon\co \Alg\to \Ground$ from the algebra to the ground
ring. Given an augmentation $\epsilon$ of a \dg algebra
$\Alg=(A,\mu,d)$, let $A_+=\ker(\epsilon)$ denote the
\emph{augmentation ideal}. There is a type \DD\ bimodule
$\lsup{\Alg}\rBarop(\Alg)^\Alg$ with underlying $\Ground$-module
$\Tensor^*(A_+)$, with basis written $[a_1|\cdots|a_k]$ for $k
\ge 0$, and structure maps
\begin{multline*}
  \delta^1[a_1|\cdots|a_k] \coloneqq a_1\otimes[a_2|\cdots|a_k]\otimes 1
  + 1\otimes [a_1|\cdots|a_{k-1}]\otimes a_k\\
  +\sum_{1\le i\le k}1\otimes
  [a_1|\cdots|d(a_i)|\cdots|a_k]\otimes 1 +
  \sum_{1\le i\le k-1}1\otimes
  [a_1|\cdots|\mu(a_i,a_{i+1})|\cdots|a_k]\otimes 1.
\end{multline*}
The bimodule $\lsup{\Alg}\rBarop(\Alg)^\Alg$ is called the
\emph{reduced bar complex of $\Alg$}. (In~\cite{LOT2}, we typically
worked with the unreduced bar complex $\Barop(\Alg)$. The canonical
inclusion $\rBarop(\Alg)\to \Barop(\Alg)$ is a homotopy equivalence.)

With this terminology in hand, we have the following reformulation of
the complex of $\Ainf$-module homomorphisms:
\begin{proposition}\label{prop:A-Mor-is-DT}
  For finite-dimensional right $\Ainf$-modules $M$ and~$N$, the chain
  complex of $\Ainf$-module homomorphisms
  from $M$ to $N$, $\Mor_{\Alg}(M,N)$, is isomorphic to
  $N \DT_\Alg \overline{\rBarop(\Alg)} \DT_\Alg \overline{M}$.
\end{proposition}

\begin{proof}
  Recall that an $\Ainf$-morphism $f\co M \to N$ consists of maps
  $f_{i+1} \co M \otimes A[1]^{i} \to N$.  By unitality we can
  restrict the algebra inputs to lie in $A_+$, and because $M$ and~$N$
  are finite-dimensional, the space of such maps (as a vector space)
  is isomorphic to $N \otimes \Hom(T^*(A_+[1]),\Field) \otimes \Hom(M,
  \Field)$, which is the underlying space of $N \DT_\Alg
  \overline{\rBarop(\Alg)} \DT_\Alg \overline{M}$.  Checking that the
  differentials on the two complexes agree is again elementary.
\end{proof}

A key property of the bar complex is that it can be used to give
resolutions of modules. In the present language, this boils down to
the following identity~\cite[\ref*{LOT2:lem:bar-res}]{LOT2}:
\begin{equation}
  \label{eq:bar-resolution}
  \lsupv{\Alg}\rBarop^\Alg \DT\Alg\simeq \lsupv{\Alg}[\Id]_\Alg.
\end{equation}
Here, $\lsupv{\Alg}[\Id]_\Alg$ is the type \DA\ bimodule whose
modulification is $\Alg$;
see~\cite[Definition~\ref*{LOT2:def:rank-1-DA-mods}]{LOT2}. In
particular, $\lsupv{\Alg}[\Id]_\Alg$ is the identity for $\DT$.

We extend the definitions of duals to bimodules as follows:
\begin{definition}\label{def:dual-bimod}
  Let $\Alg$ and $\Blg$ be finite-dimensional \dg algebras over
  $\Ground=\bigoplus_{i=1}^N\FF_2$ and
  $\Groundl=\bigoplus_{i=1}^M\FF_2$. 

  Suppose that $\lsup{\Alg,\Blg}M$ is a left-left type \DD\ module
  over $\Alg$ and $\Blg$. That is, $\lsup{\Alg,\Blg}M$ is a type $D$
  module over $\Alg\otimes\Blg$. Then $\overline{\lsup{\Alg,\Blg}M}$
  is a right type $D$ structure over $\Alg\otimes\Blg$. Interpreting
  this as a right-right type \DD\ structure
  $\overline{M}{}^{\Alg,\Blg}$, we call $\overline{M}{}^{\Alg,\Blg}$ the
  \emph{dual type \DD\ structure to $\lsup{\Alg,\Blg}M$.}

  Now, suppose that $N_{\Alg,\Blg}$ is a right-right $\Ainf$-bimodule
  over $\Alg$ and $\Blg$. Let $\overline{N}=\Hom(M,\Field)$. The
  transpose of the
  structure maps $m_{1,i,j}\co N\otimes\Alg^{\otimes
    i}\otimes\Blg^{\otimes j}\to N$ are maps $m_{1,i,j}^T\co
  \overline{N}\to \Hom(\Blg,\Field)^{\otimes
    j}\otimes\Hom(\Alg,\Field)^{\otimes i}\otimes\overline{N}$, which we
  interpret as maps 
  \[\overline{m}_{j,i,1}\co \Blg^{\otimes j}\otimes
  \Alg^{\otimes i}\otimes\overline{N}\to \overline{N}.\]
  We call
  $(\overline{N},\{\overline{m}_{j,i,1})$ the \emph{dual
    $\Ainf$-bimodule to $N_{\Alg,\Blg}$}.

  Finally, suppose $\lsup{\Alg}N_\Blg=(N,\delta^1)$ is a type \DA\
  structure. Let $\overline{N}=\Hom(N,\Field)$. The transpose of the
  map $\delta^1_{1+n}\co
  N\otimes \Blg^{\otimes n}\to \Alg\otimes M$ is a map
  $(\delta^1_{1+n})^T\co \overline{N}\otimes \Hom_{\Field}(\Alg,\Field)\to
  \Hom_{\Field}(\Blg,\Field)^{\otimes n}\overline{N}$. We can interpret
  these as maps
  \[
  \overline{\delta}^1_{1+n}\co  \Blg^{\otimes n}\otimes \overline{N}\to
  \overline{N}\otimes \Alg.
  \]
  We call $(\overline{N},\overline{\delta}^1_{1+n})$ the \emph{dual type
    \DA\ structure} to $\lsup{\Alg}N_\Blg$.
\end{definition}

The following is an easy exercise.
\begin{lemma} The dual type \DD\ structure to a type \DD\ structure
  satisfies the type \DD\ structure equation. The dual
  $\Ainf$-bimodule to an $\Ainf$-bimodule satisfies the
  $\Ainf$-structure equation. The dual type \DA\ structure to a \DA\
  structure satisfies the \DA\ structure equation.
\end{lemma}

For bimodules, we can consider the space of morphisms over either one
or both of the actions. So, Proposition~\ref{prop:D-Mor-is-DT}
corresponds to two different statements for bimodules:
\begin{proposition}\label{prop:Mor-is-DT-bimod}
  Let $\lsup{\Alg,\Blg}M$ and $\lsup{\Alg,\Blg}N$ be type \DD\
  structures, at least one of which is finite-dimensional. Then
  \[
  \Mor^{\Alg\otimes\Blg}(M,N)\cong
  \overline{M}\DT_{\Alg\otimes\Blg}(\Alg\otimes\Blg)\DT_{\Alg\otimes\Blg} N.
  \]

  Similarly, let $\lsup{\Alg,\Blg}M$ be a type \DD\ structure and let
  $\lsup{\Blg}N$ a type $D$ module. Then
  \[
  \Mor^\Blg(M,N)\cong \overline{M} \DT_{\Blg} \Blg \DT_\Blg N,
  \]
  as type $D$ structures.  A corresponding statement holds if $N$ is a
  type \DD\ or \DA\ module.
\end{proposition}
\begin{proof}
  Like Proposition~\ref{prop:D-Mor-is-DT}, this is immediate from the
  definitions. (In fact, the first half also follows from
  Proposition~\ref{prop:D-Mor-is-DT}.)
\end{proof}

\begin{lemma}\label{lem:dual-resp-DT}
  Taking duals respects the operation $\DT$ in the following sense: if
  $M_{\Alg}$ is an $\Ainf$-module and $\lsup{\Alg}N$ is a
  type $D$ structure then
  \begin{equation}
    \overline{M\DT_{\Alg} N}\cong \overline{N}\DT_{\Alg}\overline{M}.\label{eq:dual-resp-DT}
  \end{equation}
  Moreover, if $M$ is an $\Ainf$-bimodule or type \DA\ structure; and
  $N$ is a type \DA\ or \DD\ structure then the isomorphism in
  Equation~\eqref{eq:dual-resp-DT} is an isomorphism of \AAm, \DA, or
  \DD\ structures (as appropriate).
\end{lemma}
\begin{proof}[Proof sketch]
  To get the structure maps on $\overline{M\DT N}$, take
  the appropriate diagram from
  \cite[Figure~\ref*{LOT2:fig:bimod-on-DT}]{LOT2}, rotate it
  $180^\circ$, and modify the diagram so all the algebra arrows point
  down.  Then the module arrows are pointing up, as appropriate for
  diagrams involving the dual bimodule, and the diagrams are the same
  as those for the $\overline{N} \DT \overline{M}$, as desired.
\end{proof}

In this paper, we will pass freely between left modules over an algebra and
right modules over the opposite algebra. Specifically,
a right module $N_{\Alg}$ over $\Alg$ can be viewed as
a left module over $\Alg^{\op}$, which we write as $\lsub{\Alg^{\op}}N$,
and a 
left type $D$ structure $\lsup{\Alg}M$ can be naturally viewed as 
a right type $D$ structure over $\Alg^{\op}$, which we will write as
$M^{\Alg^{\op}}$. The following is straightforward.
\begin{lemma}
  \label{lem:Opposites}
  There is an isomorphism
  $$(\Alg\DT \lsup{\Alg}M)_{\Alg^{\op}}\cong (M^{\Alg^{\op}}\DT \Alg^{\op})$$
\end{lemma}


\section{\textalt{$\alpha$-$\beta$}{alpha-beta}-bordered Heegaard Diagrams}\label{sec:ab}
\subsection{\textalt{$\beta$}{beta} pointed matched circles}
\label{subsec:BetaMatched}

In~\cite{LOT1} (as reviewed in Section~\ref{sec:HF-basics}), we gave a
convention for how a pointed matched circle specifies a surface. In
that paper, we considered exclusively the case where it was the
$\alpha$-curves which ran out into the boundary. In the present
paper, we will need to place the $\alpha$- and $\beta$-curves on a
more equal footing.  Consequently, we would like to give a less biased
construction, in the spirit of Zarev's work on sutured manifolds~\cite{Zarev09:BorSut}.

In this more symmetric construction, we make the following cosmetic revision to the notion of a pointed matched circle:
\begin{definition}\label{def:decorated-pmc}
  A {\em decorated pointed
    matched circle} consists of the following data:
  \begin{itemize}
  \item a circle $Z$;
  \item a decomposition of $Z$ into two closed oriented intervals,
    $Z_\alpha$ and $Z_\beta$, whose intersection consists of two
    points, 
    and with $Z_\alpha$ and $Z_\beta$ oriented
    opposingly;
  \item a collection of $4k$ points ${\mathbf p}=\{p_1,\dots,p_{4k}\}$
    in $Z$, so that either ${\mathbf p} \subset Z_\alpha$ or ${\mathbf
      p} \subset Z_\beta$;
  \item a fixed-point-free involution $M$ on the points in ${\mathbf
      p}$; and
  \item a decoration by the letter $\alpha$ or the letter $\beta$,
    which indicates whether the points ${\mathbf p}$ lie in $Z_\alpha$
    or $Z_\beta$.
  \end{itemize}
  We require that the points $\mathbf{p}$ and involution $M$ satisfy
  the condition that performing surgery on $Z$ along the zero-spheres
  specified by ${\mathbf p}/M$ gives a connected $1$-manifold.  We
  abbreviate the data $(Z=Z_\alpha\cup Z_\beta,{\mathbf p},M)$ by
  $\PMC^{\alpha}$ or $\PMC^{\beta}$ (depending on the decoration); we
  will call the resulting object an \emph{$\alpha$-pointed matched
    circle} or a \emph{$\beta$-pointed matched circle}.
\end{definition}

\begin{figure}
  \centering
  \includegraphics[width=\textwidth]{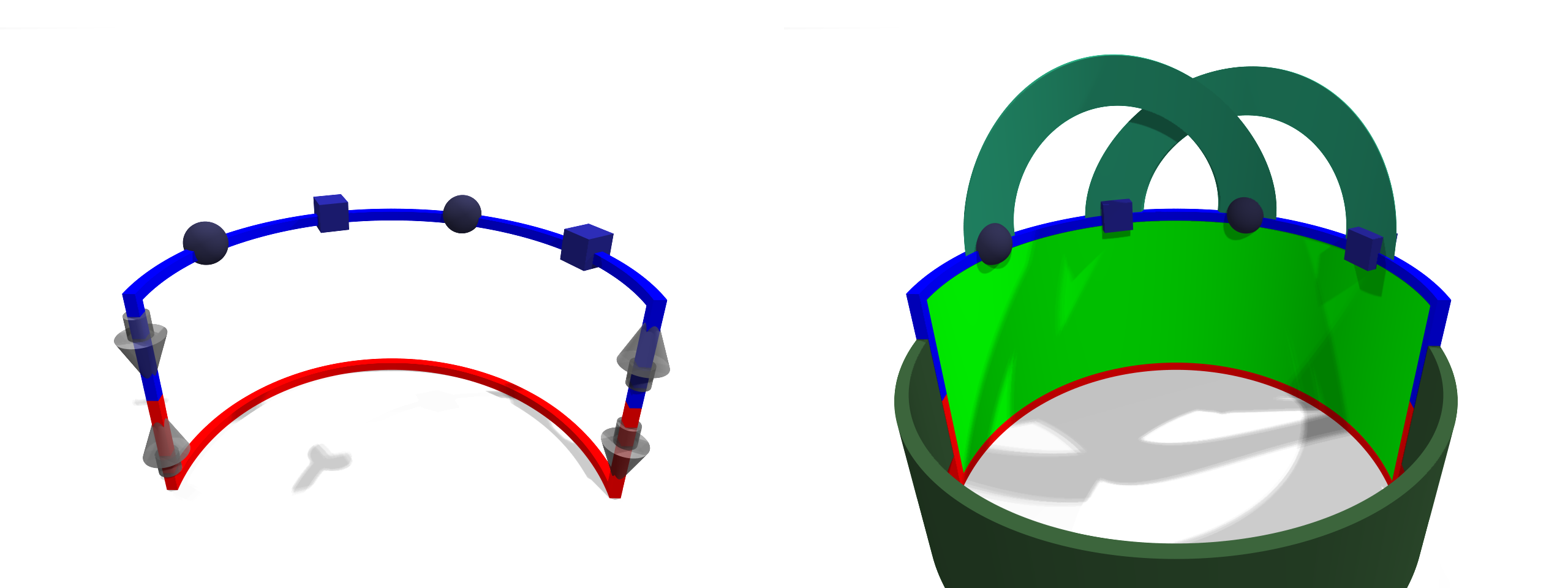}
  \caption{\textbf{Constructing a surface from a $\beta$-pointed
      matched circle.} Left: a $\beta$-pointed matched circle $\PMC^\beta$. The
    bottom half is the arc $Z_\alpha$ and the top half is
    the arc $Z_\beta$; translucent arrows indicate the orientations of
    $Z_\alpha$ and $Z_\beta$. Here, $k=1$, and the matching is indicated by
    the different shapes along the top. Right: part of the resulting surface
    $F(\PMC^\beta)$. The disk $D_0$ is in the back, and the $1$-handle
    $s$ wraps around the front.  $D_\alpha$ and $D_\beta$ are not
    shown.}
  \label{fig:beta-pmc-to-surf}
\end{figure}

\begin{construction}\label{construct:surf-from-beta-pmc}
  A decorated pointed matched circle~$\PMC^\epsilon$ (where $\epsilon$
  is $\alpha$ or $\beta$) gives rise to a surface $F$ as follows.
  Consider the disk $D_0$ with boundary $Z$. We orient $D_0$
  so that the specified orientation of $Z_\alpha$ agrees with its induced orientation
  from $\partial D_0$. Add a one-handle $s$
  along the pair of points $\{h,t\}$; then add
  one-handles along all the $M$-matched points in $\mathbf{p}$; and finally attach
  two two-handles to fill the two remaining boundary components.
  Call the resulting surface $F(\PMC^\epsilon)$.
  Each of the two disks attached at the last step meets exactly one of
  $Z_\alpha$ or $Z_\beta$: we call these disks $D_\alpha$ and
  $D_\beta$, respectively.  This surface has a preferred embedded disk
  $D_\alpha\cup s\cup D_\beta$, which is decomposed into three pieces.
  The surface inherits an orientation from $D_0$. (See
  Figure~\ref{fig:beta-pmc-to-surf} for an illustration.)
\end{construction}

The data of a decorated pointed matched circle $\PMC^{\alpha}$ is
equivalent to the earlier data of a pointed matched circle: we
contract the interval $Z^{\beta}$ to give the basepoint. The
underlying surfaces can also be identified.

\begin{definition}
  \label{def:MCGoid}
  The \emph{strongly based $\alpha\Hyph\beta$ mapping class groupoid} is the groupoid
  whose objects
  are decorated pointed matched circles.
  The morphisms from $\PMC_1^{\epsilon_1}$ to $\PMC_2^{\epsilon_2}$
  are isotopy classes of homeomorphisms 
  $$\phi\co F(\PMC_1^{\epsilon_1})\to F(\PMC_2^{\epsilon_2})$$
  so that
  \[
    \phi(D_\alpha(\PMC_1^{\epsilon_1}))=D_\alpha(\PMC_2^{\epsilon_2}) \qquad
    \phi(D_\beta(\PMC_1^{\epsilon_1}))=D_\beta(\PMC_2^{\epsilon_2}) \qquad
    \phi(s(\PMC_1))=s(\PMC_2).
    \]
\end{definition}

In this paper, we sometimes consider homeomorphisms which are not in
the mapping class groupoid as above. Specifically, we will sometimes be
interested in homeomorphisms $\phi\co
F(\PMC^\epsilon_1)\to F(\PMC^{\epsilon}_2)$ which preserve
the preferred disk $D_\alpha\cup s\cup D_\beta$, but which switch
$D_\alpha$ and $D_\beta$. 
It is equivalent to consider the surface $F^\circ=F\setminus
\interior(D_\alpha\cup s\cup D_\beta)$, and think of the induced map
on $F^{\circ}$, which exchanges $Z_\alpha=(\partial F^\circ) \cap
D_\alpha$ and $Z_\beta=(\partial F^\circ) \cap D_\beta$.

One relevant class of such homeomorphisms is the following:

\begin{definition}
  \label{def:Twins}
  Let $\PMC^\alpha$ and $\PMC^\beta$ be two pointed matched circles
  which differ only in the $\alpha$- or $\beta$-labels on the two
  intervals in the decomposition of $Z$. That is, $Z_\alpha$ for
  $\PMC^\alpha$ coincides with $Z_\beta$ for $\PMC^\beta$.  In this
  case, we say that $\PMC^\alpha$ and $\PMC^\beta$ are {\em twin
    pointed matched circles}.  For twin pointed matched circles, there
  are canonical orientation-reversing homeomorphisms
  \[
K_{\alpha,\beta}\co F(\PMC^\alpha)\to F(\PMC^\beta)
\qquad 
K_{\beta,\alpha}\co F(\PMC^\beta)\to F(\PMC^\alpha).
\]
Note that although these homeomorphisms send preferred disks to
preferred disks, they do not preserve the decorations on those
disks: $K_{\alpha,\beta}$ and $K_{\beta,\alpha}$ map $D_\alpha$ to $D_\beta$ and
vice versa.
\end{definition}

Another relevant class of such homeomorphisms are the ``half boundary
Dehn twists'' defined below.
Before defining them, we introduce some more terminology.

\begin{definition}
  Let $A$ be an (oriented) annulus with one boundary component marked
  as the ``inside boundary'' and the other as the ``outside
  boundary.'' A {\em radial curve} is any embedded curve in
  $A$ which connects the inside and outside boundary of $A$. Suppose that
  $r$ and $r'$ are two oriented, radial curves which intersect the
  inside boundary of $A$ at the same
  point, but which are otherwise disjoint. We say that $r'$ is
  {\em to the right} of $r$ if $r$ has a regular neighborhood $U$
  with an orientation-preserving identification with $(-\epsilon,\epsilon)\times [0,1]$, so that $r$ is
  identified with $\{0\}\times [0,1]$, the inside boundary meets $U$ in
  $(-\epsilon,\epsilon)\times\{0\}$, and $r'\cap U$ is contained in
  $[0,\epsilon)\times [0,1]$.
\end{definition}

\begin{definition}
  \label{def:HalfDehn}
  Let $(F,D_\alpha\cup s \cup D_\beta)$ be a surface with a preferred disk
  decomposed into three parts. A {\em positive half Dehn twist along
    the boundary}, denoted $\tau_{\bdy}^{1/2}$, is a homeomorphism with the
  following properties:
  \begin{itemize}
  \item there is a disk neighborhood $N$ of $D_\alpha\cup s \cup D_\beta$
    so that
    $\tau_{\bdy}^{1/2}$  fixes the complement of~$N$;
  \item $\tau_{\bdy}^{1/2}$ maps the preferred disk to itself,
    but switches $D_\alpha$ and $D_\beta$; and
  \item there is a radial arc $r$ in the annulus $A=N\setminus \interior(D_\alpha\cup s \cup D_\beta)$
    (oriented so that it terminates at $\partial (D_\alpha\cup s \cup
    D_\beta)$) which
    is mapped under $\tau_\bdy^{1/2}$ to 
    a new  arc $r'$, which is to the right of $r$. Here, we view
    $\bdy(D_\alpha\cup s\cup D_\beta)$ as the outside boundary of $A$.
  \end{itemize}
  \begin{figure}
    \centering
    \input{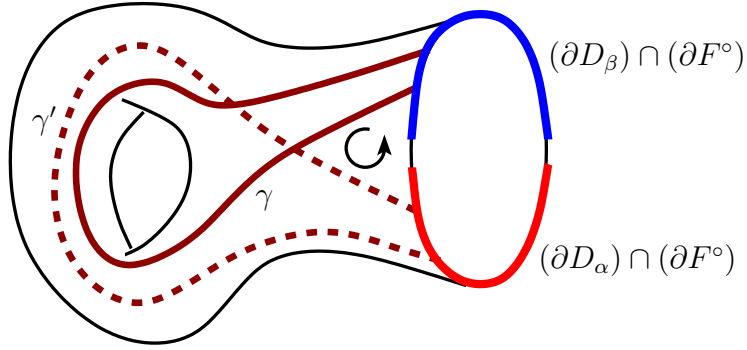}
    \caption{\textbf{An example of a
        positive half Dehn twist along the boundary.}
    \label{fig:HalfDehnTwist}
    The boundary of the genus one surface
    $\PunctF=F\setminus(D_\alpha\cup s\cup D_\beta)$ has a pair of
    distinguished arcs $(\bdy D_\alpha\cap\bdy \PunctF) \subset
    Z_\alpha$ and $(\bdy D_\alpha\cap\bdy \PunctF)\subset Z_\beta$. We
    have illustrated an arc $\gamma$ with boundary on $\bdy D_\beta$
    whose image is (up to isotopy relative to the boundary) the curve
    $\gamma'$ represented by the dashed line.}
  \end{figure}
  See Figure~\ref{fig:HalfDehnTwist} for an illustration.
  A {\em negative half Dehn twist along the boundary},
  denoted $\tau_{\bdy}^{-1/2}$, is the inverse to a positive half Dehn
  twist along the boundary. For a surface $F$ with one boundary
  component, the isotopy class of the surface homeomorphism specified
  by a half (positive or negative) Dehn twist (among homeomorphisms
  preserving the division of the disk) is uniquely specified by the
  above properties.  Therefore we sometimes refer to ``the'' (rather
  than ``a'') positive (respectively negative) half Dehn twist along
  the boundary.

  Similarly, a full Dehn twist along the boundary, $\tau_\bdy$, is the
  composite of two positive half Dehn twists.
\end{definition}

\begin{lemma}
  \label{lem:CharHalfDehnTwist}
  Let $f\co F \to F$ be a homeomorphism preserving $D_\alpha\cup s\cup
  D_\beta$, and let $\PunctF=F\setminus(D_\alpha\cup s\cup
  D_\beta)$. Suppose that the following conditions hold.
  \begin{enumerate}
    \item 
      \label{item:IdMap}
      The induced automorphism ${\overline f}$ of
      the closed surface $F^\circ/\bdy F^\circ$
      is isotopic, relative to the basepoint $[\bdy F^\circ]$, to the identity map.
    \item\label{item:ExchangeDaDb}
      The map $f$ exchanges $D_\alpha$ and $D_\beta$.
    \item 
      \label{item:HalfTwist}
      There is a curve $\gamma\subset \PunctF$ and a curve $\gamma'$ isotopic
      relative to endpoints to $f(\gamma)$, so that:
      \begin{enumerate}
      \item either the boundary of $\gamma$ is contained in $\bdy
        D_\beta\cap\bdy\PunctF$ or the boundary of $\gamma$ is
        contained in $\bdy D_\alpha\cap\bdy\PunctF$,
      \item $[\gamma]\neq 0\in H_1(\PunctF,\bdy\PunctF)$,
      \item\label{item:TwoPoints} $\gamma$ intersects $\gamma'$ transversely in exactly two points,
      \item $\gamma\cup\gamma'$ has one component $T$ which is a disk, and
      \item in the cyclic order induced by the orientation on
        $\PunctF$, the boundary of $T$ consists of an arc in
        $\gamma$, an arc in $\gamma'$, and an arc in $\partial
        F^\circ$.
      \end{enumerate}
  \end{enumerate}
  (See Figure~\ref{fig:HalfDehnTwist}.)
  Then, $f$ is isotopic (as maps preserving $D_\alpha\cup s \cup
  D_\beta$ but switching $D_\alpha$ and $D_\beta$) to a positive half
  Dehn twist along the boundary.  
\end{lemma}

\begin{proof}
  Recall that the mapping class group of a surface $F$ fixing a disk
  $D\subset F$ is a $\ZZ$-central extension of the mapping class group
  of $F$ preserving $D$ set-wise (but not point-wise); and the Dehn
  twist $\tau_\bdy$ along $\bdy(F\setminus D)$ is a generator for this
  distinguished central $\ZZ$. Thus, Conditions~(\ref{item:IdMap})
  and~(\ref{item:ExchangeDaDb}) imply that $f$ is isotopic to 
  $\tau_\partial^{m+1/2}$, for some $m\in\ZZ$. If $m\not\in\{-1,0\}$
  then the minimal number of intersection points between $\gamma$ and
  $f(\gamma)$ is greater than two, contradicting
  Condition~\eqref{item:TwoPoints}.  The condition on the boundary of
  the triangle~$T$
  ensures that, in fact, $m$ is $0$, not $-1$.
\end{proof}

\subsection{\textalt{$\beta$}{beta}-bordered Heegaard diagrams}\label{sec:beta-bordered}
In previous work~\cite{LOT1,LOT2}, the $\alpha$-curves consisted of arcs
and circles, while the $\beta$-curves were always circles. To prove
the $\Hom$ pairing theorem, we will also want to work with diagrams
where $\beta$'s, instead of $\alpha$'s, go out to the boundary.
\begin{definition}
  A \emph{$\beta$-bordered Heegaard diagram} is a quadruple
  $\HD^\beta=(\overline{\Sigma},\alphas,\overline{\betas},{\mathbf z})$ where:
  \begin{itemize}
  \item $\overline{\Sigma}$ is a compact surface of genus $g$ with one
    boundary component;
  \item $\alphas$ is a $g$-tuple of pairwise disjoint
    circles in the interior $\Sigma$ of $\overline{\Sigma}$;
  \item $\overline{\betas}$ is
    \[
    \overline{\betas}=\{\overbrace{\overline{\beta}_1^a,\dots,\overline{\beta}_{2k}^a}^{\overline{\betas}^a},\overbrace{\beta_1^c,\dots,\beta_{g-k}^c}^{\betas^c}\},
    \]
    a collection of pairwise disjoint embedded arcs (the
    $\overline{\beta}_i^a$) with boundary on $\bdy\overline{\Sigma}$ and circles
    (the $\beta_i^c$) in the interior $\Sigma$ of $\overline{\Sigma}$; and
  \item ${\mathbf z}$ is an arc in
    $\bdy\overline{\Sigma}\setminus\overline{\betas}^a$.
  \end{itemize}
We require that
$\Sigma\setminus\betas$ and $\Sigma\setminus\alphas$ both be
connected; this translates to the condition that the $\beta$-
(respectively $\alpha$-) curves be linearly independent in
$H_1(\overline{\Sigma},\bdy\overline{\Sigma})$.

The diagram ${\mathcal H}^\beta$ gives a natural $\beta$-pointed matched circle in 
the following way. Take two copies of $\partial \overline \Sigma\setminus
{\bf z}$, call them $Z_\alpha$ and $Z_\beta$, and orient them both
with the
orientation they inherit from
$\partial{\overline\Sigma}$. Let $Z$ be the result of gluing
$Z_\alpha$ to $Z_\beta$ head-to-head and tail-to-tail. Let
${\mathbf p}=\overline{\betas}^a\cap\bdy\overline{\Sigma}$, thought of
as a subset of
$Z_\beta$. Let $M$ be the involution exchanging
the endpoints of each $\overline{\beta}_i^a$.  We call this
pointed matched circle $\PMC(\HD^\beta)$.
\end{definition}
We will sometimes call bordered Heegaard diagrams as defined
in~\cite{LOT1} \emph{$\alpha$-bordered Heegaard
  diagrams}, and denote them $\HD^\alpha$, to
distinguish them from $\beta$-bordered Heegaard diagrams.  Note that
we have a slight shift in point of view from~\cite{LOT1}: we now think
of ${\mathbf z}$ as an interval, rather than
just a point $z$, and our circle $Z$ is no longer $\partial\overline{\Sigma}$, but
rather two copies of an interval in $\overline{\Sigma}$.
An $\alpha$-bordered Heegaard diagram specifies an
$\alpha$-pointed matched circle in exactly the same way as a
$\beta$-bordered Heegaard diagram specifies a $\beta$-pointed
matched circle, but placing
the points in $Z_\alpha$ rather than $Z_\beta$.

\begin{figure}
  \centering
  \includegraphics{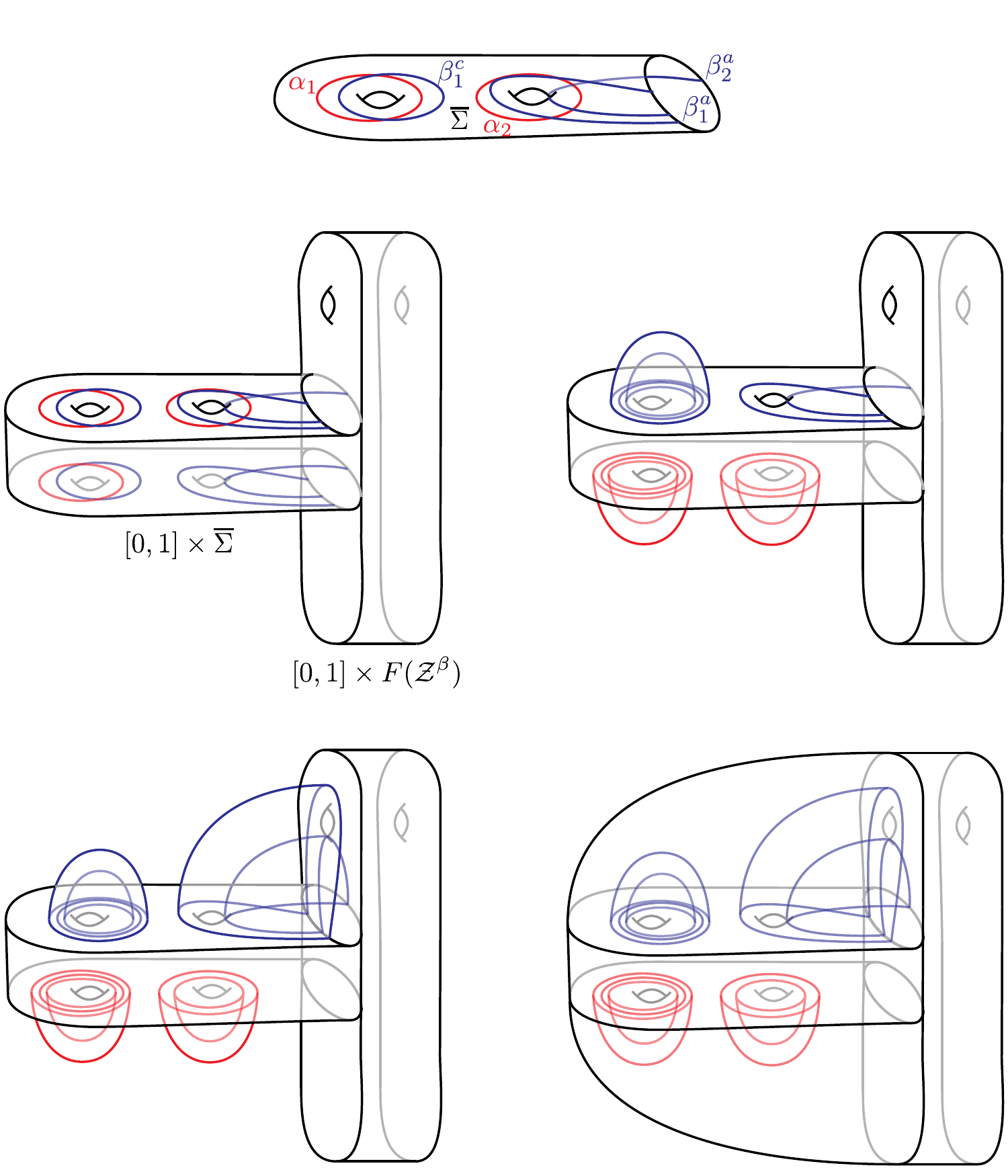}
  \caption{\textbf{Building a bordered 3-manifold from a
      $\beta$-bordered Heegaard diagram.} The analogous figure in the
    $\alpha$-bordered case is~\cite[Figure~\ref*{LOT1:fig:bord-to-3mfld}]{LOT1}.}
  \label{fig:BuildYbeta}
\end{figure}

\begin{construction}\label{construct:beta-bord-to-Y}
 Let $\PMC^\beta$ denote the pointed matched
 circle specified by $\HD^\beta$. 
 There is an associated bordered three-manifold $Y(\HD^\beta)$,
 constructed in the following four steps.
\begin{enumerate}
  \item 
 Glue $[0,1]\times\overline{\Sigma}$ to
 $[0,1]\times F(\PMC^\beta)$, by identifying 
 $[0,1]\times(\partial
     \overline{\Sigma})$ with $\{0\}\times (D_0 \cup s)$ so that
 \begin{itemize}
   \item  $({\mathbf z}\subset
     \bdy\overline{\Sigma})\times[0,1]$ is identified with the one-handle $\{0\}\times (s\subset
     F(\PMC^\beta))$
   \item $[0,1]\times (\partial\overline{\Sigma}\setminus
     \mathbf{z})$ is identified with $\{0\}\times D_0$
     \item 
 $\{1\}\times {\mathbf p} \subset [0,1]\times \overline{\Sigma}$ is
 identified with $\{0\}\times {\mathbf p} \subset [0,1]\times F(\PMC^\beta)$.
 \end{itemize}
 (The result is naturally a manifold
 with corners.) 
 \item Attach $3$-dimensional $2$-handles to
 the $\{0\}\times \alpha_i\subset [0,1]\times \overline{\Sigma}$ and to
 the $\{1\}\times \beta_i^c\subset[0,1]\times \overline{\Sigma}$. 
 \item Let
 $\eta_i$ denote the core of the $1$-handle in $F(\PMC^\beta)$ attached
 along $\beta_i^a\cap\bdy\overline{\Sigma}$. Then $(\{1\}\times
 \beta_i^a)\cup (\{0\}\times \eta_i)$ is a closed circle; attach thickened
 disks (3-dimensional 2-handles) along these circles.

 The result of the attaching so far is a $3$-manifold with three boundary
 components: two copies of $S^2$ 
 (one containing $D_\alpha$ and the other containing $D_\beta$)
 and a copy of $F(\PMC^\beta)$. 
\item Fill in the
 two $S^2$ boundary components with $3$-balls.
 \end{enumerate}
 See Figure~\ref{fig:BuildYbeta}.

 The boundary of
 $Y(\HD^\beta)$ is naturally (orientation-preserving) identified with
 $F(\PMC^\beta)$.
\end{construction}

\begin{remark}\label{rmk:why-this-construction}
  The construction of $Y(\HD^\beta)$ is convenient for describing the
  effect of gluing Heegaard diagrams.  For the $\alpha$-bordered case
  (which is exactly analogous), a shorter description of an equivalent
  bordered three-manifold is given
  in~\cite[Construction~\ref*{LOT2:construct:bordered-HD-to-bordered-Y-one-bdy}]{LOT2};
  see also
  Construction~\ref{constr:AlphaBetaBorderedDiagramToThreeManifold}
  below. In particular, this construction identifies a regular
  neighborhood of the $\alpha$-arcs union the boundary in $\HD^\alpha$
  with $\PunctF(\bdy\HD)$, in an orientation-reversing way. By
  contrast, for a $\beta$-bordered Heegaard diagram, the corresponding
  identification is orientation-preserving.
\end{remark}

We now wish to define $\CFDa(\HD^\beta)$ and $\CFAa(\HD^\beta)$.  As
for $\alpha$-bordered diagrams,
these will be defined by counting moduli spaces
$\cM^B(\x,\y; \vec{\rhos})$ between two generators $\x$ and $\y$
asymptotic to certain sequences of sets of Reeb chords on
$\PMC^\beta$.  (See \cite[Definition~\ref*{LOT1:def:emb-ind-emb-chi}]{LOT1}.)  However, the fact the diagram is $\beta$-bordered leads
to some reversals.  It
is easiest to see what happens by reference to a corresponding
$\alpha$-bordered diagram.

\begin{definition}\label{def:dual-HD-1}
  Given an $\alpha$-bordered Heegaard diagram
  $\HD^\alpha=(\Sigma,\alphas,\betas,\mathbf{z})$ there is an associated
  $\beta$-bordered Heegaard diagram
  $\overline{\HD}{}^\beta=(\Sigma,\alphas^\beta,\betas^\beta,\mathbf{z})$ obtained by
  setting $\beta_i^{\beta,c}=\alpha_i^c$,
  $\beta_i^{\beta,a}=\alpha_i^a$, $\alpha_i^\beta=\beta_i$.
\end{definition}

\begin{lemma}\label{lem:beta-bord}
  Let $Y$ be a $3$-manifold and $\phi\co F(\PMC^\alpha)\to \bdy Y$ a
  parameterization of its boundary. 
  Let $\HD^\alpha$ be an $\alpha$-bordered Heegaard diagram for
  $(Y,\phi)$. Then $\overline{\HD}{}^\beta$ is a $\beta$-bordered Heegaard diagram
  for $(-Y,\phi\circ K_{\beta,\alpha}\co {F(\PMC^\beta)}\to -\bdy Y)$.
\end{lemma}
\begin{proof}
  Recall
  from~\cite[Construction~\ref*{LOT1:constr:bordered-hd-gives-mfld}]{LOT1}
  that if $\HD^\alpha$ is an $\alpha$-bordered Heegaard diagram then
  to construct $Y(\HD^\alpha)$ one thickens $\overline{\Sigma}$ to
  $\overline{\Sigma}\times[0,1]$; glues the boundary
  $(\bdy\overline{\Sigma})\times [0,1]$ to 
  $(D_0\cup s)\subset F(\PMC^\alpha)$, and then one glues thickened disks
  to the following objects:
  \begin{itemize}
  \item  the $\alpha$-circles in $\overline{\Sigma}\times\{0\}$,
  \item  $\beta$-circles in $\overline{\Sigma}\times\{1\}$, and
  \item the unions
    of the $\alpha$-arcs in $\overline{\Sigma}\times\{0\}$ and the cores
    of the $1$-handles of $F(\PMC^\alpha)$.
  \end{itemize}
  Finally, one caps off the two $S^2$
  boundary components with $3$-balls.

  This process results in a manifold $Y(\HD^\alpha)$ which is the mirror
  image, in an obvious sense, of the manifold $Y(\overline{\HD}{}^\beta)$ from
  Construction~\ref{construct:beta-bord-to-Y}; reflecting across
  $\Sigma\times\{1/2\}$ gives an orientation-reversing homeomorphism
  between the two bordered manifolds.
\end{proof}

For each generator $\x$ of $\HD^\alpha$, there is an obvious corresponding
generator $\ol{\x}$ of~$\ol{\HD}^\beta$.  Similarly, for a homology class
$B\in\pi_2(\x,\y)$ let $\ol{B}\in\pi_2(\ol{\y},\ol{\x})$ denote the
homology class with the same local multiplicities in $\Sigma$ as $B$.

\begin{lemma}\label{lem:mod-beta-rev}
  For $\HD$ an $\alpha$-bordered Heegaard diagram, $\x, \y \in
  \Gen(\HD)$, $B \in \pi_2(\x,\y)$, and $\vec\rhos$ any sequence of
  sets of Reeb chords, there is a homeomorphism
  \[
  \cM^B(\x,\y;\vec{\rhos})\cong \cM^{\ol{B}}(\ol{\y},\ol{\x},\vec{\rhos}^\op)
  \]
  where $\vec{\rhos}^\op$ is $\vec{\rhos}$ read in the opposite order.
\end{lemma}

\begin{proof}
  Both moduli spaces are defined by counting pseudoholomorphic curves
  in $\Sigma \times [0,1] \times \RR$.  Reflecting in both the $[0,1]$
  and $\RR$ directions gives an identification between the two moduli
  spaces.
\end{proof}

For a $\beta$-bordered Heegaard diagram $\HD^\beta$, we can therefore use
all the techniques of~\cite{LOT1} to define $\CFDa(\HD^\beta)$ and
$\CFAa(\HD^\beta)$, except that the order of Reeb chords on the boundary is
reversed.  We can achieve this algebraically either by reversing all
the chords or (as we prefer) by viewing our modules as defined over
the opposite algebra.  (These are equivalent, by
Equation~\eqref{eq:OpAlg}.)  Thus we view $\CFDa(\HD^\beta)$ as a (left)
type $D$ structure over $\Alg(-\PMC)^\op$, i.e., a right type $D$
structure over $\Alg(-\PMC)$, and $\CFAa(\HD^\beta)$ as a left
$\Ainf$-module over $\Alg(\PMC)$.

We digress briefly to discuss the identification between
$\SpinC$-structures on~$Y$ and $-Y$. As usual, we view a
$\SpinC$-structure on~$Y$ as a homology class of non-vanishing vector
fields. Given a $\SpinC$-structure $\spinc$ on~$Y$, induced by a
vector field~$v$, let $-\spinc$ be the $\SpinC$-structure on~$-Y$
induced by the vector field~$-v$. (Note that the oriented $2$-plane
fields $v^\perp$ on $Y$ and $(-v)^\perp$ on $-Y$ are the same.) To
avoid confusion, recall that the conjugate $\SpinC$-structure
$\overline{\spinc}$ is also represented by $-v$, but viewed as a
vector field on~$Y$.  Thus, $-\overline{\spinc}$ is represented by $v$
as a vector field on~$-Y$.

\begin{proposition}\label{prop:beta-is-dual}
  Let $\HD^\alpha$ be an $\alpha$-bordered Heegaard diagram with boundary
  $\PMC=\PMC^\alpha$. Then $\CFDa(\ol{\HD}^\beta)$ and
  $\CFAa(\ol{\HD}^\beta)$ are duals (in the senses of
  Definitions \ref{def:dual-type-D} and~\ref{def:dual-type-A}) of the
  corresponding structures for $\HD^\alpha$:
  \begin{align*}
    \CFDa(\ol{\HD}^\beta,-\spinc)^{\Alg(-\PMC)} &\cong \lsup{\Alg(-\PMC)}\overline{\CFDa(\HD^\alpha,\spinc)}\\
    \lsub{\Alg(\PMC)}\CFAa(\ol{\HD}^\beta,-\spinc) &\cong \overline{\CFAa(\HD^\alpha,\spinc)}_{\Alg(\PMC)}.
  \end{align*}
\end{proposition}
\begin{proof}
  We will prove the duality result for $\CFDa$; the proof for $\CFAa$
  is analogous.  Let
  $X(\HD^\alpha)$ (respectively $X(\ol{\HD}^\beta)$) denote the $\Ground$-module
  generated by $\Gen(\HD^\alpha)$ (respectively $\Gen(\ol{\HD}^\beta)$).  We have
  an isomorphism $X(\ol{\HD}^\beta)=\Hom(X(\HD^\alpha),\Ground)$ by
  setting, for generators $\x,\y\in \Gen(\HD^\alpha)$,
  \[
  \ol{\x}(\y)=
  \begin{cases}
    \iota_{\x} & \text{if }\x=\y\\
    0 &\text{otherwise,}
  \end{cases}
  \]
  where $\iota_{\x}\in\Ground$ is the primitive idempotent
  corresponding to $\x$, so $\iota_\x\x=\x$.

  The type $D$ structure on $\CFDa(\HD^\alpha)$ is given by
  \begin{align*}
    \delta^1&\co X(\HD^\alpha)\to \Alg(-\PMC)\otimes X(\HD^\alpha)\\
    \delta^1(\x)&=\sum_{\y}\sum_{\substack{B\in\pi_2(\x,\y)\\\vec{\rho}\text{
          compatible with }B\\
        \ind(B,\vec{\rho})=1}}a(-\vec{\rho})\otimes \y \cdot \#(\Mod^B(\x, \y; \vec{\rho})).
  \end{align*}
  The operation $\delta^1_\beta\co X(\ol{\HD}^\beta)\to
  X(\ol{\HD}^\beta)\otimes\Alg(-\PMC)$ is defined similarly, but using the
  moduli spaces on~$\ol{\HD}^\beta$.  By Lemma~\ref{lem:mod-beta-rev},
  terms of the form $a(-\vec{\rho})\otimes \y$ in $\delta^1(\x)$ correspond
  to terms of the form $\ol{\x} \otimes a(-\vec\rho)$ in
  $\delta^1_\beta(\ol{\y})$.  (We are considering the moduli space in
  which the chords appear in reverse order, but we also multiply the
  algebra elements in the reverse order, so it is again $a(-\vec\rho)$
  that is relevant.)
  This is exactly the statement that $\CFDa(\ol{\HD}^\beta)$ is the dual of
  $\CFDa(\HD^\alpha)$.

  The behavior on the $\SpinC$ structures comes from the observation
  that, when $\x$ is viewed as
  a generator for the bordered Floer homology of $\HD$, its
  corresponding vector field points in the opposite direction from
  that of $\ol{\x}$ when viewed as a generator for the bordered Floer homology
  of $\ol{\HD}$.
\end{proof}

Note that $\ol{\HD}^\beta$ represents $-Y(\HD^\alpha)$.  There is
another way of creating a Heegaard diagram for $-Y(\HD^\alpha)$,
namely by considering $-\HD^\alpha$, which is the same Heegaard
diagram but with the orientation on the underlying surface~$\Sigma$
reversed. This operation also has the effect of dualizing modules:

\begin{proposition}\label{prop:minus-hd-mod}
  Suppose that $\HD$ is an $\alpha$-bordered Heegaard
  diagram with boundary~$\PMC$. Let $-\HD$
  denote the same Heegaard diagram but with the orientation
  of $\Sigma$ reversed. Then $\CFDa(-\HD)$ is a left
  $\Alg(\PMC)$-module. If we view
  $\CFDa(-\HD)$ as a right
  $\Alg(-\PMC)$-module, then
  $\CFDa(-\HD,-\overline{\spinc})$ is dual to $\CFDa(\HD,\spinc)$. 
  Similarly, $\CFAa(-\HD,-\overline{\spinc})$ is dual to $\CFAa(\HD,\spinc)$.
\end{proposition}
\begin{proof}
  This follows from a similar argument to Proposition~\ref{prop:beta-is-dual}.
\end{proof}

\subsection{\textalt{$\alpha\Hyph\beta$}{alpha-beta}-bordered Heegaard diagrams}
\label{sec:abDiagrams}

We will need to generalize the notion of arced bordered Heegaard
diagrams with two components and the associated 3-manifolds to cases
where the boundaries meet
$\alpha$- or $\beta$-circles.
The following is an extension
of~\cite[Definition~\ref*{LOT2:def:StronglyBordered}]{LOT2} to our
more symmetric language.
\begin{definition}
A \emph{strongly bordered three-manifold with two boundary
components} is specified by the following data:
\begin{itemize}
\item A $3$-manifold $Y$ with two boundary components $\partial_LY$
  and $\partial_RY$.
\item Orientation-preserving homeomorphisms $\phi_L\co
  F(\PMC^{\epsilon_L}_L)\to \bdy_LY$ and $\phi_R\co
  F(\PMC^{\epsilon_R}_R)\to \bdy_RY$ for some pointed matched circles
  $\PMC^{\epsilon_L}_L$ and $\PMC^{\epsilon_R}_R$.
\item A tunnel $[0,1]\times D$ connecting the two boundary
  components, which is  divided into three balls
  $$
  ([0,1]\times D_\alpha) \cup ([0,1]\times s) \cup ([0,1]\times D_\beta)
  $$
  where $\{0\}\times (D_\alpha\cup s\cup D_\beta)$ coincides with the
  corresponding part of $\partial_L Y=\phi_L(F(\PMC^{\epsilon_L}_L))$ and $\{1\}\times (D_\alpha\cup s\cup
  D_\beta)$ coincides with the corresponding part of
  $\partial_R Y=\phi_R(F(\PMC^{\epsilon_R}_R))$.
\end{itemize}

\end{definition}

For $3$-manifolds with two boundary components, one can consider
\emph{arced $(\alpha,\alpha)$-bordered Heegaard diagrams} (which are
the only kind we considered in~\cite{LOT2}), \emph{arced
  $(\alpha,\beta)$-bordered Heegaard diagrams}, \emph{arced
  $(\beta,\alpha)$-bordered Heegaard diagrams}, and \emph{arced
  $(\beta,\beta)$-bordered Heegaard diagrams}. We call all of these
types of diagrams, collectively, {\em arced bordered
  Heegaard diagrams}.  The definitions of the latter three types are
trivial adaptations of the definition of an arced bordered Heegaard
diagram from~\cite[Section~\ref*{LOT2:sec:Diagrams}]{LOT2}; except
that in keeping with our
present conventions (where we thicken basepoints), rather than drawing
simply an arc connecting the two basepoints on the boundary, now
${\mathbf z}$ denotes a rectangle in the Heegaard surface which
connects the arc of basepoints on one boundary component with the arc
of basepoints in the other.  For $(\beta,\beta)$-bordered Heegaard
diagrams the $\beta$-curves but not the $\alpha$-curves go out to the
boundaries. For $(\alpha,\beta)$-bordered diagrams,
the $\alpha$-curves go out to $\bdy_L\overline{\Sigma}$
and the $\beta$-curves go out to $\bdy_R\overline{\Sigma}$.
For $(\beta,\alpha)$-bordered Heegaard diagrams, the
$\beta$-curves go out to $\bdy_L\overline{\Sigma}$ and the
$\alpha$-curves go out to $\bdy_R\overline{\Sigma}$.

More explicitly, if $\HD$ is an arced bordered Heegaard diagram,
then there are two intervals $Z_L$ and $Z_R$ in $\partial{\overline \Sigma}$,
with the property that 
${\mathbf p}= (\partial{\overline\Sigma})\cap ({\overline \alphas}^a\cup{\overline\betas}^a)$ is contained in $Z_L\cup Z_R$, so that all points
${\mathbf p}\cap Z_L$ are all of the same type (i.e., all are either 
boundary points of $\alpha$-arcs or $\beta$-arcs) and 
${\mathbf p}\cap Z_R$ are also of the same type. Correspondingly
we let $\partial_L\HD$ denote the pointed matched circle gotten by doubling
$Z_L$ and marking it with $\alpha$ or $\beta$, according to the type of 
points in ${\mathbf p}\cap Z_L$, and we obtain $\partial_R \HD$ analogously.

An arced bordered diagram gives rise to a strongly bordered
three-manifold via the following generalization of
Construction~\ref{construct:beta-bord-to-Y}.

\begin{construction}
  \label{constr:BorderedDiagramToThreeManifold}
  Let $\PMC_L$
  and $\PMC_R$ denote the pointed matched circles specified by
  $\partial_L\HD$ and $\partial_R\HD$.
  There is an associated strongly bordered three-manifold 
  $Y(\HD)$ with two boundary components constructed in the following four steps:
  \begin{enumerate}
  \item 
  Glue $[0,1]\times \overline{\Sigma}$ to 
  $[0,1] \times F(\PMC_L)$, by gluing
  $[0,1]\times \partial_L \bSigma$ to $\{0\}\times F(\PMC_L)$ and
  $[0,1]\times \partial_R \bSigma$ to $\{0\}\times F(\PMC_R)$
  following Construction~\ref{construct:beta-bord-to-Y}.
  \item 
  Next, glue $2$-handles along
  $\{0\}\times \alpha_i^c\subset \overline[0,1]\times {\Sigma}$ and to
  $\{1\}\times \beta_i^c\subset[0,1]\times \overline{\Sigma}$.  
  \item
    \label{step:GlueAlongArcs}
  Consider the arcs $\alpha_i^a$ and $\beta_i^a$, thought of as
  supported in $\{0\}\times {\overline\Sigma}$ and
  $\{1\}\times {\overline\Sigma}$ respectively. These are completed
  into closed circles by following the endpoints through the
  one-handles in $F(\PMC_L)$ or $F(\PMC_R)$ (wherever they go).
  Add two-handles along these circles.

  We end up with a
  three-manifold with four boundary components: $F(\PMC_L)$,
  $F(\PMC_R)$ and a pair of two-spheres. 
  \item Fill in the two-spheres
  with three-balls $B_\alpha$ and $B_\beta$, to obtain the desired
  three-manifold $Y(\HD)$.
  \end{enumerate}

  Morally, the tunnel $[0,1]\times D$ needed to make $Y(\HD)$ into a
  strongly bordered 3-manifold is given by $B_\alpha\cup [0,1]\times {\mathbf z}\cup B_\beta$. In fact, since we have attached copies of $[0,1]\times F(\PMC_L)$
  and $[0,1]\times F(\PMC_R)$ at the boundary, this tunnel is actually given by
  \begin{align*}
    \big(([0,1]\times D_\alpha)\cup B_\alpha\cup ([0,1]\times D_\alpha)\big)
    &\cup \big(([0,1]\times s)\cup ([0,1] \times {\mathbf z})\cup ([0,1]\times s)\big) \\
  &\cup \big(([0,1]\times D_\beta)\cup B_\beta\cup ([0,1]\times D_\beta)\big).
  \end{align*}
\end{construction}

\begin{lemma}
  If $\HD$ and $\HD'$ are two arced bordered Heegaard diagrams with
  $\bdy_R \HD = -\bdy_L \HD' = \PMC$ (agreeing as pointed matched circles,
  including decoration), then
  \[
  Y(\HD \sos{\bdy_R}{\cup}{\bdy_L} \HD') = Y(\HD) \cup_{F(\PMC)} Y(\HD').
  \]
\end{lemma}

\begin{proof}
  This is immediate from the definitions.
\end{proof}

We give the following simpler description in the $\alpha$-$\beta$-bordered
case which will be used in the proof of Proposition~\ref{prop:Denis-reps}:
\begin{construction}
  \label{constr:AlphaBetaBorderedDiagramToThreeManifold}
  An arced ($\alpha$,$\beta$)-bordered Heegaard diagram $\HD$ specifies a strongly
  bordered three-manifold as follows. Attach disks $D_L$ and $D_R$ to
  $\partial_L{\overline\Sigma}$ and $\partial_R{\overline\Sigma}$
  respectively obtain a closed surface $\Sigma$. Let $Y$ be
  $[0,1]\times \Sigma$ with three-dimensional two-handles attached
  along the $\{0\}\times \alpha_i^c$ and $\{1\}\times
  \beta_i^c$.

  $F(\PMC_L)$ is identified with $\partial_L Y$ as follows:
  \begin{itemize}
  \item Identify $s$ 
    with $\{0\}\times {\mathbf z}$.
  \item Identify $D_\beta$ with $\{0\}\times D_L$.
  \item Identify $D_\alpha$ with $\{0\}\times D_R$.
  \item Identify $\PunctF(\PMC_L)$
    with the complement of $\{0\}\times D_L\cup
    \{0\}\times{\mathbf z}\cup \{0\}\times D_R\subset \partial_L
    Y$ using the $\alpha$-arcs.
  \end{itemize}

  $F(\PMC_R)$ is identified with $\partial_R Y$ as follows:
  \begin{itemize}
  \item Identify $s$ 
    with $\{1\}\times {\mathbf z}$. 
   \item Identify $D_\beta$ with $\{1\}\times D_L$.
  \item Identify $D_\alpha$ with $\{1\}\times D_R$
  \item Identify $\PunctF(\PMC_R)$
    with the complement of $\{1\}\times D_L\cup
    \{1\}\times{\mathbf z}\cup \{1\}\times D_R\subset \partial_L
    Y$ using the $\beta$-arcs.
  \end{itemize}
  
  The tunnel is, of course, $[0,1]\times (D_L\cup{\mathbf z}\cup D_R)$. 
\end{construction}

\begin{lemma}
  If $\HD$ is an $\alpha$-$\beta$-bordered diagram, the
  strongly bordered three-manifolds specified in Constructions~\ref{constr:AlphaBetaBorderedDiagramToThreeManifold} and~\ref{constr:BorderedDiagramToThreeManifold} are canonically isomorphic.
\end{lemma}

\begin{proof}
  Let $Y_1(\HD)$ (respectively $Y_2(\HD)$) be the strongly bordered
  $3$-manifold given by
  Construction~\ref{constr:BorderedDiagramToThreeManifold}
  (respectively
  Construction~\ref{constr:AlphaBetaBorderedDiagramToThreeManifold}). Observe
  that $Y_2(\HD)$ is a subspace of $Y_1(\HD)$ in an obvious way.  The
  two-handles attached in Step~(\ref{step:GlueAlongArcs}) of
  Construction~\ref{constr:BorderedDiagramToThreeManifold} specify a
  deformation retraction of $Y_1(\HD)$ to $Y_2(\HD)$ (by folding the
  boundary along the $2$-handles to the Heegaard surface), respecting
  the strong bordering.
\end{proof}

\begin{construction}
  \label{constr:AssocHomeo}
  Let $Y$ be a strongly bordered $3$-manifold with boundary components parameterized by $F(\PMC_L)$ and $F(\PMC_R)$.
  Suppose that $Y$ is homeomorphic to the
  product of an interval with a surface.  Then we can define a map
  $$\phi_Y\co {-F(\PMC_L)} \to F(\PMC_R)$$
  in the strongly based mapping class groupoid
  as follows.
  
  First, fix a homeomorphism
  $\Phi\co [0,1] \times F(\PMC_R)\to Y$
  so that
  \begin{itemize}
  \item $\Phi|_{\{1\}\times F(\PMC_R)}=\phi_R$ and
  \item the images under $\Phi$ of $[0,1]$ times the three distinguished
    regions $D_\alpha$, $s$, and $D_\beta$ (in $F(\PMC_R)$)
    are mapped to the three corresponding
    distinguished regions in $Y$.
  \end{itemize}
  Then, let 
  $$
  \phi_{Y}= (\Phi|_{\{0\}\times F(\PMC_R)})^{-1}\circ (-\phi_L)\co
  F(\PMC_L)\rightarrow F(\PMC_R).
  $$
  We call the strongly bordered manifold $Y$ the \emph{mapping
    cylinder of $\phi$.}
\end{construction}

An adaptation of the argument
from~\cite[Lemma~\ref*{LOT2:lem:WellDefinedStronglyBordered}]{LOT2}
shows that the above construction gives a well-defined element of the
mapping class groupoid in the sense of Definition~\ref{def:MCGoid}.

\begin{definition}\label{def:mcg-action}
  Given a strongly based homeomorphism $\phi\co F(\PMC_1)\to
  F(\PMC_2)$ and a bordered $3$-manifold $(Y,\psi\co F(\PMC_1)\to \bdy
  Y)$ we can twist the parameterization of $Y$ by $\phi$ to give
  a new bordered $3$-manifold
  \[
  \phi(Y)=(Y,\psi\circ \phi^{-1}\co F(\PMC_2)\to \bdy Y).
  \]
  Equivalently, we can define $\phi(Y)$ by gluing the mapping cylinder
  of $\phi$ to $Y$:
  \[
  \phi(Y)=Y\cup_{F(\PMC_1)} M_\phi.
  \]
\end{definition}

For arced bordered Heegaard diagrams, one can define
bimodules as in~\cite{LOT2}. The case of $(\alpha,\alpha)$-bordered
diagrams is discussed there, and the story for
$(\beta,\beta)$-bordered diagrams is an entirely straightforward
adaptation (although with reversed algebras), analogous to the 
relation between $\alpha$-bordered and $\beta$-bordered
diagrams as discussed in
Section~\ref{sec:beta-bordered}.
Some aspects of the
$(\beta,\alpha)$-bordered case are new, however, and we discuss them now.

For a $(\beta,\alpha)$-bordered Heegaard diagram $\lsupv{\beta}\HD^\alpha$, a
\emph{generator} is a tuple $\x=\{x_i\}$ of intersection points
between $\alpha$- and $\beta$-curves, such that there is exactly one
$x_i$ on each $\alpha$- or $\beta$-circle, and no $\alpha$- or
$\beta$-curve contains more than one $x_i$. Note that, unlike the
case of an $(\alpha,\alpha)$- or $(\beta,\beta)$-bordered Heegaard diagram, for a fixed
diagram $\lsupv{\beta}\HD^\alpha$, different generators can have different cardinalities.

The type \AAm\ module $\CFAAa(\lsupv{\beta}\HD^\alpha)$ associated to $\lsupv{\beta}\HD^\alpha$ is generated
over $\Ground$ by the generators $\x$. The boundary of $\lsupv{\beta}\HD^\alpha$ consists of
pointed matched circles $\PMC_L^\beta$ (coming from the $\beta$-arcs) and
$\PMC_R^\alpha$ (coming from the $\alpha$-arcs). We define an action of the
idempotents of $\Alg(\PMC_L)$ and $\Alg(\PMC_R)$ on $\CFAAa(\lsupv{\beta}\HD^\alpha)$ as
follows. Let $s$ be a subset of the $\beta$-arcs and $t$ a subset of
the $\alpha$-arcs. Then define
\[
I(s)\cdot\x\cdot I(t)=
\begin{cases}
  \x & \text{if $s$ (resp.~$t$) are exactly the $\beta$-arcs (resp.~$\alpha$-arcs)
    occupied by $\x$}\\
  0 & \text{otherwise.}
\end{cases}
\]
This action extends to actions of the rest of $\Alg(\PMC_L)$ and
$\Alg(\PMC_R)$ by counting holomorphic curves in the usual way.

The type \DD\ and \DA\ modules associated to $\lsupv{\beta}\HD^\alpha$ are defined
similarly.

For bimodules as for modules, there are two different geometric
versions of algebraic duality.
\begin{definition}
  Given a bordered Heegaard diagram with two boundary components
  $\HD=(\Sigma,\alphas,\betas,\mathbf{z})$, let $-\HD$ denote the bordered
  Heegaard diagram obtained from $\HD$ by reversing the orientation on
  $\Sigma$, and $\overline{\HD}$ the bordered Heegaard diagram
  obtained from $\HD$ by calling the old $\alpha$-curves the new $\beta$ curves,
  and the old $\beta$-curves the new $\alpha$-curves.
\end{definition}
(Compare Definition~\ref{def:dual-HD-1}.)
Lemma~\ref{lem:beta-bord} has obvious analogues for bordered Heegaard
diagrams with two boundary components:
\begin{lemma}\label{lem:beta-bord-bimod}
  Let $(Y,\phi_L\co F(\PMC_L^\alpha)\to \bdy_LY,\phi_R\co
  F(\PMC_R^\alpha)\to \bdy_RY)$ be a strongly bordered $3$-manifold
  with two boundary components.  Let $\lsup{\alpha}\HD^\alpha$ be an
  $\alpha$-$\alpha$-bordered Heegaard diagram for
  $(Y,\phi_L,\phi_R)$. Then $\lsup{\beta}\overline{\HD}{}^\beta$ is a
  $\beta$-bordered Heegaard diagram for $(-Y,\phi_L\circ
  K_{\beta,\alpha},\phi_R\circ K_{\beta,\alpha})$. Similar statements
  hold in the cases that $Y$ is $\alpha$-$\beta$ bordered or
  $\beta$-$\beta$-bordered.
\end{lemma}

\begin{proof}
  This follows similarly to Lemma~\ref{lem:beta-bord}.
\end{proof}

\begin{proposition}\label{prop:beta-bimod-is-dual}
  If $\HD$ is a bordered Heegaard diagram with two boundary
  components then $\CFDDa(\ol{\HD},-\spinc)$ is dual to $\CFDDa(\HD,\spinc)$;
  $\CFAAa(\ol{\HD},-\spinc)$ is dual to $\CFAAa(\HD,\spinc)$; and $\CFADa(\ol{\HD},-\spinc)$ is
  dual to $\CFDAa(\HD,\spinc)$.
\end{proposition}

\begin{proposition}\label{prop:minus-hd-bimod}
  If $\HD$ is a bordered Heegaard diagram with two boundary
  components then $\CFDDa(-\HD,-\overline{\spinc})$ is dual to
  $\CFDDa(\HD,\spinc)$;
  $\CFAAa(-\HD,-\overline{\spinc})$ is dual to $\CFAAa(\HD,\spinc)$; and $\CFDAa(-\HD,-\overline{\spinc})$ is
  dual to $\CFADa(\HD,\spinc)$.
\end{proposition}
(See Definition~\ref{def:dual-bimod} for the definitions of the duals
of various kinds of bimodules.)
\begin{proof}[Proof of Propositions~\ref{prop:beta-bimod-is-dual}
  and~\ref{prop:minus-hd-bimod}]
  These propositions follow from similar arguments to
  Proposition~\ref{prop:beta-is-dual}.  See also
  Proposition~\ref{prop:minus-hd-mod}.
\end{proof}


\section{An interpolating piece}\label{sec:denis}
\begin{figure}
  \centering
  \includegraphics{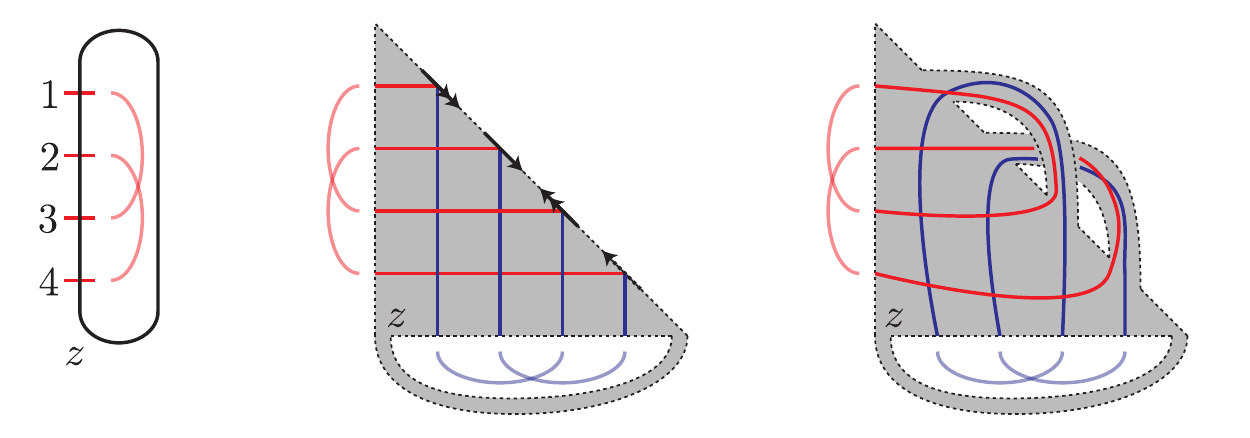}
  \caption{\textbf{An example of $\Denis(\PMC)$.} The example shown is
    for the genus $1$ pointed matched circle. Left: the pointed
    matched circle $\PMC$. Center and Right: two different depictions
    of $\Denis(\PMC)$.}
  \label{fig:Denis}
\end{figure}

Fix a pointed matched circle $\PMC$.  Our discussion of orientation
reversal and proof of the $\Hom$ pairing theorem will rely on a
particular arced $\alpha\Hyph\beta$-bordered Heegaard diagram
associated to $\PMC$, first introduced by
Auroux~\cite{AurouxBordered}.
The diagram,
which we will denote $\Denis(\PMC)$, is constructed as illustrated in
Figure~\ref{fig:Denis} and described below.

Let $k$ denote the genus of $F(\PMC)$. Let $T$ denote the triangle in
$\RR^2$ bounded by the $x$-axis, the $y$-axis, and the line
$y+x=4k+1$. Let $e_y$ (respectively $e_x$) denote the edge of $T$
along the $y$-axis (respectively $x$-axis) and $e_D$ the edge of $T$
along the line $x+y=4k+1$. Let $\Sigma'$ denote the quotient space of
$T$ in which one identifies a small neighborhoods in $e_D$ of the
points $(i,4k+1-i)$ and $(j,4k+1-j)$ if $i$ and $j$ are matched in
$\PMC$, in such a way that $\Sigma'$ is an orientable surface of genus
$k$ with one boundary component.

If $i$ and $j$ are matched in $\PMC$ then the two vertical line segments
$T\cap\{x=i\}$ and $T\cap\{x=j\}$ descend to give a single
arc in $\Sigma'$. Similarly, if
$i$ and $j$ are matched then the horizontal line segments
$T\cap\{y=4k+1-i\}$ and $T\cap\{y=4k+1-j\}$ descend to give a single arc in
$\Sigma'$. Let 
\begin{align*}
  \betas&=\bigcup_{i=1}^{4k}\{x=-i\}\subset\Sigma', &
  \alphas&=\bigcup_{i=1}^{4k}\{y=i-4k-1\}\subset\Sigma'.
\end{align*}

Finally, attach a $1$-handle to $\bdy \Sigma'$ between the points
$(0,0)$ and $(4k+1,0)$. Call the result $\Sigma$. Let $\z$
denote a neighborhood of the
core of this $1$-handle. So, $\z$ is a rectangle in $\Sigma$ connecting
the two boundary components.

Let $\Denis(\PMC)$ denote the diagram $(\Sigma,\alphas,\betas,\z)$.
We let $\bdy_R\Denis(\PMC)$ denote the boundary component of
$\Denis(\PMC)$ which intersects the $\beta$-arcs (so $\Denis(\PMC)$ is
a $(\alpha,\beta)$-bordered Heegaard diagram). Note that
$\bdy_L\Denis(\PMC)$ and $\bdy_R\Denis(\PMC)$ are twin pointed matched
circles (in the sense of Definition~\ref{def:Twins}).
Sometimes, we denote this diagram by $\lsupv{\alpha}\Denis(\PMC)^\beta$
to call attention to the fact that it is $(\alpha,\beta)$-bordered, to 
distinguish it from 
$\lsupv{\beta}\Denis(\PMC)^\alpha$, which is the same diagram but with
the roles of $\partial_L$ and $\partial_R$ reversed. In particular,
just as for $\lsupv{\alpha}\Denis(\PMC)^\beta$, 
$\bdy(\lsupv{\beta}\Denis(\PMC)^\alpha)=\PMC^\beta\amalg\PMC^\alpha$;
and the bimodules $\CFAAa(\lsupv{\alpha}\Denis(\PMC)^\beta)$ and
$\CFAAa(\lsupv{\beta}\Denis(\PMC)^\alpha)$ are canonically isomorphic.

The following was proved by Auroux~\cite{AurouxBordered}.  We recall the proof
briefly here.
\begin{proposition}\label{prop:CFAA-of-Denis}
  The type \AAm\ module $\CFAAa(\Denis(\PMC))$ associated to the
  diagram 
  $\Denis(\PMC)$, viewed as a left-right
  $\Alg(\PMC)$-$\Alg(\PMC)$-bimodule, is isomorphic to the bimodule
  $\Alg(\PMC)$.
\end{proposition}
\begin{proof}[Proof sketch.]
  First, observe that the generators $\Gen(\Denis(\PMC))$ are in
  one-to-one correspondence with the standard basis for
  $\Alg(\PMC)$ by strand diagrams. Indeed, numbering the
  $\alpha$-circles from the bottom
  and the $\beta$-circles from the left, notice that the number of
  points in $\alpha_i\cap \beta_j$ is $2$ if $i=j$, and otherwise the
  number of points is exactly the number of Reeb chords
  in $\PMC$ starting at an endpoint of $\alpha_i$ and ending at an
  endpoint of~$\alpha_j$. These intersections correspond to individual
  strands in a strand diagram: the intersection of $\alpha_i \cap
  \beta_i$ on the diagonal~$e_D$ corresponds to a smeared horizontal
  strand, and other intersections correspond to Reeb chords or
  upward-sloping strands.  An arbitrary generator of
  $\CFAAa(\Denis(\PMC))$ is a set of such intersection
  points, which thus correspond naturally to a strand diagram in the
  standard basis for $\Alg(\PMC)$.

  To compute the $\Ainf$-bimodule structure, one observes that
  $\Denis(\PMC)$ is nice (see~\cite{SarkarWang07:ComputingHFhat}), so
  the differential on
  $\CFAAa(\Denis(\PMC))$ comes
  entirely from counting rectangles (there are no interior bigons), the only
  multiplications are $m_2$'s, and these multiplications count
  half-strips. With this explicit description, it is straightforward
  to identify the differential and algebra actions on $\CFAAa(\Denis(\PMC))$
  with those on $\Alg(\PMC)$.  We refer the reader
  to~\cite{AurouxBordered} for more details.
\end{proof}

Note that since $\Denis(\PMC)$ has no closed $\alpha$- or $\beta$-circles,
$Y(\Denis(\PMC))$ as an unparametrized $3$-manifold is
$F(\PMC)\times[0,1]$. Consequently, $\Denis(\PMC)$ determines a
homeomorphism from $F(\PMC)$ to itself preserving the preferred disk
$D_\alpha\cup s\cup D_\beta$, well-defined up to isotopy
fixing $D_\alpha\cup s\cup D_\beta$. We understand this map as follows.

\begin{figure}
  \centering
  \includegraphics{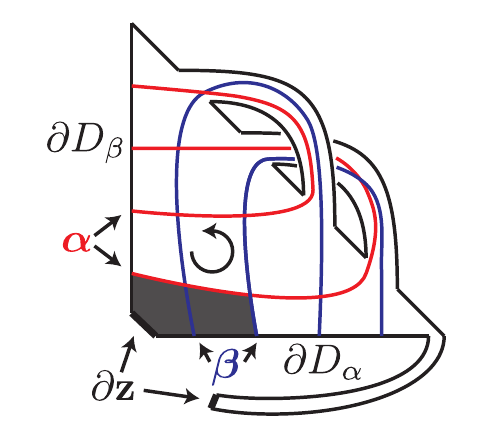}
  \caption{\textbf{Identification of the homeomorphism associated to
      $\Denis(\PMC)$.} The diagram $\Denis(\PMC)\setminus{\bf z}$ for $\PMC$
    the pointed matched circle of genus $1$ is shown. The small
    triangle identifying the homeomorphism as a positive half boundary
    Dehn twist is shaded.}
  \label{fig:DenispunctF}
\end{figure}
\begin{proposition}\label{prop:Denis-reps}
  The diagram $\Denis(\PMC)$ represents a positive half Dehn twist, in
  the following sense.
  Let $\phi_{\lsupv{\alpha}\Denis^\beta}\co {-F(\PMC^\alpha)}\to
  F(\PMC^\beta)$ denote the homeomorphism associated (as in
  Construction~\ref{constr:AssocHomeo}) to the diagram
  $\lsupv{\alpha}\Denis(\PMC)^\beta$, and
  $K_{\beta,\alpha}\co F(\PMC^\beta)\to {-F(\PMC^\alpha)}$ denote the
  canonical homeomorphism of twins (as in Definition~\ref{def:Twins}).
  Then $K_{\beta,\alpha}\circ \phi_{\lsupv{\alpha}\Denis^\beta}=\tau_\bdy^{1/2}\co
  {-F(\PMC^\alpha)}\to -F(\PMC^\alpha)$.
  Likewise,
  $\phi_{\lsupv{\alpha}\Denis^\beta} \circ K_{\beta,\alpha} \co
    F(\PMC^\beta) \to F(\PMC^\beta)$,
  $K_{\alpha,\beta} \circ \phi_{\lsupv{\beta}\Denis^\alpha} \co
    {-F(\PMC^\beta)} \to {-F(\PMC^\beta)}$, and 
  $\phi_{\lsupv{\beta}\Denis^\alpha} \circ K_{\alpha,\beta}\co
    F(\PMC^\alpha) \to F(\PMC^\alpha)$ all represent positive half
  Dehn twists on the respective surfaces.
\end{proposition}

\begin{proof}
  We concentrate on the first case, $K_{\beta,\alpha}\circ
  \phi_{\lsupv{\alpha}\Denis^\beta}$.
  Let $Y$ be the result of applying
  Construction~\ref{constr:AlphaBetaBorderedDiagramToThreeManifold}
  to the diagram $\Denis(\PMC)=(\bSigma,\alphas,\betas,\arcz)$. Since
  there are no closed $\alpha$- or $\beta$-circles, $Y$ is given by
  \[
  [0,1]\times(D_\alpha\cup_{\bdy_L\bSigma}\cup \bSigma\cup_{\bdy_R\bSigma} D_\beta).
  \]
  With notation as in Construction~\ref{constr:AssocHomeo}:
  \begin{enumerate}
  \item\label{item:phiL} The map $\phi_L\co F(\PMC^\alpha)\to
    \{0\}\times (D_\alpha\cup
    \bSigma\cup D_\beta)$ sends $D_\alpha$ to
    $D_\alpha$; $D_\beta$ to $D_\beta$; and
    $F(\PMC^\alpha)$ to $\Sigma$, sending the cores of the
    $1$-handles in $F(\PMC^\alpha)$ to the $\alpha$-arcs.
  \item\label{item:phiR} The map $\phi_R\co F(\PMC^\beta)\to
    \{1\}\times (D_\alpha\cup
    \bSigma\cup D_\beta)$ sends $D_\alpha$ to
    $D_\alpha$; $D_\beta$ to $D_\beta$; and
    $\PunctF(\PMC^\beta)$ to $\Sigma$, sending the cores of the
    $1$-handles in $\PunctF(\PMC^\beta)$ to the $\beta$-arcs.
  \item\label{item:Phi} The map $\Phi\co [0,1]\times F(\PMC^\beta)\to [0,1]\times(D_\alpha\cup \bSigma\cup D_\beta)$
    is given by $\Phi(t,x)=(t,\phi_R(x))$.
  \end{enumerate}
  In particular, by~(\ref{item:Phi}), the map $\phi_\Denis=\phi_Y$ is given by
  $\phi_R^{-1}\circ (-\phi_L)\co {-F(\PMC^\alpha)\to F(\PMC^\beta)}$.

  So, we have a commutative diagram:
  \[
  \xymatrix{
    -F(\PMC^\alpha) \ar[rr]^{\phi_{\Denis}=\phi_R^{-1}\circ\phi_L}\ar[d]_{-\phi_L}& & F(\PMC^\beta)\ar[rr]^{K_{\beta,\alpha}}\ar[d]^{\phi_R}
   & & -F(\PMC^\alpha)\ar[d]^{-\phi_L}\\
   D_\alpha\cup \bSigma\cup D_\beta \ar[rr]_{\Id} & &  D_\alpha\cup \bSigma\cup D_\beta\ar[rr]_{g} & & D_\alpha\cup \bSigma\cup D_\beta.
  }
  \]
  The map $K_{\beta,\alpha}$ takes the cores of the $1$-handles in
  $F(\PMC^\beta)$ to the cores of the $1$-handles in
  $F(\PMC^\alpha)$. So, $-\phi_L\circ K_{\beta,\alpha}\circ (\phi_R)^{-1}=g$
  is the map from $D_\alpha\cup\bSigma\cup D_\beta$ to itself
  exchanging $D_\alpha$ and $D_\beta$ and taking each $\beta$-arc to
  the corresponding $\alpha$-arc.  Hence, by
  Lemma~\ref{lem:CharHalfDehnTwist}, the map $g$ is a positive half
  boundary Dehn twist of $D_\alpha\cup\bSigma\cup D_\beta$, and so
  $K_{\beta,\alpha}\circ \phi_\Denis$ is a positive half boundary Dehn
  twist of $-F(\PMC^\alpha)$.

  The other cases can be proved by similar commutative diagrams, or
  alternatively follow from the first case, using the observations
  that, on the one hand, $K_{\alpha,\beta} \circ \tau_\bdy^{1/2} \circ
  K_{\beta,\alpha} = \tau_\bdy^{1/2}$ and, on the other hand,
  switching the left and right sides of $M_\phi$ yields
  $M_{-\phi^{-1}}$ and $(-\tau_\bdy^{1/2}) = \tau_\bdy^{1/2}$.  (In
  both cases, the surface on which we apply $\tau_\bdy^{1/2}$
  changes.)
\end{proof}

\begin{figure}
  \centering
  \includegraphics{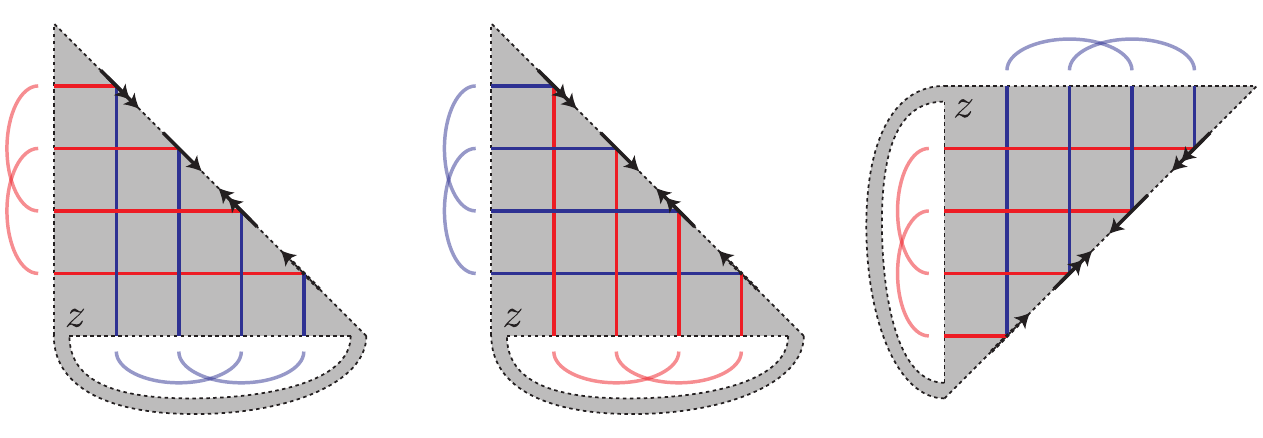}
  \caption{\textbf{The diagram $\MirrorDenis(\PMC)$.} Left: the diagram $\Denis(\PMC)$, for the genus $1$ pointed matched circle. Center: the result of exchanging the $\alpha$- and $\beta$-curves in $\Denis(\PMC)$. Right: the diagram $\MirrorDenis(\PMC)$.}
\label{fig:mirror-denis}
\end{figure}

There is another interpolating piece that will be important, a kind of
mirror image of $\Denis(\PMC)$. More precisely, we can consider the
diagram $\MirrorDenis(\PMC)$ obtained from
$\Denis(\PMC)$ by switching the $\alpha$- and $\beta$-curves and rotating the
diagram clockwise $90$ degrees.
The boundary components of
$\MirrorDenis(\PMC)$ are naturally identified with
$\PMC^\alpha$ and $\PMC^\beta$. Alternatively, $\MirrorDenis(\PMC)$ is $-\Denis(-\PMC)$, 
obtained from $\Denis(-\PMC)$ by reversing the orientation on the
diagram
(e.g., by reflecting across the $x$-axis). 
See Figure~\ref{fig:mirror-denis}. By default, as with $\Denis$, we view $\PMC^\alpha$
as the left boundary of $\MirrorDenis$ and $\PMC^\beta$ as the right
boundary of $\MirrorDenis$, but when we want to make this convention
or the opposite one explicit we will write
$\lsupv{\alpha}\MirrorDenis{}^\beta$ or
$\lsupv{\beta}\MirrorDenis{}^\alpha$, respectively.

With no additional work, we get the following algebraic description
of the bordered invariants for $\MirrorDenis(\PMC)$:
\begin{proposition}\label{prop:CFAA-of-MirrorDenis}
  The type \AAm\ module associated to $\MirrorDenis(\PMC)$ is
  isomorphic to $\overline{\Alg(\PMC)}$, the dual module to
  $\Alg(\PMC)$.
\end{proposition}

\begin{proof}
    This follows from Propositions~\ref{prop:CFAA-of-Denis}
    and~\ref{prop:beta-bimod-is-dual}.
\end{proof}

\begin{proposition}\label{prop:mdenis-reps}
  Let $\phi_{\lsupv{\alpha}\MirrorDenis^\beta}\co {-F(\PMC^\alpha)}\to
  F(\PMC^\beta)$ denote the homeomorphism associated (as in
  Construction~\ref{constr:AssocHomeo}) to the diagram
  $\lsupv{\alpha}\MirrorDenis^\beta$, and
  $K_{\beta,\alpha}\co F(\PMC^\beta)\to {-F(\PMC^\alpha)}$ denote the
  canonical homeomorphism of twins (as in Definition~\ref{def:Twins}).
  Then $K_{\beta,\alpha}\circ
  \phi_{\lsupv{\alpha}\MirrorDenis^\beta}=\tau_\bdy^{-1/2}\co
  {-F(\PMC^\alpha)}\to -F(\PMC^\alpha)$.
  Likewise,
  $\phi_{\lsupv{\alpha}\MirrorDenis^\beta} \circ K_{\beta,\alpha} \co
    F(\PMC^\beta) \to F(\PMC^\beta)$,
  $K_{\alpha,\beta} \circ \phi_{\lsupv{\beta}\MirrorDenis^\alpha} \co
    {-F(\PMC^\beta)} \to {-F(\PMC^\beta)}$, and 
  $\phi_{\lsupv{\beta}\MirrorDenis^\alpha} \circ K_{\alpha,\beta}\co
    F(\PMC^\alpha) \to F(\PMC^\alpha)$ all represent negative half
  Dehn twists on the respective surfaces.
\end{proposition}
\begin{proof}
  This follows along the lines of Proposition~\ref{prop:Denis-reps}:
  combine the construction of the associated homeomorphism with
  Lemma~\ref{lem:CharHalfDehnTwist}. Alternatively, note that exchanging
  the $\alpha$- and $\beta$-circles has the effect of reversing the
  orientation on the three-manifold. This, in turn, exchanges positive
  and negative Dehn twists, so the result follows from
  Proposition~\ref{prop:Denis-reps}.
\end{proof}

In Section~\ref{sec:Consequences}, the following corollary of
Propositions~\ref{prop:Denis-reps} and~\ref{prop:mdenis-reps} will
be useful, particularly as sanity checks on the
signs / presence of boundary Dehn twists.
\begin{corollary}\label{cor:denis-glue-reps}
  The diagram
  $\lsupv{\alpha}\Denis(-\PMC)^\beta\cup\lsupv{\beta}\Denis(\PMC)^\alpha$
  represents $\tau_\bdy\co F(\PMC^\alpha)\to F(\PMC^\alpha)$.
  The diagram
  $\lsupv{\alpha}\Denis(-\PMC)^\beta\cup\lsupv{\beta}\MirrorDenis(\PMC)^\alpha$
  represents the identity map of $F(\PMC^\alpha)$.
\end{corollary}

See Figure~\ref{fig:Identities}.
More identities of this kind are given at the end of
Appendix~\ref{sec:conventions}.

\begin{proof}
  For the first statement, we have
  \begin{align*}
    \phi_{\lsupv{\alpha}\Denis(-\PMC)^\beta\cup\lsupv{\beta}\Denis(\PMC)^\alpha}&=\phi_{\lsupv{\beta}\Denis(\PMC)^\alpha}\circ
    \phi_{\lsupv{\alpha}\Denis(-\PMC)^\beta}\\
    &=(\phi_{\lsupv{\beta}\Denis(\PMC)^\alpha}\circ
    K_{\alpha,\beta})\circ (K_{\beta,\alpha}\circ
    \phi_{\lsupv{\alpha}\Denis(-\PMC)^\beta})\\
    &= (\tau_\bdy^{1/2}\co F(\PMC^\alpha)\to
    F(\PMC)^\alpha)\circ (\tau_\bdy^{1/2}\co
    F(\PMC^\alpha)\to F(\PMC^\alpha))\\
    &=\tau_\bdy\co F(\PMC^\alpha)\to F(\PMC^\alpha).
  \end{align*}
  The second statement is similar, except the two half boundary twists
  go in opposite directions and so cancel.
\end{proof}

  \begin{figure}
    \centering
    \input{Identities}
    \caption{\textbf{Gluing $\Denis(\PMC)$ to itself and its inverse.}
      \label{fig:Identities}
      The left picture represents the Heegaard diagram
      $\lsupv{\alpha}\Denis(-\PMC)^\beta\cup\lsupv{\beta}\Denis(\PMC)^\alpha$
      (which
      in turn represents $\tau_\bdy^{-1}$), while the right illustrates
      $\lsupv{\alpha}\Denis(-\PMC)^\beta\cup\lsupv{\beta}\MirrorDenis(\PMC)^\alpha$
      (which represents the identity).  }
  \end{figure}

\begin{lemma}
  \label{lem:DenisCommutes}
  Let $\lsupv{\alpha}\HD^\alpha$ be an $\alpha$-$\alpha$-bordered Heegaard diagram,
  with $\partial_L\HD = \PMC_L$ and $\partial_R \HD=\PMC_R$.
  Then the two $(\alpha,\beta)$-bordered Heegaard diagrams
  \[
  \lsupv{\alpha}\HD^\alpha\sos{\partial_R}{\cup}{\partial_L}
  \lsupv{\alpha}\Denis(-\PMC_R)^\beta
  \qquad\text{and}\qquad
  \lsupv{\alpha}\Denis(\PMC_L)^\beta
  \sos{\partial_R}{\cup}{\partial_L}
  \lsupv{\beta}(-\ol{\HD})^\beta
  \]
  represent the same strongly bordered three-manifold.
\end{lemma}

The intuition here is that the negative half Dehn twist represented by
$\Denis$ can be ``pulled through'' from one side of $Y(\HD)$ to the
other, but we have to turn over the Morse function on $Y(\HD)$ in the
process.  (Note that $Y(-\ol{\HD})$ is orientation-preserving
homeomorphic to $Y(\HD)$.)

\begin{proof}
  Let $Y$ be a three-manifold with two boundary components. Then $Y$
  can be factored as a product of elementary cobordisms, each of which
  corresponds to attaching a
  one-handle or two-handle to $\partial_L Y$. Moreover, a strongly
  bordered three-manifold can be factored into the following
  simple pieces:
  \begin{itemize}
    \item mapping cylinders for homeomorphisms, and
    \item elementary cobordisms from $\PMC\#\PMC_1$ to $\PMC$, or from
      $\PMC$ to  $\PMC\#\PMC_1$, where
      $\PMC_1$ denotes the genus $1$ pointed matched circle, obtained
      by attaching a two-handle along the $\infty$-framed curve in
      $\PunctF(\PMC_1)\subset F(\PMC\#\PMC_1)$ (as in
      Figure~\ref{fig:CommuteDenisElem}).
  \end{itemize}
  In particular, any
  Heegaard diagram $\HD_0$ is equivalent to a diagram $\HD$ which can be
  written as a juxtaposition of pieces of these two forms, and obviously
  $-\ol{\HD_0}$ is then equivalent to $-\ol{\HD}$. Thus,
  it suffices to check the result for these two kinds of elementary pieces.

  For the first simple piece,
  let $\HD_\phi$ be the standard Heegaard diagram
  for a mapping class $\phi\co F(\PMC_1) \to F(\PMC_2)$, and let
  $\Denis(\PMC_1)(\alphas,\phi(\betas))$ (for instance) be
  the diagram defined like $\lsupv{\alpha}\Denis(\PMC_1)^\beta$ but
  using the image of the
  $\beta$-curves on $\PMC_2$ under the mapping class $\phi$, thought
  of as acting
  on the complement of $\z$ in the 
  Heegaard surface.  
  (Thus $\Denis(\PMC_1)(\alphas, \betas) = \lsupv{\alpha}\Denis(\PMC_1)^\beta$.)
  Then we have
  equivalences of Heegaard diagrams:
  \begin{align*}
    \HD_\phi
    \cup
    \Denis(-\PMC_2)(\alphas,\betas)
    & \simeq
    \Denis(-\PMC_2)(\phi(\alphas),\betas) \\
    &\cong
    \Denis(-\PMC_1)(\alphas,\phi^{-1}(\betas))\\
    & \simeq 
    \Denis(-\PMC_1)(\alphas,\betas)
    \cup (-\ol{\HD_\phi}),
  \end{align*}
  where we are using the the fact that gluing mapping cylinders
  corresponds to twisting the
  parametrization~\cite[Lemma~\ref{LOT2:lem:BorderedForDiffeo}]{LOT2},
  a homeomorphism of diagrams, and the gluing property again.

  For the second simple piece (elementary cobordisms of the specified
  form), the needed sequence of handleslides and destabilizations is
  easy to find; see Figure~\ref{fig:CommuteDenisElem}, where we have
  illustrated the case where $\PMC$ is empty.
\end{proof}

  \begin{figure}
    \centering
    \input{CommuteDenisElem}
    \caption{\textbf{Commuting $\Denis$ past an elementary cobordism.}
      Heegaard moves exhibiting the identification
      $\lsup{\alpha}\HD^\alpha\cup \Denis(-\PMC_R)\simeq
      \Denis(\PMC_R)\cup \lsup{\beta}(-\ol{\HD})^\beta$, in the case
      where $\HD$ represents an elementary cobordism from the
      two-sphere to a genus one surface (so $\Denis(\PMC_L)$ is
      empty).  The first diagram is $\lsup{\alpha}\HD^\alpha\cup
      \Denis(-\PMC_R)$, the second is gotten by a sequence of
      handleslides, the third by a destabilization, and the fourth by
      an isotopy (and a homeomorphism of diagrams).}
    \label{fig:CommuteDenisElem}
  \end{figure}


\section{Consequences}\label{sec:Consequences}
\subsection{Orientation reversal and the \textalt{$\Hom$}{Hom} pairing
  theorem for modules}
As a warm-up for the proofs of our main theorems, we start with the
module case, in which the notation is
a little simpler to follow, and the Dehn twists disappear.
\begin{proof}[Proof of Theorem~\ref{thm:or-rev}]
  Fix an $\alpha$-bordered Heegaard diagram $\HD^\alpha$ for $(Y,\phi\co
  F(\PMC)\to Y)$.
  By Lemma~\ref{lem:DenisCommutes} in the case when one boundary is
  empty (which is essentially the fact that boundary Dehn twists have
  no effect on $3$-manifolds with just one boundary component),
  the $\alpha$-bordered
  Heegaard diagram
  $\ol{\HD}^\beta\sos{\bdy}{\cup}{\bdy_R}\Denis(-\PMC)$
  represents $(-Y,\phi\co F(-\PMC)\to -Y)$.  We have now that
  \begin{align*}
    \CFAa(-Y) &\cong \CFAa(\ol{\HD}^\beta\sos{\bdy}{\cup}{\bdy_R}\Denis(-\PMC)) \\
    &\simeq 
    \CFDa(\ol{\HD}^\beta)\DT\CFAAa(\Denis(-\PMC))\\
    &\cong \CFDa(\ol{\HD}^\beta)\DT\Alg(-\PMC)\\
    &\cong \overline{\CFDa(\HD^\alpha)}\DT\Alg(-\PMC)\\
    &\simeq \Mor_{\Alg(-\PMC)}(\CFDa(Y),\Alg(-\PMC)).
  \end{align*}
  Here, the second line uses the pairing theorem (in a form using bimodules,
  see~\cite[Theorem~\ref*{LOT2:thm:GenReparameterization}]{LOT2}), the third uses
  Proposition~\ref{prop:CFAA-of-Denis}, the fourth uses
  Proposition~\ref{prop:beta-is-dual}, and the last uses
  Proposition~\ref{prop:D-Mor-is-DT}.  This proves
  Equation~\eqref{eq:ReverseTypeD}.

  Equation~\eqref{eq:ReverseTypeA} follows from Equation~\eqref{eq:ReverseTypeD}
  and Theorem~\ref{thm:DDId}. That is, 
  \begin{align*}
    \CFAa(-Y)&\simeq \Mor^{\Alg(-\PMC)}(\CFDa(Y),\lsupv{\Alg(-\PMC)}[\Id]_{\Alg(-\PMC)})\\
    &\simeq \Mor_{\Alg(\PMC)}(\CFAa(Y),\CFAAa(\Id)),
  \end{align*}
  where the first equivalence is Equation~\eqref{eq:ReverseTypeD} and
  the second follows by tensoring both domain and range with
  $\CFAAa(\Id)$, which
  is an equivalence of categories (Theorem~\ref{thm:DDId}). Tensoring
  both sides with $\CFDDa(\Id)$ gives Equation~\eqref{eq:ReverseTypeA}.
\end{proof}

\begin{figure}
  \centering
  \includegraphics{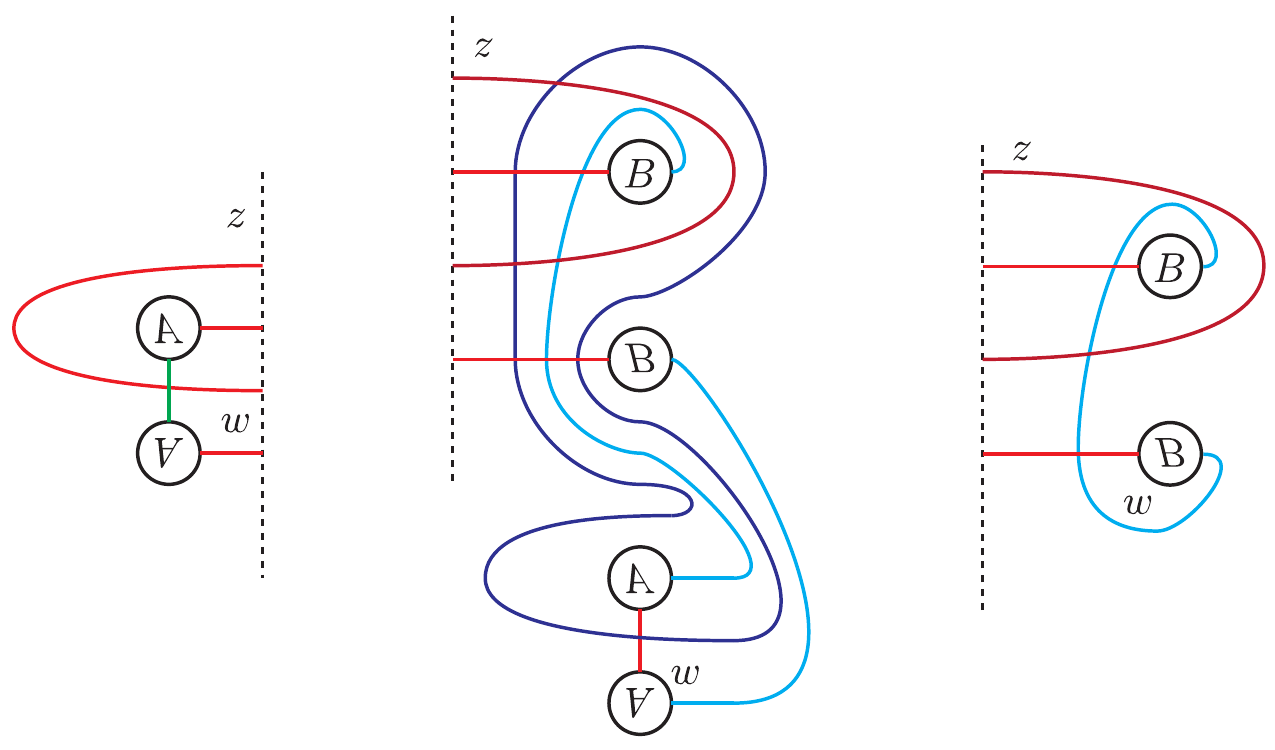}
  \caption{\textbf{Gluing $\Denis(\PMC)$ to a doubly-pointed Heegaard
      diagram.} Left: A doubly-pointed $\alpha$-bordered Heegaard
    diagram $\HD$ for the core of a $0$-framed solid torus. Center:
    the result $\Denis(\PMC)\cup_\bdy \HD^\beta$ of gluing the
    interpolating piece to $\HD^\beta$. Right: a destabilization of
    $\Denis(\PMC)\cup_\bdy \HD^\beta$.}
  \label{fig:double-denis}
\end{figure}

\begin{remark}\label{rmk:borsut}
  In the context of bordered sutured manifolds~\cite{Zarev09:BorSut}
  with upper and lower pieces of the boundary $R^+$ and $R^-$,
  the operation 
  $\HD\leadsto \ol{\HD}\sos{\bdy}{\cup}{\bdy_R}\Denis(\PMC)$ 
  corresponds to a
  slightly more complicated operation than mere orientation reversal.
  When one switches the $\alpha$- and $\beta$-curves on a diagram for
  a bordered sutured manifold,
  the roles of $R_+$ and $R_-$ are exchanged.  Attaching
  $Y(\Denis(\PMC))$ then corresponds to introducing a half Dehn twist
  along the preferred disk in $F(\PMC)$, as in
  Proposition~\ref{prop:Denis-reps}.  (As above, Dehn twists around
  the preferred disk in $F(\PMC)$ disappear if one instead glues to a
  $3$-manifold with one boundary component.)  See
  Figure~\ref{fig:double-denis} for an example where this operation is
  non-trivial.
\end{remark}

\begin{proof}[Proof of Theorem~\ref{thm:hom-pair}]
  Using Proposition~\ref{prop:D-Mor-is-DT}, Theorem~\ref{thm:or-rev},
  and the usual version of the pairing
  theorem~\cite[Theorem~\ref*{LOT:thm:TensorPairing}]{LOT1} in turn
  gives:
  \begin{align*}
    \Mor_{\Alg(-\PMC)}(\CFDa(Y_1),\CFDa(Y_2))
    &\cong \Mor_{\Alg(-\PMC)}(\CFDa(Y_1),\Id)\DT\CFDa(Y_2) \\
    &\simeq \CFAa(-Y_1)\DT \CFDa(Y_2) \\
    &\simeq \CFa(-Y_1\cup_{\partial} Y_2),
  \end{align*}
  Taking homology gives the first isomorphism from
  Theorem~\ref{thm:hom-pair}. The second isomorphism then follows from
  Corollary~\ref{cor:D-hom-is-A-hom}.
\end{proof}

\subsection{Conjugation invariance}

\begin{proposition}
  \label{prop:JMap}
  Suppose that $\lsupv{\alpha}\HD^\alpha$ is an $\alpha\Hyph\alpha$-bordered Heegaard
  diagram with boundary components $\PMC_1$ and $\PMC_2$.
  By the isomorphisms $\Alg(-(-\PMC_i))^{\op}\cong \Alg(\PMC_i)^{\op}\cong \Alg(-\PMC_i)$,
  we can view both $\CFDDa(\lsupv{\alpha}\HD^\alpha)$ and $\CFDDa(-\lsupv{\beta}\ol{\HD}^\beta)$ as 
  left-left
  $\Alg(-\PMC_1)\otimes \Alg(-\PMC_2)$-modules. Under this
  identification and similar ones for other modules, we have
  isomorphisms
  \begin{align*}
  \CFDDa(\HD,\spinc)&\cong \CFDDa(-\ol{\HD},{\overline\spinc})\\
  \CFAAa(\HD,\spinc)&\cong\CFAAa(-\ol{\HD},{\overline\spinc})\\
  \CFDAa(\HD,\spinc)&\cong\CFDAa(-\ol{\HD},{\overline\spinc}).
  \end{align*}
\end{proposition}

\begin{proof}
  The identification of the complexes $\CFDDa(\HD)$ with
  $\CFDDa(-\ol{\HD})$ is supplied by combining
  Propositions~\ref{prop:beta-bimod-is-dual} and~\ref{prop:minus-hd-bimod}.  The
  conjugation on the $\SpinC$ structures comes from the observation
  that, when $\x$ is viewed as
  a generator for the bordered Floer homology of $\HD$, its
  corresponding vector field points in the opposite direction from
  that of $\ol{\x}$ when viewed as a generator for the bordered Floer homology
  of $-\ol{\HD}$.
\end{proof}

\begin{proof}[Proof of Theorems~\ref{thm:D-is-A}
  and~\ref{thm:DA-bimod}]
  We will prove Theorem~\ref{thm:DA-bimod};
  Theorem~\ref{thm:D-is-A} can be viewed as a special case where one
  of the boundary components is empty (or, if one prefers, $S^2$),
  after noting that $\tau_\bdy$ acts trivially on bordered
  $3$-manifolds with only one boundary component.

  To establish Equation~\eqref{eq:DA-bimod-1}, we must show that
  \begin{multline}\label{eq:Reformulate}
  (\Alg(-\PMC_1)_{\Alg(\PMC_1),\Alg(-\PMC_1)} \otimes
  \Alg(-\PMC_2)_{\Alg(\PMC_2),\Alg(-\PMC_2)})\DT \lsup{\Alg(-\PMC_1),\Alg(-\PMC_2)}\CFDDa(Y, \spinc) \\
  \simeq \CFAAa(\tau_\bdy(Y),\overline{\spinc})_{\Alg(\PMC_1),\Alg(\PMC_2)},
  \end{multline}
  where $\Alg(-\PMC_i)_{\Alg(\PMC_i),\Alg(-\PMC_i)}$ denotes the
  $\Alg(-\PMC_i)$-bimodule $\Alg(-\PMC_i)$, viewed as a module with two
  right actions. 

  Fix an $\alpha\Hyph\alpha$-bordered Heegaard diagram $\lsupv{\alpha}\HD^\alpha$ for
  $Y$ (with $\partial_L\HD=\PMC_1$ and $\partial_R\HD=\PMC_2$), so that 
  $\CFAAa(Y,\spinc)$ is given by the bimodule $\CFAAa(\HD,\spinc)$ associated to $\HD$ and $\spinc$, and also (according to Proposition~\ref{prop:JMap}) by the bimodule
  $\CFAAa(-\overline{\HD},\overline{\spinc})$ associated to the
  $\beta\Hyph\beta$-bordered version $\lsup{\beta}(-\ol{\HD})^\beta$
  and the conjugate $\SpinC$-structure $\overline{\spinc}$.
  Now, glue a copy of $\lsup{\alpha}\Denis(\PMC_1)^{\beta}$ and a copy of
  $\lsup{\beta}\Denis(\PMC_2)^{\alpha}$ to the $\partial_L$ and $\partial_R$
  boundary components of
  $-\ol{\HD}$ respectively. 

  Combining Proposition~\ref{prop:CFAA-of-Denis}, the
  pairing theorem, and Proposition~\ref{prop:JMap}, we get
  \begin{multline}    \label{eq:Conjugate}
    \CFAAa(\lsup{\alpha}\Denis(\PMC_1)^{\beta}\sos{\partial_R}{\cup}{-\PMC_1^{\beta}}
    (-\ol{\HD})\sos{-\PMC_2^{\beta}}{\cup}{\partial_L}\lsup{\beta}\Denis(\PMC_2)^{\alpha},{\overline\spinc})
    \\
    \begin{aligned}
    &\simeq \bigl(\Alg(\PMC_1)_{\Alg(\PMC_1),\Alg(-\PMC_1)}\otimes
    \Alg(\PMC_2)_{\Alg(\PMC_2),\Alg(-\PMC_2)}\bigr)\DT\CFDDa(-\ol{\HD},{\overline\spinc})
    \\  &\simeq
    \bigl(\Alg(\PMC_1)_{\Alg(\PMC_1),\Alg(-\PMC_1)}\otimes
    \Alg(\PMC_2)_{\Alg(\PMC_2),\Alg(-\PMC_2)}\bigr)\DT\CFDDa(\HD,\spinc)
    \\ &\cong 
    \bigl(\Alg(-\PMC_1)_{\Alg(\PMC_1),\Alg(-\PMC_1)}\otimes
    \Alg(-\PMC_2)_{\Alg(\PMC_2),\Alg(-\PMC_2)}\bigr)\DT\CFDDa(\HD,\spinc).
    \end{aligned}
  \end{multline}
  The last isomorphism comes from the tautological identification of bimodules
  $A_{A,A^{\op}}\cong (A^\op)_{A,A^{op}}$, together with the usual
  identification of $\Alg(-\PMC)\cong \Alg(\PMC)^\op$.

  We can alternatively move $\Denis(\PMC_1)$ past $-\ol{\HD}$
  as in Lemma~\ref{lem:DenisCommutes}. This gives
  \begin{multline}
    \label{eq:IntroTwist}
    Y(\lsup{\alpha}\Denis(\PMC_1)^{\beta}\sos{\partial_R}{\cup}{
      -\PMC_1^\beta}
    (-\ol{\HD}) \sos{-\PMC_2^\beta}\cup{\partial_L}
    \lsup{\beta}\Denis(\PMC_2)^{\alpha})\\
    \begin{aligned}
    &\simeq 
    Y(\HD\sos{\partial_R}{\cup}{\partial_L}
    \lsup{\alpha}\Denis(-\PMC_2)^{\beta}\sos{\partial_R}{\cup}{\partial_L}
    \lsup{\beta}\Denis(\PMC_2)^{\alpha})\\
    &\simeq \tau_\partial(Y(\HD)).
    \end{aligned}
  \end{multline}
  Here, the equivalence is of strongly bordered 3-manifolds,
  and the last line follows from Corollary~\ref{cor:denis-glue-reps}.
  Equations~\eqref{eq:IntroTwist} and~\eqref{eq:Conjugate} combine to
  establish Equation~\eqref{eq:Reformulate}, which is equivalent to
  Equation~\eqref{eq:DA-bimod-1}.

  Equation~\eqref{eq:DA-bimod-2} is
  immediate from Equation~\eqref{eq:DA-bimod-1}, by replacing $Y$ with
  $\tau_{\bdy}^{-1}(Y)$, $\spinc$ with~$\overline{\spinc}$, and viewing the
  right actions as left actions.

  Equation~\eqref{eq:DA-bimod-3} follows from the first two parts by
  tensoring Equation~\eqref{eq:DA-bimod-1} applied to~$Y$ with
  Equation~\eqref{eq:DA-bimod-2} applied to $\Id_{\PMC_2}$:
  \begin{multline*}
      \CFDDa(Y, \spinc)_{\Alg(\PMC_1),\Alg(\PMC_2)}\DTP \lsub{\Alg(\PMC_2),\Alg(-\PMC_2)}\CFAAa(\Id_{\PMC_2}, \mathfrak{t})\\
      \simeq \CFAAa(\tau_\bdy(Y),\overline{\spinc})_{\Alg(\PMC_1),\Alg(\PMC_2)}\DTP
      \lsub{\Alg(\PMC_2),\Alg(-\PMC_2)}\CFDDa(\tau_\bdy^{-1},\overline{\mathfrak{t}}),
  \end{multline*}
  which reduces to the desired result.
\end{proof}

\begin{proof}[Proof of Corollary~\ref{cor:conj}]
  Equations~\eqref{eq:cor-conj1} and~\eqref{eq:cor-conj2} are
  immediate from Theorem~\ref{thm:D-is-A} and the observations that
  \begin{align*}
    \Alg(\PMC)_{\Alg(\PMC),\Alg(-\PMC)}\DT\lsup{\Alg(-\PMC)}\CFDa(Y,\spinc)&\simeq
    \Alg(\PMC)_{\Alg(\PMC),\Alg(-\PMC)}\DT(\CFAa(Y,\spinc)\DT\CFDDa(\Id))\\
    \CFAa(Y,\spinc)_{\Alg(\PMC)}&\simeq\CFAAa(\Id)_{\Alg(\PMC),\Alg(-\PMC)}\DT\lsup{\Alg(-\PMC)}\CFDa(Y,\spinc). 
  \end{align*}
  (both of which follow from the pairing theorem).
\end{proof}
\subsection{Orientation reversal and the \textalt{$\Hom$}{Hom} pairing
  theorems for bimodules}\label{sec:hom-pair}
\begin{proof}[Proof of Theorem~\ref{thm:bimod-rev-1}]
  We start by proving Equation~\eqref{eq:bimod-rev-1-1}.
  By \cite[Corollary~\ref*{LOT2:cor:InterpretMor}]{LOT2}, the complex
  of type~$D$ morphisms between two type~$D$ structures is
  quasi-isomorphic to the complex of $\Ainf$-morphisms between their
  modulifications, so
  Equation~\eqref{eq:bimod-rev-1-1} is
  equivalent to
  \begin{equation}
  \Mor^{\Alg_1'}(\Alg_2'\DT\lsup{\Alg_1',\Alg_2'}\CFDDa(Y),\lsup{\Alg_1'}[\Id]_{\Alg_1'})\simeq \CFAAa(-Y)_{\Alg_1',\Alg_2'}.\label{eq:bimod-rev-1-1-precise}
\end{equation}
  Using the definition of $\Mor$ in terms of $\DT$
  (Proposition~\ref{prop:Mor-is-DT-bimod}), we want to show that
  \[
  \overline{\Alg_2'\DT\lsup{\Alg_1',\Alg_2'}\CFDDa(Y)}\DT\Alg_1'\simeq \CFAAa(-Y)_{\Alg_1',\Alg_2'}.
  \]
  Since taking duals respects $\DT$ (Lemma~\ref{lem:dual-resp-DT}), this boils
  down to  
  \begin{equation}
    \overline{\lsup{\Alg_1',\Alg_2'}\CFDDa(Y)}\DT
    \overline{\Alg_2'} \DT\Alg_1'\simeq \CFAAa(-Y)_{\Alg_1',\Alg_2'}.\label{eq:boils-down}
  \end{equation}

  Now, fix a Heegaard diagram $\lsupv{\alpha}\HD^\alpha$, with
  $\partial_L\HD=\PMC_1$ and $\partial_R\HD=\PMC_2$, for $(Y,\phi_L\co
  F(\PMC_1)\to \bdy_LY,\phi_R\co F(\PMC_2)\to\bdy_RY)$. According to
  Lemma~\ref{lem:beta-bord-bimod},
  $\lsupv{\beta}\ol{\HD}^\beta$ is a Heegaard diagram for $-Y$ with
  bordering $\phi_L\circ K_{\beta,\alpha}$ and $\phi_R\circ
  K_{\beta,\alpha}$.  By Propositions~\ref{prop:Denis-reps}
  and~\ref{prop:mdenis-reps} and Lemma~\ref{lem:DenisCommutes},
  $\lsup{\alpha}\Denis(-\PMC_1)^{\beta}\sos{\partial_R}{\cup}{F(\PMC_1^\beta)}\lsupv{\beta}\ol{\HD}^\beta
  \sos{F(\PMC_2^\beta)}{\cup}{\partial_L}\lsup{\beta}\MirrorDenis(-\PMC_2)^{\alpha}$
  also represents $-Y$, with the bordering $-\phi_L$ and $-\phi_R$.
  By the pairing theorem and
  Propositions~\ref{prop:beta-bimod-is-dual},~\ref{prop:CFAA-of-Denis}
  and~\ref{prop:CFAA-of-MirrorDenis},
  \[
  \CFAAa(\lsup{\alpha}\Denis(-\PMC_1)^{\beta}\sos{\partial_R}{\cup}{F(\PMC_1^\beta)}\lsupv{\beta}\ol{\HD}^\beta
  \sos{F(\PMC_2^\beta)}{\cup}{\partial_L}\lsup{\beta}\MirrorDenis(-\PMC_2)^{\alpha})\simeq \overline{\CFDDa(\hbox{$\displaystyle\lsupv{\alpha}\HD^\alpha$})}\DT\overline{\Alg_2'}\DT\Alg_1'.
  \]
  This implies Equation~\eqref{eq:boils-down}, and hence
  Equation~\eqref{eq:bimod-rev-1-1-precise}.

  To prove Equation~\eqref{eq:bimod-rev-1-2}, start by 
  tensoring both sides of Equation~\eqref{eq:bimod-rev-1-1-precise} with
  $\CFDDa(\Id_{\PMC_1})$, to obtain
  \begin{equation*}
  \Mor^{\Alg_1'}(\Alg_2'\DT\lsup{\Alg_1',\Alg_2'}\CFDDa(Y),\lsup{\Alg_1',\Alg_1}\CFDDa(\Id_{\PMC_1})) \\
  \simeq
  \lsup{\Alg_1}\CFDAa(-Y)_{\Alg_2'}.
  \end{equation*}
  Since tensoring over $\Alg_1'$ with $\CFAAa(\Id_{\PMC_1})$ gives an equivalence of
  categories which carries type~$D$ bordered invariants to type $A$ 
  bordered invariants (a bimodule analogue of
  Theorem~\ref{thm:DDId}, see
  \cite[Lemma~\ref*{LOT2:lemma:AADDquasi-equiv}]{LOT2}), this is the same as
  \begin{equation}
    \Mor_{\Alg_1}(\Alg_2'\DT\lsup{\Alg_2'}\CFDAa(Y)_{\Alg_1},\lsup{\Alg_1}[\Id]_{\Alg_1})
    \simeq
    \lsup{\Alg_1}\CFDAa(-Y)_{\Alg_2'}.\label{eq:bimod-rev-1-2-precise}
  \end{equation}
  After lowering the index $\Alg_1$, this is exactly
  Equation~\eqref{eq:bimod-rev-1-2}.

  Equation~\eqref{eq:bimod-rev-1-3} is immediate from
  Equations~\eqref{eq:bimod-rev-1-2} and~\eqref{eq:DA-bimod-3}.

  To prove Equation~\eqref{eq:bimod-rev-1-4}, observe that, by viewing
  left actions as right actions by the opposite algebra,
  Equation~\eqref{eq:bimod-rev-1-1} is equivalent to:
  \[
  \Mor_{\Alg_1}(\CFDDa(Y)_{\Alg_1,\Alg_2},\Alg_1)\simeq \lsub{\Alg_1,\Alg_2}\CFAAa(-Y).
  \]
  Applying Theorem~\ref{thm:DA-bimod} to both sides gives:
  \[
  \Mor_{\Alg_1}(\CFAAa(\tau_\bdy(Y))_{\Alg_1,\Alg_2},\Alg_1)\simeq
  \lsub{\Alg_1,\Alg_2}\CFDDa(\tau_\bdy^{-1}(-Y)).
  \]
  Recalling that $\tau_\bdy^{-1}(-Y)=-(\tau_\bdy(Y))$, this is just
  Equation~\eqref{eq:bimod-rev-1-4} with $Y$ replaced by
  $\tau_\bdy(Y)$.
\end{proof}

\begin{proof}[Proof of Theorem~\ref{thm:bimod-rev-2}]
  Equation~\eqref{eq:bimod-rev-2-1} is equivalent to the statement that 
  \begin{equation}\label{eq:bimod-pair-2-boils-down}
  \overline{\CFDDa(Y)}{}^{\Alg_1',\Alg_2'}\DT
    \Alg_1'\DT\Alg_2'\simeq
    \CFAAa(-\tau_\bdy^{-1}(Y))_{\Alg_1',\Alg_2'}.
  \end{equation}
  Fix a Heegaard diagram $\lsupv{\alpha}\HD^\alpha$ for $Y$ with $\partial_L\HD=\PMC_1$
  and $\partial_R\HD=\PMC_2$.  The pairing theorem and
  Propositions~\ref{prop:beta-is-dual} and~\ref{prop:CFAA-of-Denis}
  identify the left hand side of
  Equation~\eqref{eq:bimod-pair-2-boils-down} with
  $\CFAAa(\lsup{\alpha}\Denis(-\PMC_1)^{\beta}\sos{\partial_R}{\cup}{\partial_L}
  \lsupv{\beta}\ol{\HD}^\beta\sos{\partial_R}{\cup}{\partial_L}
  \lsupv{\beta}\Denis(-\PMC_2)^\alpha)$. 
  On the
  other hand, by Lemma~\ref{lem:DenisCommutes} and Corollary~\ref{cor:denis-glue-reps}, $\lsup{\alpha}\Denis(-\PMC_1)^{\beta}
  \sos{\partial_R}{\cup}{\partial_L}\lsupv{\beta}\ol{\HD}^\beta\sos{\partial_R}{\cup}{\partial_L}
  \lsup{\beta}\Denis(-\PMC_2)^\alpha$ is a bordered Heegaard diagram
  for
  $\tau_\bdy(-Y)=-\tau_\bdy^{-1}(Y)$. This proves
  Equation~\eqref{eq:bimod-pair-2-boils-down}.

  For Equation~\eqref{eq:bimod-rev-2-2}, tensor both sides of
  Equation~\eqref{eq:bimod-rev-2-1} with
  $\CFDDa(\Id_{\PMC_1})\otimes\CFDDa(\Id_{\PMC_2})$ to get
  \begin{equation*}
    \Mor^{\Alg_1'\otimes\Alg_2'}(\lsup{\Alg_1',\Alg_2'}\CFDDa(Y),
    \lsup{\Alg_1'\otimes\Alg_2',
      \Alg_1\otimes\Alg_2}(\CFDDa(\Id_{\PMC_1})\otimes\CFDDa(\Id_{\PMC_2})))
    \simeq
    \CFDDa(-\tau_{\bdy}^{-1}(Y)).
  \end{equation*}
  Since $\cdot \DT_{\Alg_1'\otimes\Alg_2'}\CFAAa(\Id_{\PMC_1})\otimes\CFAAa(\Id_{\PMC_2})$ is an equivalence
  of categories, we have
  \begin{equation*}
    \Mor_{\Alg_1\otimes\Alg_2}(\CFAAa(Y)_{\Alg_1,\Alg_2},\lsup{\Alg_1'\otimes\Alg_2'}[\Id]_{\Alg_1\otimes\Alg_2})
    \simeq
    \lsup{\Alg_1'\otimes\Alg_2'}\CFDDa(-\tau_{\bdy}^{-1}(Y)).
  \end{equation*}
  Lowering indices gives Equation~\eqref{eq:bimod-rev-2-2}.

  Alternatively, we can apply Theorem~\ref{thm:DA-bimod} to both
  sides of Equation~\eqref{eq:bimod-rev-2-1}. Specifically, viewing
  right actions by $A$ as
  left actions by $A^\op$, rewrite Equation~\eqref{eq:bimod-rev-2-1} as
  \[
    \Mor_{\Alg_1\otimes\Alg_2}(\CFDDa(Y,\spinc)_{\Alg_1,\Alg_2},
    \Alg_1\otimes\Alg_2)\simeq
    \lsub{\Alg_1,\Alg_2}\CFAAa(-\tau_\bdy^{-1}(Y),-\spinc).
  \]
  (For extra precision, we have added the $\SpinC$ structure to the
  equation.)
  Now, by Theorem~\ref{thm:DA-bimod}, the left hand side is
  \[
  \Mor_{\Alg_1\otimes\Alg_2}(\CFAAa(\tau_{\bdy}(Y)_{\Alg_1,\Alg_2},
  \overline{\spinc}),
    \Alg_1\otimes\Alg_2),
  \]
  while the right hand side is
  \[
  \lsub{\Alg_1,\Alg_2}\CFDDa(\tau_{\bdy}^{-1}(-\tau_\bdy^{-1}(Y)),-\overline{\spinc})=  \lsub{\Alg_1,\Alg_2}\CFDDa(-Y,-\overline{\spinc}).
  \]
  Replacing $Y$ by $\tau_{\bdy}^{-1}(Y)$ and conjugating the $\SpinC$ structure
  gives Equation~\eqref{eq:bimod-rev-2-2}.
\end{proof}

\begin{proof}[Proof of Corollaries~\ref{cor:bimod-hom-pair}
  and~\ref{cor:bimod-mod-hom-pair}]
  Equation~\eqref{bimod-cor-1-1} follows from
  Equation~\eqref{eq:bimod-rev-1-1} of Theorem~\ref{thm:bimod-rev-1}
  by taking $Y=Y_1$, tensoring both
  sides with $\CFDDa(Y_2)$ over
  $\Alg_1'$, and applying the pairing theorem. Equation~\eqref{bimod-cor-1-2} is obtained by
  applying Equation~\eqref{eq:bimod-rev-2-1} of
  Theorem~\ref{thm:bimod-rev-2} with $Y=F(\PMC_1)\times[0,1]$
  and then tensoring over $\Alg_1\otimes\Alg_1'$ with
  $\CFDDa(-Y_1)\otimes\CFDDa(Y_2)$.

  Equation~\eqref{eq:MorDDtoD} of
  Corollary~\ref{cor:bimod-mod-hom-pair} can be viewed as a special
  case of the first part of Corollary~\ref{cor:bimod-hom-pair}, in
  which one boundary component is empty.  Equation~\eqref{eq:D-to-DD}
  follows by taking $Y_1 = M_{\psi^{-1}} = -M_\psi$ for the first
  equation and taking $Y_1 = M_\psi$ and reversing the orientation on
  $Y_2$ for the second equation,
  in view of the fact that the action of $\psi$ on a bordered
  manifold~$Y$ is realized by gluing $M_\psi$ to~$Y$,
  see \cite[Lemma~\ref*{LOT2:lem:BorderedForDiffeo}]{LOT2}.
\end{proof}

\subsection{Dualizing bimodules}
\label{sec:dualizing-bimodules}
So far, we have used the type \AAm\ module associated to
$\Denis(\PMC)$. We next observe that the type \DD\ module associated
to $\MirrorDenis(-\PMC)$ gives a finite-dimensional model for the bar
resolution of $\Alg(\PMC)$.

Specifically, endow ${\overline \Alg}=\Hom_\Field(\Alg,\Field)$ with the
structure of a type \DD\ bimodule, as follows. Let $\Chord(\PMC)$
denote the set of chords in $\PMC$, i.e., arcs in $Z\setminus \{z\}$
connecting points in $\mathbf{a}$. Recall that to each chord
$\xi\in\Chord(\PMC)$ there is an associated algebra element $a(\xi)\in
\Alg=\Alg(\PMC)$.
The map
$$\delta^1\co \overline{\Alg} \rightarrow \Alg\otimes{\overline \Alg}\otimes \Alg$$
is defined by
\begin{equation}
  \label{eq:Differential}
  \delta^1(\phi)=1\otimes \bar{d}(\phi) \otimes 1 +
  \sum_{\xi\in\Chord(\PMC)} a(\xi)\otimes (a(\xi)\cdot \phi)\otimes 1
  +\sum_{\xi\in\Chord(\PMC)} 1\otimes (\phi\cdot a(\xi)) \otimes a(\xi).
\end{equation}
Here, $\bar{d}$ denotes the differential on $\overline{\Alg}$ (the
dual type $AA$ structure to the bimodule $\lsub{\Alg}\Alg_{\Alg}$, see
Definition~\ref{def:dual-bimod}) and
$a\cdot \phi$ and $\phi\cdot a$ denote the left and right actions of $\Alg$ on ${\overline \Alg}$.
We denote this type \DD\ bimodule $\SmallBar$. (We leave it to the reader to check that this satisfies the
structure equations for a type \DD\ bimodule.)

\begin{proposition}\label{prop:CFDD-of-Denis}
  The type \DD\ structure $\SmallBar$ is isomorphic to
  $\CFDDa(\MirrorDenis(-\PMC))$.
\end{proposition}
\begin{proof}
  The Heegaard diagram $\MirrorDenis(-\PMC)$ is a nice diagram
  (see~\cite{SarkarWang07:ComputingHFhat}). As in
  Proposition~\ref{prop:CFAA-of-MirrorDenis}, the rectangles supported
  in $\MirrorDenis(-\PMC)$ correspond to differentials in $\overline{\Alg(\PMC)}$.
  These give the terms of the form $1\otimes \bar{d}\phi\otimes 1$
  as in Equation~\eqref{eq:Differential}. 
  
  We must consider also rectangles which go out to the boundary. Those
  which go out to the $\alpha$-boundary give the terms of the form
  $\sum_{\xi\in\Chord(\PMC)} a(\xi)\otimes (a(\xi)\cdot \phi)\otimes 1$,
  while those which go out to the $\beta$-boundary give the terms of
  the form $\sum_{\xi\in\Chord(\PMC)} 1\otimes (\phi\cdot a(\xi)) \otimes
  a(\xi)$.
\end{proof}

\begin{corollary}
  The type \DD\ structure $\SmallBar$ is bounded.
\end{corollary}
\begin{proof}
  This follows immediately from the
  fact that $\MirrorDenis(-\PMC)$ is an admissible
  diagram. (Alternately, it is not hard to give a purely algebraic argument.)
\end{proof}

This gives a finite-dimensional model for the bar complex:
\begin{proposition}
  \label{prop:BoundedBar}
  Let $\Alg=\Alg(\PMC)$.  There are homotopy equivalences
  \begin{align}
    \SmallBar&\simeq \lsup{\Alg}\rBarop(\Alg)^{\Alg} \label{eq:BoundedBar}\\
    \SmallBar\DT \lsub{\Alg}{\Alg}_{\Alg} &\simeq \lsup{\Alg}[\Id]_{\Alg} \label{eq:BoundedId}\\
    \SmallBar\DT \lsub{\Alg}{\overline\Alg}_{\Alg} &
    \simeq \CFDAa(\tau_\bdy^{-1}\co F(\PMC)\to
    F(\PMC)) \label{eq:SerreFunctor}\\
    \oSmallBar\DT \lsub{\Alg}{\Alg}_{\Alg} &
    \simeq \CFDAa(\tau_\bdy\co F(\PMC)\to
    F(\PMC)) \label{eq:smallbar-dual}.
  \end{align}
\end{proposition}

\begin{proof}
  To prove Equation~\eqref{eq:BoundedId} observe that:
  \begin{align*}
    \SmallBar\DT \Alg &=\CFDDa(\lsupv{\alpha}\MirrorDenis(-\PMC)^{\beta}) 
    \DT \CFAAa(\lsupv{\beta}\Denis(\PMC)^{\alpha})\\
    &\simeq  \CFDAa(\lsupv{\alpha}\MirrorDenis(-\PMC)^{\beta}
    \cup\lsupv{\beta}\Denis(\PMC)^{\alpha}) \\
    &\simeq \CFDAa(\Id) \\
    &= \lsupv{\Alg}[\Id]_{\Alg}.
  \end{align*}
  The first equation follows from
  Propositions~\ref{prop:CFDD-of-Denis} and~\ref{prop:CFAA-of-Denis},
  the second follows from the
  pairing theorem, the third is a consequence of
  Propositions~\ref{prop:Denis-reps} and~\ref{prop:mdenis-reps}, and
  the last is Theorem~\ref{thm:id-is-id}, proved below (also proved
  in \cite[Theorem~\ref*{LOT2:thm:Id-is-Id}]{LOT2}). (Note the proof
  of Theorem~\ref{thm:id-is-id} does not rely on the current
  proposition.)

  Tensoring
  both sides of Equation~\eqref{eq:BoundedId} with
  $\lsup{\Alg}\rBarop(\Alg)^{\Alg}$,
  and using the fact that $\lsub{\Alg}\Alg_{\Alg}\DT
  \lsup{\Alg}\rBarop(\Alg)^{\Alg}\simeq \lsub{\Alg}[\Id]^{\Alg}$
  (essentially the
  statement that the algebra is quasi-isomorphic to its bar
  resolution), we obtain Equation~\eqref{eq:BoundedBar}.

  Equation~\eqref{eq:SerreFunctor} follows from the pairing theorem
  and Propositions~\ref{prop:CFAA-of-MirrorDenis}, \ref{prop:CFDD-of-Denis},
  and~\ref{prop:mdenis-reps}:
  \begin{align*}
    \SmallBar\DT \lsub{\Alg}{\overline\Alg}_{\Alg} &=
    \lsup{\Alg}(\CFDDa(\lsupv{\alpha}\MirrorDenis(-\PMC)^\beta)^\Alg\DT
    \lsub{\Alg}\CFAAa(\lsupv{\beta}\MirrorDenis(\PMC)^\alpha)_{\Alg}\\
    &\simeq
    \CFDAa(\lsupv{\alpha}\MirrorDenis(-\PMC)^\beta\cup\lsupv{\beta}\MirrorDenis(\PMC)^\alpha)\\
    &\simeq \CFDAa(\tau_\bdy^{-1}\co F(\PMC)\to F(\PMC)).
  \end{align*}

  To prove Equation~\eqref{eq:smallbar-dual}, first observe that
  $\oSmallBar=\lsup{\Alg}\CFDDa(\lsupv{\alpha}\Denis(-\PMC)^\beta)^\Alg$;
  this follows from Propositions~\ref{prop:minus-hd-bimod}
  and~\ref{prop:CFDD-of-Denis} (or a direct calculation). So, 
  \begin{align*}
    \oSmallBar\DT_{\Alg(\PMC)}
    \lsub{\Alg}{\Alg}_{\Alg}&=\lsup{\Alg}\CFDDa(\lsupv{\alpha}\Denis(-\PMC)^\beta)^\Alg\DT
    \lsub{\Alg}\CFAAa(\lsupv{\beta}\Denis(\PMC)^\alpha)_{\Alg}\\
    &\simeq \CFDAa(\tau_\bdy\co F(\PMC)\to F(\PMC)).\qedhere
  \end{align*}
\end{proof}

In a similar spirit, we can use the geometry of these pieces
to determine the explicit form for $\CFDDa(\Id)$. See also~\cite{LOT4}
for a different argument.

Consider the type \DD\ bimodule $\lsup{\Alg(\PMC)}K^{\Alg(\PMC)}$, defined
as follows. Let ${\mathbf s}$ be a subset of ${\mathbf p}/M$, and ${\mathbf t}$
denote its complement. The sets $\mathbf{s}$ and $\mathbf{t}$ have associated idempotents $I({\mathbf s})$
and $I({\mathbf t})$. We call such
pairs $(I(\mathbf{s}),I(\mathbf{t}))$ {\em complementary idempotents}.
Our type \DD\ bimodule
$\lsup{\Alg(\PMC)}K^{\Alg(\PMC)}=\bigoplus_{i=-k}^k\lsup{\Alg(\PMC,i)}K^{\Alg(\PMC,-i)}$has one
generator for each complementary
pair of idempotents. Let 
\[
\bOne=\sum_{(I,J)\text{ complementary}}\!\!\!\!\!\!
I\otimes_\Field J.
\]
Then the differential on
$\lsup{\Alg(\PMC)}K^{\Alg(\PMC)}$ is given by
\begin{equation}
  \label{eq:DDid}
  \delta^1 \bOne= 
  \sum_{\xi\in\Chord(\PMC)} a(\xi) \otimes \bOne \otimes
  a(\xi)=\sum_{\substack{\xi\in\Chord(\PMC)\\(I,J)\text{
          complementary}}}\!\!\!\!\!\!
  a(\xi)I \otimes (I\otimes_\Field J) \otimes Ja(\xi).
\end{equation}
As before, $a(\xi)$ denotes the algebra element in $\Alg(\PMC)$
associated to the chord~$\xi$. 

\begin{theorem}
  \label{thm:PreciseDD}\cite[Theorem~\ref*{LOT4:thm:DDforIdentity}]{LOT4}
  The bimodule $\lsup{\Alg(\PMC)}K^{\Alg(\PMC)}$ is isomorphic to
  $\CFDDa(\Id_\PMC)$, with the left action of $\Alg(-\PMC)$ viewed
  as a right action by $\Alg(\PMC)$.
\end{theorem}

\begin{proof}
  The type $\DD$ identity bimodule can be represented by an
  $\alpha$-$\alpha$ bordered Heegaard diagram $\HD$, as explained
  in~\cite[Section~\ref*{LOT2:sec:DiagramsForAutomorphisms}]{LOT2}. Instead
  of directly computing this type \DD\ structure, we will compute its
  modulification. To this end,
  form $-\overline{\HD}$, attach
  $\Denis(\PMC)$ and $\Denis(-\PMC)$ to its boundaries, and then
  destabilize the $k$ $\alpha$-circles in $-\overline{\HD}$.
  (This is the diagram
  illustrated on the left in Figure~\ref{fig:Identities}.)
  The type \AAm\ module associated to this diagram is exactly
  $\CFDDa(\Id)_{\Alg(\PMC),\Alg(-\PMC)}$. 

  Observe that the simplified diagram
  $\lsupv{\alpha}\Denis(\PMC)^{\beta} \cup \lsupv{\beta}\Denis(-\PMC)^{\alpha}$
  is a nice diagram, so the associated $\CFAAa$ is a bimodule with no
  higher action. Moreover, the bimodule structure on
  $\CFAAa(\lsupv{\alpha}\Denis(\PMC)^{\beta}\cup\lsupv{\beta}\Denis(-\PMC)^{\alpha})$ (viewed
  as having one left and one right action)
  is isomorphic to the bimodule structure on $\lsub{\Alg}\Alg_\Alg\DT
  \lsup{\Alg}K^\Alg\DT\lsub{\Alg}\Alg_\Alg$ as in
  Proposition~\ref{prop:CFAA-of-Denis}. 
  So, it only remains to
  determine the differential on 
  $\CFAAa(\lsupv{\alpha}\Denis(\PMC)^{\beta}\cup\lsupv{\beta}\Denis(-\PMC)^{\alpha})$,
  which in turn is given by embedded rectangles.
  There are rectangles supported in $\lsupv{\alpha}\Denis(\PMC)^{\beta}$, which
  correspond to differentials in $\Alg(\PMC)$
  (Proposition~\ref{prop:Denis-reps}), those supported in
  $\lsupv{\beta}\Denis(-\PMC)^{\alpha}$, which correspond to differentials in $\Alg(-\PMC)$,
  and rectangles which go between the two, which correspond to the
  differential $\delta^1$ from Equation~\eqref{eq:DDid}.
\end{proof}

\begin{remark}
  The alert reader may notice some redundancy in the proofs in this paper.  In
  fact, it is a consequence of Lemma~\ref{lem:QuasiInverse} that for a
  Koszul algebra, with dualizing bimodule $\lsup{A}K^B$ (which, in the
  present context, is $\lsup{A}\CFDDa(\Id)^A$, and in particular
  $A=B$), the bar resolution is homotopy equivalent to $\lsup{B}{\overline K}{}^A\DT
  \lsub{A}{\overline A}_A \DT \lsup{A}K^B$.
  But this latter bimodule
  is precisely the model $\SmallBar(\PMC)$ described in
  Equation~\eqref{eq:Differential}.  Thus,
  Proposition~\ref{prop:BoundedBar} can be viewed as a consequence of
  Lemma~\ref{lem:QuasiInverse} and Theorem~\ref{thm:PreciseDD}.
\end{remark}

\subsection{Algebraic consequences}
We start with a lemma regarding Serre functors which is probably
well-known in certain circles:
\begin{lemma}\label{lemma:SerreIs}
  For an augmented \dg algebra $\Alg$ and finite-dimensional right
  $\Ainf$-modules $M_\Alg$
  and~$N_\Alg$ over~$\Alg$,
  \begin{equation}
    \label{eq:SerreIs}
    \Mor_\Alg(M_\Alg, N_\Alg)^* \simeq \Mor_\Alg(N_\Alg, M_\Alg \DT \lsupv{\Alg}\rBarop^\Alg \DT \lsub{\Alg}\overline{\Alg}_\Alg)
  \end{equation}
  in a natural way.  That is, the Serre functor on $\ModCat_\Alg$ is
  given by tensoring on the right with $\rBarop(\Alg) \DT
  \overline{\Alg}$.
\end{lemma}
\begin{proof}
  The left hand side of Equation~\eqref{eq:SerreIs} is given by
  \begin{align*}
    \Mor_\Alg(M_\Alg,N_\Alg)^* &\cong (N_\Alg \DT
       \lsupv{\Alg}\overline{\rBarop}{}^\Alg \DT \lsub{\Alg}\overline{M})^*\\
     &\cong M_\Alg \DT \lsupv{\Alg}\rBarop^\Alg \DT
         \lsub{\Alg}\overline{N}\\
     &\simeq M_\Alg \DT \lsupv{\Alg}\rBarop^\Alg \DT
         \lsub{\Alg}\overline{\Alg}_\Alg \DT \lsupv{\Alg}\overline{\rBarop}{}^\Alg \DT
         \lsub{\Alg}\overline{N}\\
     &\cong \Mor_\Alg(N_\Alg, M_\Alg \DT \lsupv{\Alg}\rBarop^\Alg \DT \lsub{\Alg}\overline{\Alg}_\Alg),
  \end{align*}
  where we have used Proposition~\ref{prop:A-Mor-is-DT} (twice) and a
  dualized version of Equation~\eqref{eq:bar-resolution},
  $\overline{\Alg}\DT\lsupv{\Alg}\overline{\rBarop}{}^\Alg\simeq
  \lsub{\Alg}[\Id]^\Alg$.
\end{proof}

\begin{proof}[Proof of Theorem~\ref{thm:serre}]
  This is Lemma~\ref{lemma:SerreIs} plus the observation
  from Proposition~\ref{prop:BoundedBar} that
  \begin{align*}
    \lsupv{\Alg}\rBarop{}^\Alg \DT \lsub{\Alg}\overline{\Alg}_\Alg
      &\simeq {\SmallBar} \DT
        \lsub{\Alg}\overline{\Alg}_\Alg\\
      &\simeq \CFDAa(\tau_\bdy^{-1}\co F(\PMC)\to F(\PMC)).\qedhere
  \end{align*}
\end{proof}

Given two bimodules, we can also consider the complex of
$\Ainf$-bimodule morphisms between them.
The homology of this complex,
$H_*(\Mor_{\Alg,\Alg}(\lsub{\Alg}M_\Alg,\lsub{\Alg}N_\Alg))$, is also called
the \emph{Hochschild cohomology} of $M$ with $N$ and denoted
$\HH^*(M,N)$. 
The special case that $M=\Alg$ gives the Hochschild cohomology of $N$,
the derived functor associated to the functor of invariants in
$N$. With these observations, we are now ready to prove
Corollary~\ref{cor:Hochschild}.
\begin{proof}[Proof of Corollary~\ref{cor:Hochschild}]
  The Hochschild cohomology in question is the homology of the complex
  \begin{align*}
  \Mor_{\Alg,\Alg'}(\Alg,\Alg\DT\CFDAa(Y))&=(\lsupv{\Alg}\overline{\rBarop}{}^\Alg\DT\overline{\Alg}\DT\lsupv{\Alg}\overline{\rBarop}{}^\Alg)\DT_{\Alg\otimes\Alg'}(\Alg\DT\CFDAa(Y))\\
    &\simeq(\lsupv{\Alg}\rBarop^\Alg\DT\Alg\DT\lsupv{\Alg}\overline{\rBarop}{}^\Alg)\DT_{\Alg\otimes\Alg'}(\Alg\DT\CFDAa(Y)),
  \end{align*}
  using Proposition~\ref{prop:A-Mor-is-DT}),
  a dual version of Equation~\eqref{eq:bar-resolution}, and
  Equation~\eqref{eq:bar-resolution} itself.
  Rearranging the tensor products, we obtain
  \begin{align*}
    \Mor_{\Alg,\Alg'}(\Alg,\Alg\DT\CFDAa(Y))&\simeq(\lsupv{\Alg}\rBarop^\Alg)\DT_{\Alg\otimes\Alg'}(\Alg\DT\lsupv{\Alg}\overline{\rBarop}{}^\Alg\DT\Alg\DT\CFDAa(Y))\\
    &=
    \HC_*(\Alg\DT\lsupv{\Alg}\overline{\rBarop}{}^\Alg\DT\Alg\DT\CFDAa(Y))\\
    &\simeq \HC_*(\Alg\DT\CFDAa(\tau_\bdy)\DT \CFDAa(Y)) \\
    &\simeq \HC_*(\Alg\DT\CFDAa(\tau_\bdy(Y))) \\
    &\simeq \CFKa((\tau_\bdy(Y))^\circ,K).
  \end{align*}
  Here, $\HC_*$ denotes the Hochschild chain complex (whose homology
  is Hochschild homology).  The second line is the definition of
  $\HC_*$, the
  third uses Proposition~\ref{prop:BoundedBar}, the fourth uses the
  pairing theorem, and the last line uses
  \cite[Theorem~\ref*{LOT2:thm:DoublePairing}]{LOT2}. Taking homology
  gives the result.
\end{proof}

Finally, we give a simple proof that $\CFDAa(\Id)\simeq
\lsupv{\Alg}[\Id]_\Alg$.
\begin{proof}[Proof of Theorem~\ref{thm:id-is-id}]
  After raising the left index, we want to show
  $\CFDAa(\Id_\PMC)=\lsupv{\Alg(\PMC)}[\Id]_{\Alg(\PMC)}$. By the
  pairing theorem, $\CFAAa(\lsupv{\beta}\Denis(\PMC)^{\alpha})\DT
  \CFDAa(\Id_{\PMC})\simeq
  \CFAAa(\lsupv{\beta}\Denis(\PMC)^{\alpha})$. By Proposition~\ref{prop:CFAA-of-Denis}, this
  says that 
  \[
  \Alg(\PMC)\DT\CFDAa(\Id_{\PMC})\simeq \Alg(\PMC).
  \]
  But this implies
  $\CFDAa(\Id_{\PMC})\simeq\lsupv{\Alg(\PMC)}[\Id]_{\Alg(\PMC)}$, as desired.
\end{proof}
\begin{remark}
  For the proof of Theorem~\ref{thm:id-is-id}, it is irrelevant what
  $\Denis(\PMC)$ represents geometrically. All we need to know is
  that $\Alg(\PMC)$ is the bordered invariant of some Heegaard
  diagram. We also did not use any of the other theorems in this
  paper, many of which depend on the invertibility of $\CFDDa(\Id)$,
  which itself is equivalent, via the pairing theorem, to
  Theorem~\ref{thm:id-is-id}.
\end{remark}


\section{Gradings}
\label{sec:Gradings}

Bordered Heegaard Floer homology can be equipped with gradings, and
the pairing theorems described here are compatible with these gradings
in a natural way. We review these notions briefly.  For more details,
see \cite[Chapter~\ref*{LOT:chap:gradings}]{LOT1} and
\cite[Section~\ref*{LOT2:sec:algebras-gradings}]{LOT2}.

Given a pointed matched circle $\PMC$, there is a certain Heisenberg
group $\smallGroup(\PMC)$ equipped with a distinguished central
element $\lambda$, which has the property that $\Alg(\PMC)$ is graded
by $G=\smallGroup(\PMC)$.
It makes sense to talk about the category of differential graded modules
over this algebra. Objects in this category consist of pairs $(S,M)$,
where $S$ is a $G$-set, and $M$ is a module graded by $S$ in a way which is
compatible with the $G$-grading on~$\Alg$.

Given $G$-sets $S$ and $T$, we can form the space $\Hom(S,T)$, which
is orbit space of $S\times T$, divided out by its diagonal $G$ action.
(Note that this is not the same as the collection of $G$-set maps
$S\to T$.)
Now, given differential graded modules $(S,M)$ and $(T,N)$, the
morphism complex $\Mor((S,M),(T,N))$ is a $\ZZ$-set graded chain complex,
where the grading set is $\Hom(S,T)$, and the underlying chain complex
is as described earlier. Note that $\Hom(S,T)$ still admits an action
by $\ZZ$ (generated by the action of $\lambda$ on~$T$ or
$\lambda^{-1}$ on~$S$).  In
particular, the homology of the morphism space is also graded by
$\Hom(S,T)$.  (For generalities on these matters,
see~\cite[Section~\ref*{LOT2:sec:G-set-mod-cats}]{LOT2}.)

Bordered Heegaard Floer homology  modules are graded in the
above sense.  For example, given a $\PMC$-bordered three-manifold
$Y_1$ and a compatible bordered Heegaard diagram $\HD$, there is a
grading set $S=S(\HD)$ with the property that $\CFDa(Y_1)$ and $\CFAa(Y_1)$
are $S$-graded.

A graded version of Theorem~\ref{thm:hom-pair} (for $\CFDa$) reads as follows:

\begin{theorem}\label{thm:hom-pair-graded}
  Let $Y_1$ and $Y_2$ be bordered $3$-manifolds with $\bdy Y_1=\bdy
  Y_2=F(\PMC)$. Let $S_1$ and $S_2$ denote the grading sets for $Y_1$ and $Y_2$
  respectively. Then, there is an identification of the grading set for $\CFa(-Y_1\cup_{\bdy} Y_2)$
  with the $\ZZ$-set
  $\Hom(S_1,S_2)$, in such a manner that there is a graded isomorphism
  \[
  \HFa(-Y_1\cup_\bdy Y_2)\cong \Ext_{\Alg(-\PMC)}(\CFDa(Y_1),\CFDa(Y_2))
  \]
  which respects the identification of grading sets.
\end{theorem}

Theorem~\ref{thm:hom-pair-graded} follows from the proof of
Theorem~\ref{thm:hom-pair}, with two additional observations. The
first is that the grading set for $\Denis(\PMC)$ is naturally
identified with $\smallGroup$ in such a manner that
Proposition~\ref{prop:CFAA-of-Denis} holds in its graded form
(i.e., $\CFAAa(\Denis(\PMC))$ is isomorphic to $\Alg(\PMC)$ as a
$\smallGroup(\PMC)$-graded bimodule); and the second observation is
that the traditional pairing
theorem~\cite[Theorem~\ref*{LOT2:thm:GenReparameterization}]{LOT2}
used in establishing Theorem~\ref{thm:or-rev} also holds
in a graded form,
see~\cite[Theorem~\ref*{LOT2:thm:GradedPairing}]{LOT2}).

If we keep track of $\SpinC$-structures, the isomorphism in
Theorem~\ref{thm:hom-pair-graded} is given by
\[
\bigoplus_{\substack{\spinc\in\SpinC(-Y_1\cup_\bdy Y_2)\\
    \spinc|_{-Y_1}=-\spinc_1,\ \spinc|_{Y_2}=\spinc_2}}\!\!\!\!
\HFa(-Y_1\cup_\bdy Y_2,\spinc)\cong \Ext_{\Alg(-\PMC)}(\CFDa(Y_1,\spinc_1),\CFDa(Y_2,\spinc_2)).
\]
Similarly, a version of Theorem~\ref{thm:or-rev} keeping track of
$\SpinC$-structures is:
\begin{align*}
  \Mor_{\Alg(-\PMC)}(\lsub{\Alg(-\PMC)}\CFDa(Y,\spinc),\Alg(-\PMC))&\simeq \CFAa(-Y,-\spinc)_{\Alg(-\PMC)}\\
  \Mor_{\Alg(\PMC)}(\CFAa(Y,\spinc)_{\Alg(\PMC)},\Alg(\PMC))&\simeq \lsub{\Alg(\PMC)}\CFDa(-Y,-\spinc).
\end{align*}
(Compare Proposition~\ref{prop:beta-is-dual}.)

Gradings can also be added in a straightforward way for
Theorem~\ref{thm:hom-pair} for $\CFAa$, and to
the rest of the theorems from the introduction. In particular,
the gradings in Corollary~\ref{cor:Hochschild} are obtained from a straightforward
adaptation of~\cite[Theorem~\ref*{LOT2:thm:DoublePairing}]{LOT2}.


\section{Examples}\label{sec:Examples}
In this section, we compute a few simple examples with the $\Hom$
pairing theorem, and compare them with the results of the original,
tensor product pairing theorem.

\subsection{Review of the torus algebra}
For simplicity, all of our examples will have torus boundary, and we
will work in the central $\SpinC$-structure, so we start by reviewing
the algebra $\Alg(T^2)=\Alg(T^2,0)$ associated to the (unique) genus
$1$ pointed matched circle. The algebra $\Alg(T^2)$ has an
$\FF_2$-basis with $8$ elements $\iota_0$, $\iota_1$, $\rho_1$,
$\rho_2$, $\rho_3$, $\rho_{12}$, $\rho_{23}$ and $\rho_{123}$. The
elements $\iota_0$ and $\iota_1$ are orthogonal idempotents. The other
relations on the algebra are:
\begin{align*}
\iota_0\rho_1\iota_1&=\rho_1 & \iota_1\rho_2\iota_0&=\rho_2&
\iota_0\rho_3\iota_1&=\rho_3\\
\rho_1\rho_2&=\rho_{12} &\rho_2\rho_3&=\rho_{23} &
\rho_1\rho_{23}&=\rho_{123}\\
\rho_{12}\rho_3&=\rho_{123} &\rho_3\rho_2&=0& \rho_2\rho_1&=0.
\end{align*}
(See also~\cite[Section~\ref*{LOT1:sec:torus-algebra}]{LOT1}.)

\subsection{\textalt{$\Hom$}{Hom} pairing theorem for some solid tori}
\begin{figure}
  \centering
  \includegraphics{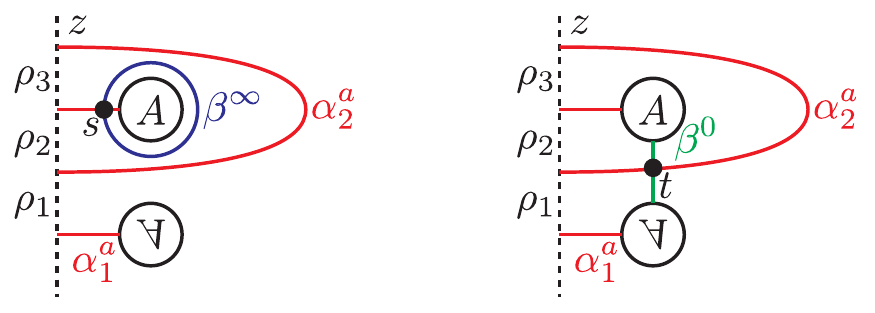}
  \caption{\textbf{Standard bordered Heegaard diagrams for the
      $\infty$- and $0$-framed solid tori.} Left: the diagram
    $\HD_\infty$ for the $\infty$-framed solid torus. Right: the
    diagram $\HD_0$ for the $0$-framed solid torus.}
  \label{fig:0-infty-tori}
\end{figure}

We start by gluing together some solid tori. Consider the standard
diagrams $\HD_\infty$ and $\HD_0$ for the $\infty$- and $0$-framed
solid tori, shown in Figure~\ref{fig:0-infty-tori}
(compare~\cite[Section~\ref*{LOT1:sec:surg-exact-triangle}]{LOT1}). 

The module $\CFDa(\HD_\infty)$ has a single generator $s$ with
$\iota_1 s=s$ and differential
\[
\bdy(s)=\rho_{23}s.
\]

(We have used module notation for $\CFDa(\HD_\infty)$, so that this is
another way of saying $\delta^1(s) = \rho_{23} \otimes s$.)
The module $\CFDa(\HD_0)$ has a single generator $t$ with $\iota_0t=t$
and differential
\[
\bdy(t)=\rho_{12}t.
\]

So, for instance, the chain complex of homomorphisms 
$\Mor(\CFDa(\HD_\infty),\CFDa(\HD_\infty))$ is generated by elements
$f$ and $g$ where
\begin{align*}
  f(s)&= s\\
  g(s)&=\rho_{23}s.
\end{align*}
The differential of $f$ is
\begin{equation}\label{eq:diff-on-S2S1}
(\bdy f)(s)=\bdy(f(s))+f(\bdy s)=\bdy s+f(\rho_{23}s)=\rho_{23}s+\rho_{23}s=g(s)+g(s)=0.
\end{equation}
Similarly, $\bdy g=0$. So, $\Ext(\CFDa(\HD_\infty),\CFDa(\HD_\infty))$
is two-dimensional. This is
consistent with Theorem~\ref{thm:hom-pair}, since $-\HD_\infty
\cup_{\bdy} \HD_\infty$ represents $S^1 \times S^2$ and
 $\HFa(S^1\times
S^2)\cong \FF_2\oplus\FF_2$.

It is clear from Equation~\eqref{eq:diff-on-S2S1} that $\gr(g)$ is $1$
lower than $\gr(f)$. To illustrate the behavior of the gradings in the
$\Hom$ pairing theorem, we compute this directly.  We use the notation
from~\cite[Section~\ref*{LOT1:sec:torus-algebra}]{LOT1} for the
grading groups of the torus, and the grading refinement there, taking values
in a somewhat larger group $G \superset \smallGroup(T^2)$.
We have:
\[
  \DBigGrSet(\HD_\infty)=\bigGroup(T^2)/\langle(-1/2;0,-1,-1)\rangle
  \qquad\qquad 
  \DSmallGrSet(\HD_\infty)=G/\langle(-1/2;0,-1)\rangle.
\]
Declaring arbitrarily that $\gr(s)=[(0;0,0)]$, it follows that
$\gr(\rho_{23}s)=\gr(\rho_{23})\gr(s)=[(-1/2;0,1)]$. So,
\begin{align*}
  \gr(f)&=[(0;0,0)]\times[(0;0,0)]\in \bigl(\DSmallGrSet(\HD_\infty)\times
  \DSmallGrSet(\HD_\infty)\bigr)/G\\
  \gr(g)&=[(0;0,0)]\times[(-1/2;0,1)]=[(0;0,0)]\times[(-1;0,0)]\in \bigl(\DSmallGrSet(\HD_\infty)\times
  \DSmallGrSet(\HD_\infty)\bigr)/G.
\end{align*}
In particular, $\gr(g)=\lambda^{-1}\gr(f)$, as claimed.

As another simple example, $\Mor(\CFDa(\HD_0),\CFDa(\HD_\infty))$ is
generated by the three maps
  $t\mapsto \rho_1 s$, 
  $t\mapsto \rho_3 s$, and
  $t\mapsto \rho_{123}s$
with differentials
\begin{align*}
  \bdy(t\mapsto \rho_1 s)&=(t\mapsto\rho_{123}s)\\
  \bdy(t\mapsto \rho_3 s)&=(t\mapsto\rho_{123}s),
\end{align*}
so $\Ext(\CFDa(\HD_0),\CFDa(\HD_\infty))$ is $1$-dimensional, in
agreement with $Y(-\HD_0 \cup_\bdy \HD_\infty) = S^3$.

Similar computations show:
\begin{align*}
\Ext(\CFDa(\HD_0),\CFDa(\HD_0))&\cong \FF_2\oplus \FF_2 \\
\Ext(\CFDa(\HD_\infty),\CFDa(\HD_0))&\cong\FF_2.
\end{align*}

\subsection{\textalt{$\CFAa$}{CFA\textasciicircum} is \textalt{$\CFDa$}{CFD\textasciicircum}: an example}
Next, we illustrate Theorem~\ref{thm:D-is-A} for the $0$-framed solid
torus $\HD_0$ discussed above. The module
$\lsub{\Alg(T^2)}\CFDa(\HD_0)$ has elements $t$, $\rho_2t$ and
$\rho_{12}t$. The differential of $t$ is $\rho_{12}t$, so
$H_*(\CFDa(\HD_0))=\Field\langle \rho_2t\rangle$. This agrees with the
rank of $\CFAa(\HD_0)$ (on which the differential is trivial).

Moreover, we can reconstruct the $\Ainf$-module structure on
$H_*(\CFDa(\HD_0))$. We record the \dg module structure on
$\CFDa(\HD_0)$ as:
\[
\begin{tikzpicture}
  \node at (0,0) (t) {$t$};
  \node at (2,0) (rho2t) {$\rho_2t$};
  \node at (2,-1) (rho12t) {$\rho_{12}t$};
  \draw[->] (t) to node[above]{\lab{\rho_2}} (rho2t);
  \draw[->] (rho2t) to node[right]{\lab{\rho_1}} (rho12t);
  \draw[->] (t) to node[below]{\lab{1+\rho_{12}}} (rho12t);
\end{tikzpicture}
\]
Let $x=\rho_2t$. Cancelling the differential from $t$ to $\rho_{12}t$ gives us
$\Ainf$-structure on $\Field\langle x\rangle$ given by
the expansion of $\rho_2,(1+\rho_{12})^{-1},\rho_1$. That is, in $H_*(\CFDa(\HD_0))$,
$m_3(\rho_2,\rho_1,x)=x$, $m_4(\rho_2,\rho_{12},\rho_1,x)=x$,
$m_5(\rho_2,\rho_{12},\rho_{12},\rho_1,x)=x$, and so on. By contrast,
the (right) $\Ainf$-module structure on $\CFAa(\HD_0)$ is given by
$m_3(x,\rho_3,\rho_2)=x$, $m_4(x,\rho_3,\rho_{23},\rho_2)=x$, and so
on. Under the identification of left modules over $\Alg(T^2)$ with
right modules over $\Alg(-T^2)$, these two $\Ainf$-structures agree.

\subsection{Hochschild cohomology of $\Alg(T^2)$}
In this section, we compute the Hochschild cohomology of $\Alg(T^2)$;
in light of Corollary~\ref{cor:Hochschild}, this is the same as
computing $\HFKa$ of $-1$ surgery on the Borromean knot.

Let $\Id$ denote the identity map of the torus.
Since tensoring with $\CFDDa(\Id)$ gives an equivalence of categories,
it is equivalent to compute 
\[
\HH^*(\CFDDa(\Id),\CFDDa(\Id))=H_*\bigl(\Mor_{\Alg(T^2),\Alg'(-T^2)}(\CFDDa(\Id),\CFDDa(\Id))\bigr),
\]
and this is what we will do.

Recall from Theorem~\ref{thm:PreciseDD} that
$\lsup{\Alg(\PMC,i)}\CFDDa(\Id)^{\Alg(\PMC,-i)}$ is generated by pairs
of complementary idempotents. In the case under consideration,
$\PMC=T^2$ and $i=0$, so there are two pairs of complementary
idempotents $\iota_0\otimes\iota_1$ and $\iota_1\otimes\iota_0$. The
differential on $\CFDDa(\Id)$ is given by
\begin{align*}
  \bdy(\iota_0\otimes\iota_1)&=\rho_1\otimes(\iota_1\otimes\iota_0)\otimes\rho_1+\rho_3\otimes(\iota_1\otimes\iota_0)\otimes\rho_3+\rho_{123}\otimes(\iota_1\otimes\iota_0)\otimes\rho_{123}\\
  \bdy(\iota_1\otimes\iota_0)&=\rho_2\otimes(\iota_0\otimes\iota_1)\otimes\rho_2.
\end{align*}

A basis of $\Mor_{\Alg(T^2),\Alg'(-T^2)}(\CFDDa(\Id),\CFDDa(\Id))$ is
given by the maps $f$ sending $\iota_i\otimes \iota_{1-i}$ to
$\rho\otimes (\iota_j\otimes \iota_{1-j})\otimes \sigma$ and
$\iota_{1-i}\otimes\iota_i$ to zero.  Here $\rho$ and $\sigma$ are
chords in $\Alg(T^2)$ respectively, with $\iota_i\rho\iota_j=\rho$ and
$\iota_{1-j}\sigma\iota_{1-i}=\sigma$.  Without loss of information,
we will denote the map $f$ as $\langle \rho\otimes\sigma\rangle$.
Then the generators of
$\Mor_{\Alg(T^2),\Alg'(-T^2)}(\CFDDa(\Id),\CFDDa(\Id))$ are:
\[
\begin{matrix}
  \langle \iota_1\otimes\iota_0\rangle &  
  \langle \rho_{23}\otimes\iota_0\rangle & 
  \langle \iota_1\otimes\rho_{12}\rangle & 
  \langle \rho_2\otimes\rho_2\rangle & 
  \langle \rho_{23}\otimes\rho_{12}\rangle \\ 
  \langle \iota_0\otimes\iota_1\rangle &   
  \langle \rho_{12}\otimes\iota_1\rangle &    
  \langle \iota_0\otimes\rho_{23}\rangle &    
  \langle \rho_1\otimes\rho_3\rangle & 
  \langle \rho_{12}\otimes\rho_{23}\rangle &   
  \langle \rho_1\otimes\rho_1\rangle &   
  \langle \rho_1\otimes\rho_{123}\rangle \\   
  \langle \rho_{123}\otimes\rho_1\rangle & 
  \langle \rho_{123}\otimes\rho_{123}\rangle &  
  \langle \rho_{123}\otimes\rho_3\rangle &   
  \langle \rho_3\otimes\rho_1\rangle &   
  \langle \rho_3\otimes\rho_{123}\rangle &   
  \langle \rho_3\otimes\rho_3\rangle. & 
\end{matrix}
\]
(The maps in the first row send $\iota_1\otimes\iota_0$ to the specified
element; the maps in the second and third rows send $\iota_0\otimes\iota_1$ to the
specified element.)

The nontrivial differentials are given by
\begin{align*}
  \bdy\langle \iota_1\otimes\iota_0\rangle &=  \langle
  \rho_2\otimes\rho_2\rangle+  \langle \rho_1\otimes\rho_1\rangle+
  \langle \rho_{123}\otimes\rho_{123}\rangle+  \langle
  \rho_3\otimes\rho_3\rangle\\ 
  \bdy   \langle \rho_{23}\otimes\iota_0\rangle&=  \langle
  \rho_{123}\otimes\rho_1\rangle\\ 
  \bdy  \langle \iota_1\otimes\rho_{12}\rangle &=  \langle
  \rho_3\otimes\rho_{123}\rangle\\
  \bdy  \langle \iota_0\otimes\iota_1\rangle &=   \langle \rho_2\otimes\rho_2\rangle+ \langle \rho_1\otimes\rho_1\rangle+
  \langle \rho_{123}\otimes\rho_{123}\rangle+  \langle
  \rho_3\otimes\rho_3\rangle\\
  \bdy  \langle \rho_{12}\otimes\iota_1\rangle&=   \langle
  \rho_{123}\otimes\rho_3\rangle\\
  \bdy   \langle \iota_0\otimes\rho_{23}\rangle  &=  \langle
  \rho_1\otimes\rho_{123}\rangle\\
  \bdy   \langle \rho_1\otimes\rho_3\rangle&=  \langle \rho_{12}\otimes\rho_{23}\rangle.
\end{align*}

A straightforward computation shows that the homology is
$4$-dimensional, generated by 
\begin{gather*}
\bOne=\langle \iota_1\otimes\iota_0\rangle +   \langle
\iota_0\otimes\iota_1\rangle \\
w = \langle \rho_2\otimes\rho_2\rangle \qquad  
x = \langle \rho_1\otimes\rho_1\rangle \qquad  
y =\langle \rho_3\otimes\rho_3\rangle \qquad 
z =\langle \rho_{123}\otimes\rho_{123}\rangle 
\end{gather*}
with the relation
\[
w+x+y+z=0.
\]
The element $\langle \iota_1\otimes\iota_0\rangle + \langle
\iota_0\otimes\iota_1\rangle$ acts as a unit for the multiplication on
$\HH^*(\CFDDa(\Id),\penalty600\CFDDa(\Id))$ and all other products vanish. The
grading of $\bOne$ is one lower than the grading of $w$, $x$, $y$ and
$z$ (with the convention that the grading on $\HH^*$ is of
cohomological type).

\begin{remark}
Recall~\cite{Gerstenhaber63} that the the Hochschild cohomology of an
algebra also inherits a Lie bracket, called the  ``Gerstenhaber bracket''.
Whereas the algebra structure on the Hochschild cohomology is convenient to describe in terms
of automorphisms of $\CFDDa(\Id)$ as above, the Gerstenhaber bracket is not 
transparent from this perspective. Nonetheless, with a little
more work, one can identify the generators of the homology of the
standard Hochschild cochain complex as
\[
\iota_0[]+\iota_1[],
\rho_1[\rho_1^*]+\rho_{123}[\rho_{123}^*]+\rho_{12}[\rho_{12}^*],
\rho_3[\rho_3^*]+\rho_{123}[\rho_{123}^*]+\rho_{23}[\rho_{23}^*], \rho_{123}[\rho_3^*|\rho_2^*|\rho_1^*].
\]
From this, it is straightforward to verify that the Gerstenhaber
bracket vanishes.
\end{remark}

\subsection{\textalt{$\CFAAa$}{CFAA\textasciicircum} is \textalt{$\CFDDa$}{CFDD\textasciicircum} with
  a negative boundary Dehn twist: an example}
We illustrate Theorem~\ref{thm:DA-bimod}
by verifying that for the standard Heegaard diagram for $\Id_{T^2}$,
the rank of $H_*(\CFDDa(\Id_{T^2}))$ agrees with the rank of
$H_*(\CFAAa(\tau_\bdy))$.  Of course, the theorem asserts much more
than this: the bimodule structures agree. Even in this simple case,
computing the $\Ainf$-bimodule structure on the homology is somewhat tedious,
and we will not record the details here.

From Theorem~\ref{thm:PreciseDD}, or alternatively~\cite[Section~\ref*{LOT2:subsec:AAId1}]{LOT2} or~\cite{LOT4}, as a type
\DD\ structure, $\lsup{\Alg(T^2),\Alg(-T^2)}\CFDDa(\Id)$ has two
generators $x$ and $y$, with
\begin{equation}\label{eq:DD-identity}
\begin{aligned}
  \bdy
  x&=(\rho_1\sigma_3+\rho_3\sigma_1+\rho_{123}\sigma_{123})\otimes y\\
  \bdy y&=(\rho_2\sigma_2)\otimes x.
\end{aligned}
\end{equation}
(Here, we use $\rho$'s to denote elements of $\Alg(T^2)$ and
$\sigma$'s to denote elements of $\Alg(-T^2)$.) Expanding this, as a
bimodule (rather than type \DD\ structure), $\lsub{\Alg(T^2),\Alg(-T^2)}\CFDDa(\Id)$ has $34$
generators, with differentials as shown in Figure~\ref{fig:cfdd-id}.
The homology is $16$-dimensional.

\begin{figure}
\[
\begin{tikzpicture}
  \node at (0,1) (x) {$x$};
  \node at (0,-1) (s2x) {$\sigma_2x$};
  \node at (0,-3) (s12x) {$\sigma_{12}x$};
  \node at (2,1) (r2x) {$\rho_2x$};
  \node at (2,-2) (r2s2x) {$\rho_2\sigma_2x$};
  \node at (2,-4) (r2s12x) {$\rho_2\sigma_{12}x$};
  \node at (4,1) (r12x) {$\rho_{12}x$};
  \node at (4,-2) (r12s2x) {$\rho_{12}\sigma_2x$};
  \node at (4,-4) (r12s12x) {$\rho_{12}\sigma_{12}x$};
  \node at (5,2) (y) {$y$};
  \node at (5,0) (s1y) {$\sigma_1y$};
  \node at (6,-2) (s3y) {$\sigma_3y$};
  \node at (6,-4) (s23y) {$\sigma_{23}y$};
  \node at (6,-6) (s123y) {$\sigma_{123}y$};
  \node at (7,2) (r1y) {$\rho_1y$};
  \node at (7,0) (r1s1y) {$\rho_1\sigma_1y$};
  \node at (8,-2) (r1s3y) {$\rho_1\sigma_3y$};
  \node at (8,-4) (r1s23y) {$\rho_1\sigma_{23}y$};
  \node at (8,-6) (r1s123y) {$\rho_1\sigma_{123}y$};
  \node at (9,2) (r3y) {$\rho_3y$};
  \node at (10,0) (r3s1y) {$\rho_3\sigma_1y$};
  \node at (10,-2) (r3s3y) {$\rho_3\sigma_3y$};
  \node at (10,-4) (r3s23y) {$\rho_3\sigma_{23}y$};
  \node at (10,-6) (r3s123y) {$\rho_3\sigma_{123}y$};
  \node at (11,2) (r23y) {$\rho_{23}y$};
  \node at (12,0) (r23s1y) {$\rho_{23}\sigma_1y$};
  \node at (12,-2) (r23s3y) {$\rho_{23}\sigma_3y$};
  \node at (12,-4) (r23s23y) {$\rho_{23}\sigma_{23}y$};
  \node at (12,-6) (r23s123y) {$\rho_{23}\sigma_{123}y$};
  \node at (13,2) (r123y) {$\rho_{123}y$};
  \node at (14,0) (r123s1y) {$\rho_{123}\sigma_1y$};
  \node at (14,-2) (r123s3y) {$\rho_{123}\sigma_3y$};
  \node at (14,-4) (r123s23y) {$\rho_{123}\sigma_{23}y$};
  \node at (14,-6) (r123s123y) {$\rho_{123}\sigma_{123}y$};
  \draw[->, color=red] (x) to (r1s3y);
  \draw[->, color=red] (s2x) to (r1s23y);
  \draw[->, color=red] (s12x) to (r1s123y);
  \draw[->, color=blue] (x) to (r3s1y);
  \draw[->, color=blue] (r2x) to (r23s1y);
  \draw[->, color=blue] (r12x) to (r123s1y);
  \draw[->, color=purple] (x) to (r123s123y);
  \draw[->, color=brown] (y) to (r2s2x);
  \draw[->, color=brown] (r1y) to (r12s2x);
  \draw[->, color=brown] (s1y) to (r2s12x);
  \draw[->, color=brown] (r1s1y) to (r12s12x);
\end{tikzpicture}
\]
  \caption{\textbf{The \dg bimodule for $\CFDDa(\Id)$.}  Arrows coming
    from the same term in $\bdy x$ or $\bdy y$ are
    parallel.}
\label{fig:cfdd-id}
\end{figure}

\begin{figure}
  \[
  \includegraphics{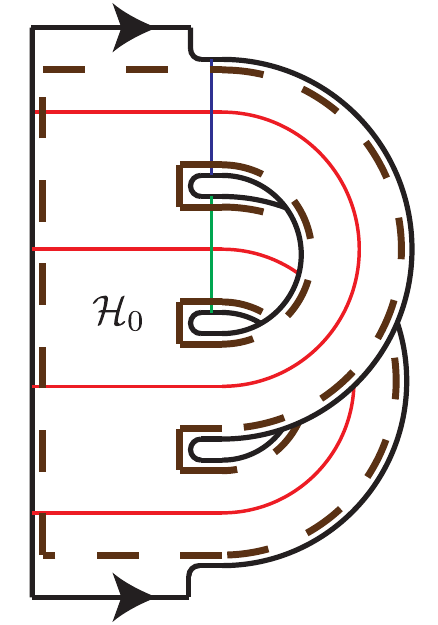} \qquad\qquad
  \includegraphics{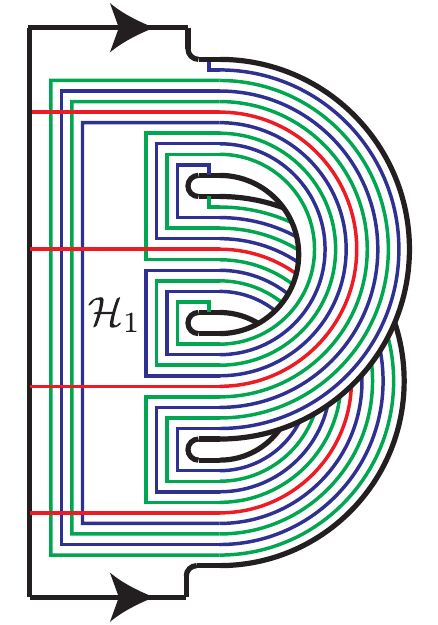}
  \]
  \caption{\textbf{Standard Heegaard diagram for the negative boundary
      Dehn twist.} On the left is half of a Heegaard diagram for the
    identity homeomorphism; if this surface is $\HD_0$ then the result
    $\HD_0 \sos{\bdy_R}{\cup}{\bdy_L} -\HD_0$ of gluing $\HD_0$ to its
    mirror image (along $\bdy_R$ and $\bdy_L$, respectively) is a
    Heegaard diagram for the identity map. Applying a negative Dehn
    twist to the $\beta$ curves around the dashed curve (and isotoping
    away a bigon) gives the diagram $\HD_1$ on the right. The result
    $\HD_1\sos{\bdy_R}{\cup}{\bdy_L} -\HD_0$ of gluing $\HD_1$ to the
    mirror of $\HD_0$ gives the standard Heegaard diagram for the
    positive boundary Dehn twist.}
  \label{fig:DehnTwistDiagram}
\end{figure}

Figure~\ref{fig:DehnTwistDiagram} shows part of the standard Heegaard
diagram
for the negative boundary Dehn twist. Inspecting the diagram, we see
that there
are $16$ generators in the middle $\SpinC$ structure
and no provincial domains, hence no differential
on $\CFAAa$. Thus, the rank of $\CFAAa(\tau_\bdy)$ agrees with the
rank of $H_*(\CFDDa(\Id))$, as claimed.

\subsection{A surgery on the trefoil}
We  conclude by computing $\HFa$ of the three-manifold $Y$ obtained
as $-2$ Dehn filling of the
left-handed trefoil complement, using Theorem~\ref{thm:hom-pair}. It follows from~\cite[Theorem~\ref*{LOT1:thm:HFKtoHFDframed}]{LOT1} that the
type $D$ structure $\CFDa(T)$ associated to the $-2$-framed trefoil
complement is:
\[
\begin{tikzpicture}
  \node at (0,0) (x3) {$x_3$};
  \node at (0,-1) (y2) {$y_2$};
  \node at (0,-2) (x2) {$x_2$};
  \node at (2,-2) (y1) {$y_1$};
  \node at (4,-2) (x1) {$x_1.$};
  \draw[->] (x3) to node[left]{$\rho_1$} (y2);
  \draw[->] (x2) to node[left]{$\rho_{123}$} (y2);
  \draw[->] (y1) to node[below]{$\rho_2$} (x2);
  \draw[->] (x1) to node[below]{$\rho_3$} (y1);
  \draw[->, bend right=30] (x1) to node[above right,inner sep=2pt]{$\rho_{12}$} (x3);
\end{tikzpicture}
\]
(compare~\cite[Figure~\ref*{LOT1:fig:CableTrefoil}]{LOT1}). Here, the $x_i$ are generators with
$\iota_{0}x_i=x_i$ and the $y_i$ are generators with
$\iota_{1}y_i=y_i$. The arrow from $x_1$ to $y_1$, say, denotes the
fact that $\rho_3y_1$ occurs in $\bdy(x_1)$.

We
compute $\Mor(\CFDa(\HD_0),\CFDa(T))$. Each $x_i$ is replaced by the
maps $t\mapsto x_i$ and $t\mapsto \rho_{12}x_i$.  Each $y_i$ is
replaced by the three maps $t\mapsto \rho_1 y_i$, $t\mapsto \rho_3
y_i$, and $t\mapsto \rho_{123}y_i$.  The differentials are:
\[
\begin{tikzpicture}
  \node at (0,0) (x3) {$t\mapsto x_3$};
  \node at (0,-1) (12x3) {$t\mapsto \rho_{12}x_3$};
  \node at (0,-3) (1y2) {$t\mapsto \rho_1y_2$};
  \node at (0,-4) (3y2) {$t\mapsto \rho_3y_2$};
  \node at (0,-5) (123y2) {$t\mapsto \rho_{123}y_2$};
  \node at (0,-7) (x2) {$t\mapsto x_2$};
  \node at (0,-8) (12x2) {$t\mapsto \rho_{12}x_2$};
  \node at (3,-6) (1y1) {$t\mapsto \rho_1y_1$};
  \node at (3,-7) (3y1) {$t\mapsto \rho_3y_1$};
  \node at (3,-8) (123y1) {$t\mapsto \rho_{123}y_1$};
  \node at (6,-7) (x1) {$t\mapsto x_1$};
  \node at (6,-8) (12x1) {$t\mapsto \rho_{12}x_1.$};
  \draw[->] (x3) to (12x3);
  \draw[->] (x2) to (12x2);
  \draw[->] (x1) to (12x1);
  \draw[->] (3y2) to (123y2);
  \draw[->] (3y1) to (123y1);
  \draw[->, bend right=60] (x3) to (1y2);
  \draw[->] (x2) to (123y2);
  \draw[->] (1y1) to (12x2);
  \draw[->] (x1) to (3y1);
  \draw[->] (12x1) to (123y1);
  \draw[->] (x1) to (12x3);
\end{tikzpicture}
\]
The homology $\Ext(\CFDa(\HD_0),\CFDa(T))$ is $2$-dimensional,
generated by $t\mapsto\rho_1 y_2$ and $t\mapsto
(\rho_1y_1+x_2+\rho_3y_2).$

In this example, after making a choice, the grading set 
is identified with a double coset space of
$\smallGroup=\smallGroup(T)$.  Projecting onto the homological
component of $\smallGroup$, the double coset space maps to
$H_1(Y)\cong \ZZ/2\ZZ$. The gradings of our two generators project to
the two different elements, showing that $\HFa(Y;\spinc)$ is $\ZZ/2\ZZ$ for
each of the two $\SpinC$ structures $\spinc$.

This result can be compared with results from~\cite{SomePlumb}
or~\cite{IntSurg}, which give alternate methods for computing the
Heegaard Floer homology groups of $-2$-surgery on the left-handed
trefoil.

The reader might find it disappointing that these calculations do not
distinguish $-2$ surgery on the trefoil from that on the unknot: we
gave an example which is an \emph{L-space}. It is easy to distinguish
these two three-manifolds, however, via the $\QQ$ grading on Heegaard
Floer homology from~\cite{AbsGraded}; see~\cite{LOTRelGrad} for a
discussion on how to extract this information  using bordered Floer homology.


\section{Relation to Koszul duality}
\label{sec:koszul}

We will now justify the reference to Koszul duality in  the \DD\ bimodule described in Equation~\eqref{eq:DDid}.  We also use Koszul
duality to explain a symmetry in the homology of the algebra of $\Alg(\PMC)$
which was observed in~\cite{LOT2}.
We start with some review.

\subsection{Formalities on Koszul duality}
\label{sec:koszul-frame}

Fix a ground ring $\Ground = \bigoplus \FF$, where $\FF$ is a field of
characteristic $2$ (to avoid sign issues).  In our applications, $\FF$
is $\FF_2 = \ZZ/2\ZZ$.

\begin{definition}
  A \emph{quadratic algebra}~$A$ over $\Ground$ is a graded algebra
  generated by a
  finitely generated $\Ground$-module $V$ of elements of degree~$1$, with relations
  $R \subset V \otimes V$ living in degree~$2$.  Its \emph{quadratic
    dual algebra} $A^!$ is the algebra generated by $V^*$ with
  relations $R^\perp \subset (V \otimes V)^* \cong V^* \otimes V^*$.
\end{definition}

Here the last isomorphism switches the tensor factors, i.e., $(V
\otimes W)^* \cong W^* \otimes V^*$.

\begin{definition}\label{def:rank-1}
  For augmented \dg algebras $A,B$ over~$\Ground$, an $A$-$B$ bimodule
  $M$ of any type (\DD, \AAm, \DA, or \AD) is \emph{rank one} if it is
  isomorphic to $\Ground$ as a $\Ground$-$\Ground$ bimodule.
\end{definition}

For any quadratic algebra~$A$, there is a rank 1 type \DD\ bimodule
$\lsup{A}K(A)^{A^!}$, defined by
\begin{equation}
\label{eq:dual-bimod}
\delta^1({\bOne}) =  \sum_{i} v_i \otimes {\bOne} \otimes (v_i^*),
\end{equation}
where ${\bOne}$ is the generator of $K(A)$,
$\{v_i\}$ is a basis for $V$\!, and $\{v_i^*\}$ is the
dual basis for $V^*$.
The fact that $A^!$ has relations given by
$R^\perp$ is exactly what is necessary to make this a \DD\ bimodule.

\begin{definition}\label{def:quadratic-koszul}
  A quadratic algebra is \emph{Koszul} if $\lsup{A}K(A)^{A^!}$ is
  quasi-invertible. In this case we say that $A^!$ is the \emph{Koszul
    dual} of~$A$.
\end{definition}
(Bimodules $K$ and $L$ over $A$ and $B$ are \emph{quasi-inverses} if
$K\DT L\simeq \lsupv{A}[\Id]_A$ and $L\DT K\simeq \lsupv{B}[\Id]_B$;
if a quasi-inverse to $K$ exists, we say $K$ is \emph{quasi-invertible}.)

We compare the Definition~\ref{def:quadratic-koszul} with the standard
definition of Koszul duality for quadratic 
algebras in Proposition~\ref{prop:StandardDefinition} below.

Equation~\eqref{eq:dual-bimod} is reminiscent of
Equation~\eqref{eq:DDid}, suggesting that we think of $\Alg(\PMC)$ as
something like a ``quadratic'' algebra where the ``linear'' elements
are $a(\xi)$ for $\xi\in \Chord(\PMC)$.  
Since $\Alg(\PMC)$ has both linear terms in the
relations (like $a(\rho_1)a(\rho_2) = a(\rho_1 \uplus \rho_2)$ when
$\rho_1$ and $\rho_2$ abut) and a differential (like
$\partial a(\rho) = a(\rho_2) a(\rho_1) + \cdots$), the notion of
quadratic dual has to be extended.  This can be done; in particular,
the dual of a linear-quadratic algebra (where the relations have
linear and quadratic terms) is a quadratic-differential algebra (with
quadratic relations, but a differential); see, e.g.,
\cite[Chapter~5]{PP05:QuadraticAlgebras}.  The case of
linear-quadratic-differential algebras, with 
both types of terms, does not seem to be standard, but is a fairly
straightforward generalization.  Unfortunately, with the easiest
definitions $\Alg(\PMC,i)$ is not quite quadratic dual to
$\Alg(\PMC,-i)$, but instead to an algebra that is only homotopy
equivalent to $\Alg(\PMC,-i)$.  Instead of pursuing this discussion we
make the following definitions.

\begin{definition}\label{def:koszul-dual}
  Let $A$ and $B$ be augmented \dg algebras over $\Ground$. A
  \emph{Koszul dualizing bimodule} between $A$ and $B$ is a type \DD\
  structure $\lsup{A}K^B$ over $A$ and $B$ which is:
  \begin{itemize}
  \item quasi-invertible,
  \item rank $1$, and
  \item such that the image of $\delta^1$ lies in $A_+\otimes K
    \otimes B_+$.
  \end{itemize}
  Two \dg algebras $A$ and $B$ over $\Ground$ are \emph{Koszul dual}
  if there is Koszul dualizing bimodule between them.
\end{definition}
Observe that if it is quasi-invertible then the type \DD\ bimodule
$\lsup{A}K^{A^!}$ associated to a quadratic algebra satisfies the
conditions of Definition~\ref{def:koszul-dual}.
Note also that Definition~\ref{def:koszul-dual} is symmetric between
$A$ and $B$, by replacing $\lsup{A}K^B$ by its dual
$\lsup{B}\overline{K}^A$.

Recall that a rank-one type \DA\ bimodule
$\lsup{B}M_A$ with $\delta^1_1=0$
is the bimodule $\lsupv{B}[f]_A$ associated
to an $\Ainf$-map $f \co A \to B$, as
in~\cite[Definition~\ref*{LOT2:def:rank-1-DA-mods}]{LOT2}.
The choice of $f$ is determined by a choice of generator for $M$ as a
$\Ground$-module, and changes by conjugation by a unit in $\Ground$ if we change the
generator.  (If $\Ground = \bigoplus \ZZ/2\ZZ$, as is the case in bordered Floer
theory, then $1$ is the unique unit in~$\Ground$.)

\begin{lemma}\label{lem:bimod-quasi-isom}
  If $f, g \co A \to B$ are $\Ainf$-homomorphisms so that
  $\lsupv{B}[f]_A$ and $\lsupv{B}[g]_A$ are homotopy equivalent type
  \DA\ structures, then $f$ and $g$ induce conjugate maps on homology,
  i.e., there is a unit $[u]\in H_*(B)$ so that for all $[a]\in H_*(A)$,
  \[
  [f(a)]=[u][g(a)][u]^{-1}.
  \]
\end{lemma}
\begin{proof}
  By hypothesis there are maps
  $\phi\co\lsupv{B}[f]_A\to\lsupv{B}[g]_A$ and
  $\psi\co\lsupv{B}[g]_A\to\lsupv{B}[f]_A$, as well as homotopies $F\co
  \lsupv{B}[f]_A\to \lsupv{B}[f]_A$ from $\psi\circ \phi$ to
  $\Id_{[f]}$ and $G\co \lsupv{B}[g]_A\to \lsupv{B}[g]_A$ from
  $\phi\circ \psi$ to $\Id_{[g]}$. The $\Ainf$-relations give:
  \begin{align*}
    \bdy f_1(a)&=f_1 (\bdy a) &
    \bdy g_1(a)&=g_1( \bdy a)\\
    \bdy \phi_1 &=0&
    \bdy \psi_1&=0\\
    f_1(a)\phi_1+\phi_1g_1(a)&=\bdy \phi_2(a)+\phi_2(\bdy a) &
    g_1(a)\psi_1+\psi_1f_1(a)&=\bdy \psi_2(a)+\psi_2(\bdy a)\\
    \bdy F_1&=\phi_1\psi_1+1 &
    \bdy G_1&=\psi_1\phi_1+1.
  \end{align*}
  (We are abusing notation: the map $\phi_1\co
  \Ground\to B\otimes\Ground$ corresponds to an element $\phi_1\in B$,
  and similarly for $\psi_1$, $\phi_2$, and so on.) The
  equations on the last line imply that $[\psi_1]$ and $[\phi_1]$ are inverses in
  $H_*(B)$.  So, either equation on the third line implies that, if
  we take $a$ to be a cycle, then there is the following equation in
  homology:
  \[
  [f_1(a)][\phi_1]=[\phi_1][g_1(a)].
  \]
  That is, the maps on homology induced by $f_1$ and $g_1$ differ  by
  conjugation by $[\phi_1]$.
\end{proof}

\begin{lemma}
  \label{lem:QuasiInverse}
  Let $\lsupv{A}K^{B}$ be a quasi-invertible type $DD$ bimodule.
  Then its quasi-inverse is given by
  $\lsub{B}\overline{B}_B\DT \lsupv{B}{\overline K}{}^{A}\DT \lsub{A}A_A$.
\end{lemma}

(Here $\lsupv{B}{\overline K}{}^{A}$ denotes the 
  dual of $\lsupv{A}K^{B}$ in the sense of Definition~\ref{def:dual-type-D}.)

\begin{proof}
  This is essentially~\cite[Proposition~\ref*{LOT2:prop:DDAA-duality}]{LOT2}.
  We repeat the proof here for the reader's convenience.  Let $\lsub{B}L_A$
  be the quasi-inverse of $\lsupv{A}K^B$. Then
  \begin{align*}
    \lsub{B}L_A &\simeq 
    \Mor_{B}(\lsub{B}B_B,\lsub{B}L_A) \\
    &\simeq \Mor^A(\lsup{A}K^B\DT \lsub{B}B_B,
    \lsup{A}K^B\DT \lsub{B}L_A) \\
    &\simeq \Mor^A(\lsup{A}K^B\DT \lsub{B}B_B,
    \lsup{A}\Id_A) \\
    &\simeq {\overline{\lsupv{A}K^B
    \DT\lsub{B}B_B}} \DT \lsub{A}A_A \DT \lsup{A}\Id_A   \\
    &\simeq \lsub{B}{\overline B}_B 
    \DT \lsupv{B}{\overline{K}}{}^A\DT \lsub{A}A_A,
  \end{align*}
  where here the first step follows from the fact that any module is
  quasi-isomorphic to its cobar resolution, the second from the fact that $\lsup{A}K^B$
  is quasi-invertible (hence $\lsup{A}K^B\DT \cdot$ induces an equivalence of homotopy
  categories), the third from the fact that $K$ and $L$ are quasi-inverses,
  the fourth follows from Proposition~\ref{prop:D-Mor-is-DT}, and the fifth is 
  straightforward.
\end{proof}

\begin{lemma}
  \label{lem:RankOneInverse}
  Let $\lsup{A}K^B$ be a Koszul dualizing bimodule between $A$ and
  $B$. Then $\lsup{A}K^B$ has a rank-one quasi-inverse. Moreover
  the left $A$-module $\lsub{A}A_A\DT \lsupv{A}K^B\DT \lsub{B}{\overline B}$
  is a projective resolution of~$\Ground$, thought of as a left
  $A$-module via the augmentation~$\epsilon$.
\end{lemma}

\begin{proof}
  Lemma~\ref{lem:QuasiInverse} gives
  $\lsup{A}K^B\DT 
  \lsub{B}{\overline B}_B\DT \lsupv{B}{\overline K}{}^A\DT \lsub{A}A_A
  \simeq \lsup{A}\Id_A,$
  which in turn ensures that
  $\lsup{A}K^B\DT 
  \lsub{B}{\overline B}_B\DT \lsupv{B}{\overline K}{}^A\simeq
  \lsup{A}\Barop(A)^A$.
  From this it follows that
  $\lsub{A}A_A\DT \lsup{A}K^B\DT 
  \lsub{B}{\overline B}_B  \DT \lsupv{B}{\overline K}{}^A
  \DT \lsub{A}\Ground$
  is a projective resolution of $\lsub{A}\Ground$.
  Because the image of $\delta^1_{\overline K}$ lies in
  $B_+\otimes\Ground\otimes A_+$, and $A_+$ acts by zero on
  $\lsub{A}\Ground$,
  the terms in the above differential coming from
  $\lsup{B}{\overline K}{}^A\DT\lsub{A}\Ground$ are trivial;
  i.e.,
  $\lsub{A}A_A\DT \lsup{A}K^B\DT 
  \lsub{B}{\overline B}$
  is a projective resolution of $\lsub{A}\Ground$.
  In particular, its homology has rank one.

  Similarly, another application of Lemma~\ref{lem:QuasiInverse} gives 
  $\lsub{B}{\overline B}_B\DT \lsupv{B}{\overline K}{}^A\DT \lsub{A}A_A
  \DT \lsup{A}K^B\DT \lsub{B}\Ground
  \simeq \Ground$.
  Because the image of $\delta^1_{K}$ lies in
  $A_+\otimes\Ground\otimes B_+$, and $B_+$ acts by zero on
  $\lsub{B}\Ground$, it follows that
  $ \lsub{B}{\overline B}_B\DT\lsup{B}\overline{K}^A\DT\lsub{A}A$ is
  quasi-isomorphic to
  $\lsub{B}\Ground$. Thus, the homology of
  $ \lsub{B}{\overline B}_B\DT \lsup{B}\overline{K}^A\DT\lsub{A}A_A$
  is one-dimensional, and is the desired rank one 
  quasi-inverse to $\lsupv{A}K^B$.
\end{proof}

We will also use the following version of the homological perturbation
lemma:
\begin{lemma}\label{lem:homological-perturbation}
  Let $\Alg$ and $\Blg$ be $\Ainf$-algebras, $\lsub{\Alg}M_\Blg$ an $\Ainf$-bimodule
  over $\Alg$ and $\Blg$ and $\lsub{\Alg}N$ a left $\Ainf$-module over~$\Alg$. 
  Assume that $\Alg$ and $\Blg$ are defined over ground rings
  $\Ground$ and $\Groundl$ respectively, which are either $\Field$ or finite direct
  sums of copies of $\Field$, and that $\Alg$ and $\Blg$ are equipped
  with augmentations $\epsilon_\Alg\co \Alg\to \Ground$, $\epsilon_\Blg\co\Blg\to\Groundl$.
  Let $f\co \lsub{\Alg}N\to \lsub{\Alg}M$ be a
  quasi-isomorphism of left $\Ainf$-modules. Then there is an
  $\Ainf$-bimodule structure $\lsub{\Alg}N_\Blg$ on $N$, extending the
  given left $\Ainf$-module structure, so that $f$ can be extended to
  an $\Ainf$-bimodule quasi-isomorphism $F\co \lsub{\Alg}N_\Blg\to
  \lsub{\Alg}M_\Blg$.
\end{lemma}
\begin{proof}
  The proof is a simple extension of standard techniques (see, for
  instance,~\cite[Section 3.3]{AinftyAlg}, and the references
  therein). In fact, we will see that the result essentially follows
  from the corresponding result for modules, as formulated
  in~\cite[Lemma~\ref*{LOT4:lem:HomologicalPerturbation}]{LOT4}, say.

  Since is $\Alg$ defined over $\Ground$ which is a direct sum of copies
  of $\Field$, any $\Ainf$-quasi-isomorphism of $\Ainf$ $\Alg$
  modules is an $\Ainf$-homotopy equivalence
  (see~\cite[Proposition~\ref*{LOT2:prop:derived-is-derived-is-derived-is-derived}]{LOT2},
  say). So, let $g\co \lsub{\Alg}M\to\lsub{\Alg}N$ be a homotopy
  inverse to $f$ and let $T\co \lsub{\Alg}M\to\lsub{\Alg}M$ be a
  homotopy between $f\circ g$ and $\Id_M$. 
  
  An $\Ainf$-bimodule structure on $N$ is a map $m^N\co
  \Tensor^*A_+\otimes N\otimes\Tensor^*B_+\to N$ satisfying a
  compatibility condition.  Similarly, the maps $F$ are given by a map
  $F\co \Tensor^*(A_+)\otimes N\otimes\Tensor^*(B_+)\to M$ satisfying
  a compatibility condition. The maps $m^N$ and $F$ are
  defined by the top and bottom of Figure~\ref{fig:induced-strs}, respectively.
  \begin{figure}
  \begin{gather*}
  \mathcenter{
    {\begin{tikzpicture}[x=1cm,y=32pt,baseline=(x.base)]
        \node at (-1,7) (x) {$m^N(\overline{a}\otimes \x\otimes{\overline{b}})=$};
      \end{tikzpicture}}
    {\begin{tikzpicture}[x=1cm,y=32pt,baseline=(x.base)]
        \node at (-1,1) (terminal) {};
        \node at (-1,2) (g) {$g$};
        \node at (-1,3) (m1) {$m$};
        \node at (-1,4) (f) {$f$};
        \node at (0,5) (mu) {$\Delta$};
        \node at (-2,5) (lmu) {$\Delta$};
        \node at (-1,6) (x) {${\mathbf x}$};
        \node at (0,6) (alg) {${\overline b}$};
        \node at (-2,6) (lalg) {${\overline a}$};
        \draw[modarrow] (g) to (terminal);
        \draw[othmodarrow] (m1) to (g);
        \draw[modarrow] (x) to (f);
        \draw[othmodarrow] (f) to (m1);
        \draw[tensorblgarrow] (alg) to (mu);
        \draw[tensorblgarrow, bend left=15] (mu) to (m1);
        \draw[tensoralgarrow] (lalg) to (lmu);
        \draw[tensoralgarrow, bend right=15] (lmu) to (m1);
        \draw[tensoralgarrow, bend right=15] (lmu) to (g);
        \draw[tensoralgarrow, bend right=15] (lmu) to (f); 
      \end{tikzpicture}}
    {\begin{tikzpicture}[x=1cm,y=32pt,baseline=(x.base)]
        \node at (-1,0) (terminal) {};
        \node at (-1,1) (g) {$g$};
        \node at (-1,2) (m2) {$m$};
        \node at (-1,3) (T1) {$T$};
        \node at (-1,4) (m1) {$m$};
        \node at (-1,5) (f) {$f$};
        \node at (-2,6) (lmu) {$\Delta$};
        \node at (0,6) (mu) {$\Delta$};
        \node at (-3,7) (plus) {$+$};
        \node at (-1,7) (x) {${\mathbf x}$};
        \node at (0,7) (alg) {${\overline b}$};
        \node at (-2,7) (lalg) {${\overline a}$};
        \draw[othmodarrow] (m2) to (g);
        \draw[modarrow] (g) to (terminal);
        \draw[othmodarrow] (T1) to (m2);
        \draw[othmodarrow] (m1) to (T1);
        \draw[modarrow] (x) to (f);
        \draw[othmodarrow] (f) to (m1);
        \draw[tensorblgarrow] (alg) to (mu);
        \draw[tensorblgarrow, bend left=15] (mu) to (m2);
        \draw[tensorblgarrow, bend left=15] (mu) to (m1);
        \draw[tensoralgarrow] (lalg) to (lmu);
        \draw[tensoralgarrow, bend right=15] (lmu) to (m2);
        \draw[tensoralgarrow, bend right=15] (lmu) to (m1);
        \draw[tensoralgarrow, bend right=15] (lmu) to (g); 
        \draw[tensoralgarrow, bend right=15] (lmu) to (T1);
        \draw[tensoralgarrow, bend right=15] (lmu) to (f); 
      \end{tikzpicture}}
    {\begin{tikzpicture}[x=1cm,y=32pt,baseline=(x.base)]
        \node at (-1,-2) (terminal) {};
        \node at (-1,-1) (g) {$g$};
        \node at (-1,0) (m3) {$m$};
        \node at (-1,1) (T2) {$T$};
        \node at (-1,2) (m2) {$m$};
        \node at (-1,3) (T1) {$T$};
        \node at (-1,4) (m1) {$m$};
        \node at (-1,5) (f) {$f$};
        \node at (-2,6) (lmu) {$\Delta$};
        \node at (0,6) (mu) {$\Delta$};
        \node at (-1,7) (x) {${\mathbf x}$};
        \node at (-3,7) (plus) {$+$};
        \node at (1,7) (plus2) {$+$};
        \node at (2,7) (dots) {$\dots$};
        \node at (-2,7) (lalg) {${\overline a}$};
        \node at (0,7) (alg) {${\overline b}$};
        \draw[othmodarrow] (m3) to (g);
        \draw[modarrow] (g) to (terminal);
        \draw[othmodarrow] (m3) to (g);
        \draw[othmodarrow] (T2) to (m3);
        \draw[othmodarrow] (m2) to (T2);
        \draw[othmodarrow] (T1) to (m2);
        \draw[othmodarrow] (m1) to (T1);
        \draw[modarrow] (x) to (f);
        \draw[othmodarrow] (f) to (m1);
        \draw[tensorblgarrow] (alg) to (mu);
        \draw[tensorblgarrow, bend left=15] (mu) to (m3); 
        \draw[tensorblgarrow, bend left=15] (mu) to (m2); 
        \draw[tensorblgarrow, bend left=15] (mu) to (m1); 
        \draw[tensoralgarrow] (lalg) to (lmu);
        \draw[tensoralgarrow, bend right=15] (lmu) to (m3); 
        \draw[tensoralgarrow, bend right=15] (lmu) to (m2); 
        \draw[tensoralgarrow, bend right=15] (lmu) to (m1); 
        \draw[tensoralgarrow, bend right=15] (lmu) to (g); 
        \draw[tensoralgarrow, bend right=15] (lmu) to (T1); 
        \draw[tensoralgarrow, bend right=15] (lmu) to (T2); 
        \draw[tensoralgarrow, bend right=15] (lmu) to (f); 
      \end{tikzpicture}}}\\
  \mathcenter{
    {\begin{tikzpicture}[x=1cm,y=32pt,baseline=(x.base)]
        \node at (-1,3) (x) {$F(\overline{a}\otimes\x\otimes{\overline{b}})=$};
        \end{tikzpicture}}
    {\begin{tikzpicture}[x=1cm,y=32pt,baseline=(x.base)]
        \node at (-1,2) (terminal) {};
        \node at (-1,3) (f) {$f$};
        \node at (-1,4) (x) {${\mathbf x}$};
        \draw[modarrow] (x) to (f);
        \draw[algarrow] (f) to (terminal);
      \end{tikzpicture}}
    {\begin{tikzpicture}[x=1cm,y=32pt,baseline=(x.base)]
        \node at (-1,1) (terminal) {};
        \node at (-1,2) (T1) {$T$};
        \node at (-1,3) (m1) {$m$};
        \node at (-1,4) (f) {$f$};
        \node at (0,5) (mu) {$\Delta$};
        \node at (-2,5) (lmu) {$\Delta$};
        \node at (-1,6) (x) {${\mathbf x}$};
        \node at (0,6) (blg) {${\overline b}$};
        \node at (-2,6) (alg) {${\overline a}$};
        \node at (-3,6) (plus) {$+$};
        \draw[algarrow] (T1) to (terminal);
        \draw[algarrow] (m1) to (T1);
        \draw[modarrow] (x) to (f);
        \draw[algarrow] (f) to (m1);
        \draw[tensorblgarrow] (blg) to (mu);
        \draw[tensorblgarrow, bend left=15] (mu) to (m1); 
        \draw[tensoralgarrow] (alg) to (lmu);
        \draw[tensoralgarrow, bend right=15] (lmu) to (m1); 
        \draw[tensoralgarrow, bend right=15] (lmu) to (f); 
        \draw[tensoralgarrow, bend right=15] (lmu) to (T1);
      \end{tikzpicture}}
    {\begin{tikzpicture}[x=1cm,y=32pt,baseline=(x.base)]
        \node at (-1,0) (terminal) {};
        \node at (-1,1) (T2) {$T$};
        \node at (-1,2) (m2) {$m$};
        \node at (-1,3) (T1) {$T$};
        \node at (-1,4) (m1) {$m$};
        \node at (-1,5) (f) {$f$};
        \node at (0,6) (mu) {$\Delta$};
        \node at (-2,6) (lmu) {$\Delta$};
        \node at (-3,7) (plus) {$+$};
        \node at (-1,7) (x) {${\mathbf x}$};
        \node at (1,7) (plus2) {$+$};
        \node at (2,7) (dots) {$\dots$};
        \node at (0,7) (blg) {${\overline b}$};
        \node at (-2,7) (alg) {${\overline a}$};
        \draw[algarrow] (m2) to (T2);
        \draw[algarrow] (T2) to (terminal);
        \draw[algarrow] (T1) to (m2);
        \draw[algarrow] (m1) to (T1);
        \draw[modarrow] (x) to (f);
        \draw[algarrow] (f) to (m1);
        \draw[tensorblgarrow] (blg) to (mu);
        \draw[tensorblgarrow, bend left=15] (mu) to (m2);
        \draw[tensorblgarrow, bend left=15] (mu) to (m1);
        \draw[tensoralgarrow] (alg) to (lmu);
        \draw[tensoralgarrow, bend right=15] (lmu) to (m1);
        \draw[tensoralgarrow, bend right=15] (lmu) to (m2);
        \draw[tensoralgarrow, bend right=15] (lmu) to (f); 
        \draw[tensoralgarrow, bend right=15] (lmu) to (T1);
        \draw[tensoralgarrow, bend right=15] (lmu) to (T2);
      \end{tikzpicture}}
  }
  \end{gather*}
  \caption{\textbf{Induced $\Ainf$-bimodule structure and
      quasi-isomorphism.} The maps $\Delta\co \Tensor^*(A_+)\to
    \Tensor_*(A_+)^{\otimes n}$ and $\Tensor^*(B_+)\to
    \Tensor^*(B_+)^{\otimes n}$ are iterates of the obvious
    comultiplication on $\Tensor^*$.}\label{fig:induced-strs}
  \end{figure}
  
  It remains to check that the operations $m^N$ satisfy the
  $\Ainf$-bimodule relations and that the maps $F$ satisfy the
  $\Ainf$-bimodule homomorphism relations. Rather than doing this
  directly, we will use the bar construction to reduce to the case of
  modules verified in~\cite{LOT4}. Recall that the (reduced) \emph{bar
  resolution} of an $\Ainf$-module $\lsub{\Alg}M$ is given by
  $\Tensor^*A_+\otimes M$ with differential 
  \begin{multline}\label{eq:bar-res}
  \bdy(a_1\otimes\dots\otimes a_n\otimes
  \x)=\sum a_1\otimes\dots\otimes
  \mu_k(a_i,\dots,a_{i+k-1})\otimes\dots\otimes a_n\otimes\x\\
  +\sum a_1\otimes\dots\otimes m_{k}(a_{n-k+2},\dots,a_n,\x).
  \end{multline}
  The $\Ainf$-relation for an $\Ainf$-module $\lsub{\Alg}M$ is the
  same as the relation $\bdy^2=0$ on the bar resolution
  $\Tensor^*A_+\otimes M$.

  Similarly, the left bar resolution of an $\Ainf$ bimodule
  $\lsub{\Alg}M_\Blg$ is the right $\Ainf$-module given by
  $\Tensor_*A_+\otimes M$ with $m_1$ given by the formula in
  Equation~\eqref{eq:bar-res} and higher $\Ainf$-operations given by
  \[
  m_{1+\ell}((a_1\otimes\dots\otimes a_n\otimes
  \x),b_1,\dots,b_\ell)=\sum_{k=0}^\infty
  a_1\otimes\dots\otimes m_{k,1,\ell}(a_{n-k+1},\dots,a_n,\x,b_1,\dots,b_\ell).
  \]
  The $\Ainf$-bimodule relations for $\lsub{\Alg}M_\Blg$ are
  equivalent to the $\Ainf$-module relations for its left bar
  resolution $\Tensor_*A_+\otimes M$.

  Turning to the case at hand, inspecting
  the top of Figure~\ref{fig:induced-strs}
  and~\cite[Formula~(\ref*{LOT4:eq:induced-mod-str})]{LOT4}, the
  following two constructions give the same right $\Ainf$-operations:
  \begin{itemize}
  \item Constructing $\lsub{\Alg}N_\Blg$ by
    the top of Figure~\ref{fig:induced-strs} and then taking the left bar resolution.
  \item Taking the left bar resolution of $\lsub{\Alg}M_\Blg$ and the
    bar resolution of $\lsub{\Alg}N$ and then
    applying~\cite[Formula~(\ref*{LOT4:eq:induced-mod-str})]{LOT4} to
    construct a right module structure on the result.
  \end{itemize}
  So, it follows
  from~\cite[Lemma~\ref*{LOT4:lem:HomologicalPerturbation}]{LOT4} that
  the left bar resolution of $\lsub{\Alg}N_\Blg$ satisfies the
  $\Ainf$-module relations. This proves that $\lsub{\Alg}N_\Blg$
  satisfies the $\Ainf$-bimodule relations.

  A similar argument,
  comparing the bottom of Figure~\ref{fig:induced-strs} to
  \cite[Formula~(\ref*{LOT4:eq:induced-quasi-iso})]{LOT4}, shows
  that the map $F$ is an $\Ainf$-bimodule homomorphism.
\end{proof}

\begin{proposition}\label{prop:koszul-dual-isom}
  If $A$ is any augmented \dg algebra and $B$ and $C$ are both Koszul
  dual to $A$, then $B$ and $C$ are quasi-isomorphic.
\end{proposition}

\begin{proof}
  Let $\lsup{B}K^A$ and $\lsup{C}K'{}^A$ be the dualizing bimodules,
  and let $\lsub{A}L_C$ and $\lsub{A}L'{}_C$ be their respective
  quasi-inverses.  Since $\lsup{B}K^A$ is right bounded, we can form
  the \DA\ bimodule $\lsup{B}K^A \DT \lsub{A}L'{}_C$.  
  According to Lemma~\ref{lem:RankOneInverse}, we can use a model
  for $\lsub{A}L'{}_C$ which has rank one. Further, since $\lsub{A}L'
  \simeq \lsub{A}\Ground$, by
  Lemma~\ref{lem:homological-perturbation}, $\lsub{A}L'_C$ is
  isomorphic to a rank $1$ bimodule $\lsub{A}L''_C$ such that
  $m_{i+1}(a_1,\dots,a_i,l)= 0$ for any $a_1,\dots,a_i\in A_+$ and
  $l\in L''$.
  From this and the fact that the image of $\delta^1_K$ lies in
  $B_+\otimes \Ground\otimes A_+$, it follows that $\delta^1_1= 0$ on the rank one
  type \DA\ bimodule $\lsup{B}K^A \DT \lsub{A}L''_C$.

  Thus, $\lsup{B}K^A\DT\lsub{A}L''_C$
  is the bimodule
  $\lsupv{B}[f]_C$ associated to an $\Ainf$-map $f \co C \to B$.
  Similarly, $\lsup{C}K'{}^A \DT \lsub{A}L_B$ (for an appropriate
  model for $\lsub{A}L_B$) is the bimodule
  $\lsupv{C}[g]_B$ of an $\Ainf$-map $g \co B \to C$.  From the fact
  that $K$ and $L$ (respectively $K'$ and $L''$) are quasi-inverses, it
  follows that $\lsupv{B}[f]_C$ and $\lsupv{C}[g]_B$ are
  quasi-inverses to each other, which in turn implies (by
  Lemma~\ref{lem:bimod-quasi-isom}) that $f$ and $g$ induce
  isomorphisms on homology, as desired.
\end{proof}

Recall that for an augmented $\Ainf$-algebra $\Alg$, the cobar
resolution $\Cobarop(\Alg)$ is $T^*(\Alg_+[1]^*)$, the tensor algebra
on the (shifted) dual to $\Alg_+$.  $\Cobarop(\Alg)$ is itself a \dg
algebra, with a product that is the tensor product in the tensor
algebra and a differential that is dual to $\overline{D}{}^{\Alg} :
T^*(\Alg[1]) \to T^*(\Alg[1])$, the operation that encodes the
multiplications $\mu_i$ on~$\Alg$ \cite[Section~\ref*{LOT2:sec:defin-Ainfty-alg}]{LOT2}.

\begin{proposition}\label{prop:cobar-koszul}
  Any augmented \dg algebra $A$ is Koszul dual to $\Cobarop(A)$.
\end{proposition}

\begin{proof}
  The dualizing bimodule $\lsup{A}K^{\Cobarop(A)}$ is defined by
  \begin{align*}
    \delta^1(\bOne) &= \sum_i a_i \otimes \bOne \otimes a_i^*[1],
  \end{align*}
  where $a_i$ runs over a basis of $A_+$ and $a_i^*$ is the dual basis
  of $A^*$.  Its quasi-inverse is the bimodule $\lsub{\Cobarop(A)}L_A$ defined by
  \[
  m_{1,1,n}(\langle b_1^* | b_2^* | \dotsb | b_n^*\rangle,
    \bOne, a_1, \dotsc, a_{n-1}, a_n) =
    b_1^*(a_n) \, b_2^*(a_{n-1}) \dotsm b_n^*(a_1)\,\bOne
  \]
  (with all other products~$0$), where $b_i^*(a_{n-i})$ is the
  canonical pairing between $b_i^* \in A_+[1]^*$ and $a_{n-i} \in A$.
\end{proof}

A Koszul dual to $A$ is, by these definitions, just a \dg algebra
that is quasi-isomorphic to
$\Cobarop(A)$; compare ~\cite{LPWZ08:KoszulEquivAinf}.  In particular,
for any augmented \dg algebra,
$\Cobarop(\Cobarop(A))$ is quasi-isomorphic to~$A$.
With luck (as in the classical case of
quadratic Koszul duality of Definition~\ref{def:quadratic-koszul}),
there is a Koszul dual that is much smaller than $\Cobarop(A)$.

We now compare our definition of Koszul duality for quadratic algebras with
the more familiar one, see for example~\cite[Definition~2.1.1]{PP05:QuadraticAlgebras}.
\begin{proposition}
  \label{prop:StandardDefinition}
  If $A$ is a Koszul, quadratic algebra in the sense of
  Definition~\ref{def:quadratic-koszul}, then
  the module $\lsub{A}A_A\DT \lsup{A}K(A)^{A^!} \DT
  \lsub{A^!}{\overline{A^!}}$ is a graded projective
  resolution of $\Ground$ whose generators in homological
  degree~$i$ are also in algebraic degree~$i$.
\end{proposition}

\begin{proof}
  Lemma~\ref{lem:RankOneInverse} guarantees that
  $\lsub{A}A_A\DT \lsup{A}K(A)^{A^!} \DT \lsub{A^!}{\overline{A^!}}$ 
  is a projective resolution of~$\Ground$, so we only need to check
  the grading property.  The quadratic algebra~$A$ is automatically
  graded, with the generators~$V$ living in degree~$1$. In order
  to extend this to a bigrading on 
  $\lsub{A}A_A\DT \lsup{A}K(A)^{A^!} \DT \lsub{A^!}{\overline{A^!}}$ 
  it is natural to think of $A$ as bigraded, so that
  $V$ lives in
  grading $(1,0)$. (Note that the second component of this bigrading
  is trivial on $A$.) 
  In order for Equation~\eqref{eq:dual-bimod} to give a differential
  on  $\lsup{A}K(A)^{A^!}$ which changes the bigrading by $(0,-1)$,
  $V^* \subset A^!$ must lie in grading $(-1,-1)$.  It follows that $A^!$
  lies in gradings $(-i,-i)$ for $i \ge 0$ and $\overline{A^!}$ lies
  in gradings $(i,i)$, again with $i \ge 0$.  This implies the proposition.
\end{proof}

Similarly, since $\Cobarop(A) \simeq A^!$ and $A^!$ lies in gradings
$(-i,-i)$, it follows that $\Ext^{ij}(\Ground, \Ground)=0$ if $i \ne j$.

\begin{remark}
  Definition~\ref{def:koszul-dual} and
  Propositions~\ref{prop:koszul-dual-isom} and~\ref{prop:cobar-koszul}
  can be extended to the case of $\Ainf$-algebras.  The only
  difficulty is defining the notion of a \DD\ bimodule over two
  $\Ainf$-algebras; see \cite[Remark \ref*{LOT2:rem:DD-difficult}]{LOT2}.
\end{remark}

\begin{remark}
  Definition~\ref{def:koszul-dual} is quite similar to those considered by
  Lef\'evre\hyp Hasegawa~\cite{LefevreAInfinity} and
  Keller~\cite{KellerLefevre}, except that they work with an algebra
  and a coalgebra rather than two algebras.  More precisely, for $C$ a
  coalgebra that is finite-dimensional in each grading, let $C^*$ be
  the graded dual, which is an algebra.  Then a \emph{twisting
    cochain} $\tau\co C \to A$ (see, e.g., \cite[Section 2.3]{KellerLefevre}) is the
  same data as a rank~$1$ \DD\ bimodule $\lsup{A}K^{C^*}$, and
  Definition~\ref{def:koszul-dual} is close to the definition of
  a Koszul-Moore triple in loc.\ cit., with some difference in the
  technical conditions.
\end{remark}

\begin{remark}
  Tensoring with the dualizing bimodule $\lsup{A}K^B$ does not give an
  equivalence of categories between derived categories of modules over
  $A$ and over $B$, but rather between $\lsub{B}\ModCat$,
  $\Ainf$-modules over $B$, and $\lsubsupv{u}{A}{\ModCat}$, the homotopy
  category of \emph{unbounded} type~$D$ structures.  This is
  presumably related to the fact that the full derived categories of
  modules over Koszul dual algebras are not equivalent in general
  (see~\cite[Section 2.12]{BGS96:KoszulReps} and \cite[Section 1]{KellerLefevre}).
  In the case of the algebras
  considered in this paper, however, the identity bimodule
  $\lsupv{\Alg}[\Id]_\Alg$ is homotopy equivalent to a bounded module
  (coming from, for instance, an admissible diagram for the identity map),
  so the categories of bounded and unbounded type~$D$
  structures are quasi-equivalent,
  $\lsubsupv{u}{\Alg}{\ModCat} \simeq \lsubsupv{b}{\Alg}{\ModCat}
  \simeq \lsub{\Alg}{\ModCat}$. 
\end{remark}

\subsection{Koszul duality in bordered Floer homology}
\label{sec:koszul-bordered}
The formulation of Koszul duality in terms of bimodules is well-suited
to bordered Floer homology: it allows us to use the combinatorics of
Heegaard diagrams to prove the desired Koszul duality for our
algebras.

\begin{proposition}\label{prop:koszul-dual-1}
  The algebra $\Alg(\PMC, i)$ is Koszul dual to $\Alg(\PMC, -i)$.
\end{proposition}

\begin{proof}
  The bimodules $\lsup{\Alg(\PMC,i)}\CFDDa(\Id)^{\Alg(\PMC,-i)}$ and
  $\lsub{\Alg(\PMC,-i)}\CFAAa(\Id)_{\Alg(\PMC,i)}$, computed with
  respect to the standard Heegaard diagram for the identity map,
  obviously satisfy the conditions of
  Definition~\ref{def:koszul-dual}, so $\Alg(\PMC,i)$ and
  $\Alg(\PMC,-i)$ are Koszul dual.
\end{proof}

For the other duality, we need to consider another diagram.

\begin{construction}
  \label{construct:HalfIdentity}
  Given an $\alpha$-pointed matched circle $\PMC^\alpha$, the
  \emph{half-identity diagram} $\cG(\PMC^\alpha)$ is the
  $\alpha\Hyph\beta$-bordered Heegaard diagram obtained as
  follows. Let $\Sigma_\dr$ be the disk with one-handles attached to
  its boundary as specified by $-\PMC$. Let $\alphas^a$ denote curves
  running through the one-handles, meeting the boundary along the
  pointed matched circle $\PMC$. Let $\betas^a$ be a collection of
  dual arcs: there is one in each one-handle, and $\beta_i^a$ meets
  only $\alpha_i^a$, transversely, and in a single point.  Finally, attach another
  one-handle to $\Sigma_{\dr}$ so as to separate the
  $\alpha$-endpoints from the $\beta$-endpoints, to obtain a Heegaard
  surface~$\Sigma$.  The resulting $\alpha\Hyph\beta$-bordered is
  ${\mathcal G}(\PMC^\alpha)$.

  The diagram $\mathcal{G}(\PMC^\alpha)$ has two boundary components,
  one of which is $\PMC^\alpha$ and the other of which is called the {\em dual pointed
    matched circle} and denoted
  $\PMC^\beta_*$. (The pointed matched circle $\PMC^\beta_*$ naturally
  corresponds to turning
  the Morse function specifying $\PMC$ upside-down.)  Let
  $\PMC^\alpha_*$, or just $\PMC_*$, denote the $\alpha$-pointed
  matched circle twin to $\PMC^\beta_*$.

  Sometimes, we
  write $\cG(\PMC^\alpha,\PMC_*^\beta)$ to indicate both
  boundaries.  $\cG(\PMC_*^\beta,\PMC^\alpha)$ is the same diagram
  with the roles of $\partial_L$ and $\partial_R$
  switched.
\end{construction}

\begin{figure}
  \centering
  \input{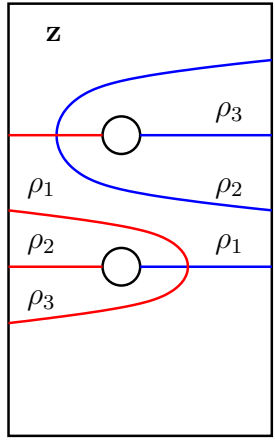}
  \caption{\textbf{Heegaard diagram for ${\mathcal G}$ for a genus one
      surface.} This is the Heegaard diagram $\mathcal{G}(\PMC^\alpha)$
    for the case where the genus is one, so 
    $\PMC=\PMC_*$, as described in
    Construction~\ref{construct:HalfIdentity}.}
  \label{fig:StrangeDuality}
\end{figure}

For a picture of the Heegaard diagram for ${\mathcal G}(\PMC^\alpha)$
for the torus,
see Figure~\ref{fig:StrangeDuality}.
The standard identity diagram is $-\cG(\PMC^\alpha, \PMC_*^\beta)
\sos{\bdy_R}{\cup}{\bdy_L} \cG(-\PMC_*^\beta,-\PMC^\alpha)$.

\begin{proof}[Proof of Theorem~\ref{thm:KoszulDual}]
  The first part is Proposition~\ref{prop:koszul-dual-1}.
  
  The bimodules
  $\lsup{\Alg(\PMC,i)}\CFDDa(\mathcal{G}(\PMC))^{\Alg(\PMC_*,i)}$ and 
  $\lsub{\Alg(\PMC_*,i)}\CFAAa(-\mathcal{G}(\PMC))_{\Alg(\PMC,i)}$
  also satisfy the conditions of
  Definition~\ref{def:koszul-dual}, so $\Alg(\PMC,i)$ and
  $\Alg(\PMC_*,i)$ are Koszul dual.
  
  Proposition~\ref{prop:koszul-dual-isom} now implies that
  $\Alg(\PMC,-i)$ is quasi-isomorphic to $\Alg(\PMC_*,i)$.
\end{proof}

It is interesting to note that in the case of the pointed matched
circle $\PMC$ for the torus algebra with $i=0$, both bimodules 
$\CFDDa(\mathcal{G}(\PMC),0)$ and $\CFDDa(\Id)$ give Koszul self-dualities
of the torus algebra $\Alg(\PMC,0)$. However, the two bimodules are different;
the bimodule
$\CFDDa(\mathcal{G}(\PMC),0)$ is given by
\[\delta({\bOne})=\rho_1\otimes {\bOne} \otimes \rho_1 +
\rho_2\otimes {\bOne} \otimes \rho_2 + \rho_3\otimes {\bOne} \otimes
\rho_3\]
(a fact which can be verifying by enumerating holomorphic
curves in Figure~\ref{fig:StrangeDuality}).  Contrast this with
$\CFDDa(\Id, 0)$, which by Equation~\eqref{eq:DD-identity} is given by
\[
\delta(\bOne) = \rho_1\otimes {\bOne} \otimes \rho_1 +
\rho_2\otimes {\bOne} \otimes \rho_2 + \rho_3\otimes {\bOne} \otimes
\rho_3 + \rho_{123}\otimes \bOne \otimes \rho_{123}.
\]
If we tensor one of these bimodules with the inverse of the other, we
get a non-trivial $\Ainf$-automorphism $f\co \Alg(T^2)\to \Alg(T^2)$,
given by
\begin{align*}
  f_1(x) &= x\\
  f_3(\rho_3, \rho_2, \rho_1) &= \rho_{123},
\end{align*}
with all other terms being~$0$.  (This automorphism can also be
computed by counting holomorphic curves in
Figure~\ref{fig:StrangeDuality} as a \DA\ bimodule.)

In a different direction, the symmetry $H_*(\Alg(\PMC,i))\cong
H_*(\Alg(\PMC_*,-i))$ explains some numerical coincidences apparent in
the homology calculations from~\cite{LOT2}.  

For $\PMC$ any pointed matched circle of genus $k>0$,
$\Alg(\PMC,-k)\cong \Field$ and
$\Alg(\PMC,-k+1)$ have no differential, so
\begin{align*}
  \dim(H_*(\Alg(\PMC,-k)))&=\dim(\Alg(\PMC,-k))=1 \\
  \dim(H_*(\Alg(\PMC,-k+1)))&=\dim(\Alg(\PMC,-k+1))=8k^2,
\end{align*}
both of which depend only on $k$, not the pointed matched circle.  It now follows from
Theorem~\ref{thm:KoszulDual} that 
\begin{align*}
  H_*(\Alg(\PMC,k-1))&=8k^2 \\
  H_*(\Alg(\PMC,k))&=1,
\end{align*}
for any genus $k$ pointed matched circle $\PMC$, as well.

In~\cite{LOT2}, we computed that if $\PMC$ is the split pointed matched
circle of genus $2$ then 
\[
\sum_{i}\dim(H_*(\Alg(\PMC,i)))=T^{-2}+32T^{-1}+98+32T+T^2,
\]
while if $\PMC$ denotes the antipodal pointed matched circle of genus
$2$ then
\[
\sum_{i}\dim(H_*(\Alg(\PMC,i)))=T^{-2}+32T^{-1}+70+32T+T^2.
\]
The coincidences in the Poincar{\'e} polynomials
in all but the middle-most term (and their symmetry) is now explained by Koszul duality.

In general,
$H_*(\Alg(\PMC),i)$ is not necessarily isomorphic to
$H_*(\Alg(\PMC),-i)$.  For instance, for the genus~$3$ pointed
matched circle~$\PMC_1$ with matched points 
\[
(1,7), (2,9), (3,5), (4,6), (8,11), (10,12),
\]
computer computation gives
\[
\sum_{i}\dim(H_*(\Alg(\PtdMatchCirc_1,i))) \cdot T^i=
T^{-3}+72\cdot T^{-2} + 600\cdot T^{-1} + 1224
  + 616 \cdot T + 72\cdot T^2 + T^3,
\]
and in particular $\dim(H_*(\Alg(\PMC_1,1)))\neq
\dim(H_*(\Alg(\PMC_1,-1)))$. Similarly, for the dual pointed matched
circle $\PMC_2 = \PMC_1^*$, with matched points
\[
(1,10), (2,4), (3,12), (5,11), (6,8), (7,9),
\]
computer computation gives
\[
\sum_{i}\dim(H_*(\Alg(\PtdMatchCirc_2,i))) \cdot T^i=
T^{-3}+72\cdot T^{-2} + 616\cdot T^{-1} + 1224
  + 600 \cdot T + 72\cdot T^2 + T^3,
\]
which is consistent with Theorem~\ref{thm:KoszulDual}.

\begin{remark}
In light of Auroux's reinterpretation of bordered Floer theory~\cite{AurouxBordered}, the referee
points out that it is interesting to compare the results of this section with~\cite[Section~(5k)]{SeidelBook}.
\end{remark}


\appendix
\section{User's guide to orientation
  conventions}\label{sec:conventions}
Type $D$ structures (Definition~\ref{def:type-D-str}) are written with
the algebra as a superscript, and modules with the algebra as a
subscript. Examples:
\begin{center}
  \begin{tabular}{ll}
    $\lsup{\Alg}M$ & Left type $D$ structure over $\Alg$.\\
    $M_\Alg$ & Right $\Ainf$-module over $\Alg$.\\
    $\lsup{\Alg}M_\Blg$ & Type \DA\ bimodule; left type $D$ over
    $\Alg$, right type $A$ over $\Blg$.
  \end{tabular}
\end{center}

The algebras $\Alg(\PMC)$ and $\Alg(-\PMC)$ are opposites:
\[
\Alg(\PMC)^\op=\Alg(-\PMC).
\]
So, there are identifications
$\ModCat_{\Alg(\PMC)}\equiv\lsub{\Alg(-\PMC)}\ModCat$ and
$\ModCat^{\Alg(\PMC)}\equiv\lsupv{\Alg(-\PMC)}\ModCat$. With respect to
these identifications,
\[
\lsup{\Alg(-\PMC)}\CFDa(\HD)\equiv \CFDa(\HD)^{\Alg(\PMC)}\qquad\qquad
\lsub{\Alg(\PMC)}\CFAa(\HD)\equiv \CFAa(\HD)_{\Alg(-\PMC)}.
\]
The following modules are associated to $\alpha$- or $\beta$-bordered
Heegaard diagrams:
\begin{center}
  \smallskip
  \begin{tabular}{lll}
    \toprule
    Diagram & Type $D$ structures & Type $A$ structures \\
    \midrule
    $\HD^\alpha$, $\bdy \HD^\alpha = \PMC^\alpha$
       & $\lsup{\Alg(-\PMC)}\CFDa(\HD)\equiv\CFDa(\HD)^{\Alg(\PMC)}$
           & $\lsub{\Alg(-\PMC)}\CFAa(\HD)\equiv\CFAa(\HD)_{\Alg(\PMC)}$\\
    $\HD^\beta$, $\bdy\HD^\beta = \PMC^\beta$
       & $\lsup{\Alg(\PMC)}\CFDa(\HD)\equiv\CFDa(\HD)^{\Alg(-\PMC)}$
           & $\lsub{\Alg(\PMC)}\CFAa(\HD)\equiv\CFAa(\HD)_{\Alg(-\PMC)}$ \\
    \bottomrule
  \end{tabular}
  \smallskip
\end{center}

The surface associated to an orientation-reversed pointed matched
circle is the orientation reverse of the surface associated to the
pointed matched circle: $F(-\PMC)=-F(\PMC)$.

If $\HD^\alpha$ is an $\alpha$-bordered Heegaard diagram there is a
corresponding $\beta$-bordered Heegaard diagram~$\overline{\HD}^\beta$
(Definition~\ref{def:dual-HD-1}). The corresponding invariants are
dual (Proposition~\ref{prop:beta-is-dual}):
\begin{align*}
  \overline{\lsup{\Alg(-\PMC)}\CFDa(\HD^\alpha)}&=\CFDa(\overline{\HD}^\beta)^{\Alg(-\PMC)} &
  \overline{\CFAa(\HD^\alpha)_{\Alg(\PMC)}}&=\lsub{\Alg(\PMC)}\CFAa(\overline{\HD}^\beta).
\end{align*}
One can also reverse the orientation of $\HD^\alpha$, giving a new
$\alpha$-bordered Heegaard diagram $-\HD^\alpha$. Again, the
invariants are dual (Proposition~\ref{prop:minus-hd-mod}):
\begin{align*}
  \overline{\lsup{\Alg(-\PMC)}\CFDa(\HD^\alpha)}&=\CFDa(-\HD^\alpha)^{\Alg(-\PMC)} &
  \overline{\CFAa(\HD^\alpha)_{\Alg(\PMC)}}&=\lsub{\Alg(\PMC)}\CFAa(-\HD^\alpha).
\end{align*}
Analogous statements hold for bimodules associated to bordered
Heegaard diagrams with two boundary components
(Propositions~\ref{prop:beta-bimod-is-dual}
and~\ref{prop:minus-hd-bimod}).

Given a strongly based mapping class $\phi\co F(\PMC_1)\to F(\PMC_2)$
there is an associated mapping cylinder $M_\phi$, with $\bdy_L M_\phi
= -F(\PMC_1)$ and $\bdy_R M_\phi = F(\PMC_2)$; see
Construction~\ref{constr:AssocHomeo} and also~\cite[Section
\ref*{LOT2:sec:DiagramsForAutomorphisms}]{LOT2}. To emphasize: maps go
from (minus) the left boundary to the right boundary. There is an
action of the mapping class groupoid on the set of bordered
$3$-manifolds, by $\phi(Y,\psi\co F(\PMC_1)\to \bdy
Y)=(Y,\psi\circ\phi^{-1}\co F(\PMC_2)\to \bdy Y)$ (where $\phi\co
F(\PMC_1)\to F(\PMC_2)$); see
Definition~\ref{def:mcg-action}. Equivalently,
$\phi(Y)=Y\cup_{F(\PMC_1)}M_\phi$.

Note that if we reverse the roles of the left and right boundary on
the mapping cylinder $M_\phi$, where $\phi\co -\bdy_L M_\phi \to
\bdy_R M_\phi$, we get the mapping cylinder of $-\phi^{-1}\co -\bdy_R
M_\phi \to \bdy_L M_\phi$.  In particular, switching the two sides of
the mapping cylinder of a positive Dehn twist gives the mapping
cylinder of another positive Dehn twist, as both taking the inverse
of~$\phi$ and reversing the orientation of the surfaces switch
positive and negative Dehn twists.

The strongly based mapping class $\tau_\bdy$ plays a special role. The
map $\tau_\bdy\co F(\PMC)\to F(\PMC)$ is a positive Dehn twist around
the boundary of the preferred disk $D_\alpha\cup s\cup D_\beta$ in
$F(\PMC)$. The effect of gluing $\tau_\bdy$ to a strongly bordered
$3$-manifold with two boundary components is to \emph{decrease} the
framing on the arc by~$1$.

Of particular importance in this paper is the bordered diagram
$\Denis(\PMC)$ and its mirror $\MirrorDenis(\PMC)$
(Section~\ref{sec:denis}). These are defined so that
\[
\bdy\Denis(\PMC)=\PMC^\alpha\amalg\PMC^\beta=\bdy\MirrorDenis(\PMC).
\]

To make it clear which boundary is being glued, we often write
$\lsup{\alpha}\Denis(\PMC)^\beta$ or
$\lsup{\beta}\Denis(\PMC)^\alpha$, to indicate whether we think of the
$\alpha$- of $\beta$-boundary of $\Denis(\PMC)$ as on the left. These
are two ways of writing the same diagram. In particular,
\[
\lsub{\Alg(\PMC)}\CFAAa(\lsup{\alpha}\Denis(\PMC)^\beta)_{\Alg(\PMC)}
=\lsub{\Alg(\PMC)}\CFAAa(\lsup{\beta}\Denis(\PMC)^\alpha)_{\Alg(\PMC)}=\lsub{\Alg(\PMC)}\Alg(\PMC)_{\Alg(\PMC)};
\]
and in both cases, the $\alpha$-boundary corresponds to the right
action. (The second equality uses
Proposition~\ref{prop:CFAA-of-Denis}.) If the $\alpha$ boundary
corresponds to the left action then the bimodules are:
\[
\lsub{\Alg(-\PMC)}\CFAAa(\lsup{\alpha}\Denis(\PMC)^\beta)_{\Alg(-\PMC)}
=\lsub{\Alg(-\PMC)}\CFAAa(\lsup{\beta}\Denis(\PMC)^\alpha)_{\Alg(-\PMC)}=\lsub{\Alg(-\PMC)}\Alg(-\PMC)_{\Alg(-\PMC)}.
\]
Similarly,
\[
\lsub{\Alg(\PMC)}\CFAAa(\lsup{\alpha}\MirrorDenis(\PMC)^\beta)_{\Alg(\PMC)}=\lsub{\Alg(\PMC)}\overline{\Alg(\PMC)}_{\Alg(\PMC)}.
\]

The following is a representative sample of the valid gluings of
$\Denis$ and $\MirrorDenis$ pieces, and what they represent, as maps
from (minus) the left boundary to the right boundary:
\begin{align*}
    \lsup{\alpha}\Denis(-\PMC)^\beta&\cup\lsup{\beta}\Denis(\PMC)^\alpha
    & \tau_\bdy&\co F(\PMC^\alpha)\to F(\PMC^\alpha)\\
    \lsup{\alpha}\Denis(\PMC)^\beta&\cup\lsup{\beta}\Denis(-\PMC)^\alpha
    & \tau_\bdy&\co F(-\PMC^\alpha)\to F(-\PMC^\alpha)\\
    \lsup{\beta}\Denis(-\PMC)^\alpha&\cup\lsup{\alpha}\Denis(\PMC)^\beta
    & \tau_\bdy&\co F(\PMC^\beta)\to F(\PMC^\beta)\\
    \lsup{\alpha}\MirrorDenis(-\PMC)^\beta&\cup\lsup{\beta}\Denis(\PMC)^\alpha
    & \Id&\co F(\PMC)\to F(\PMC)\\
    \lsup{\alpha}\Denis(-\PMC)^\beta&\cup\lsup{\beta}\MirrorDenis(\PMC)^\alpha
    & \Id&\co F(\PMC)\to F(\PMC)\\
    \lsup{\alpha}\MirrorDenis(-\PMC)^\beta&\cup\lsup{\beta}\MirrorDenis(\PMC)^\alpha
    & \tau_\bdy^{-1}&\co F(\PMC)\to F(\PMC).
\end{align*}
See Propositions~\ref{prop:Denis-reps} and~\ref{prop:mdenis-reps} and
Corollary~\ref{cor:denis-glue-reps}.


\bibliographystyle{hamsalpha}
\bibliography{heegaardfloer}
\end{document}